\documentclass[11pt]{amsart}

\usepackage{amssymb, epsfig, xcolor, upgreek, cancel, enumitem}

\usepackage{pgfplots}

\pgfplotsset{/pgf/number format/use comma, compat=newest}

\newcommand{\essinf}{\mathop{\rm ess\hspace{0.1 cm}inf}\limits}
\newcommand{\emm}{\mathpzc{M}}
\newcommand{\enn}{\mathpzc{N}}
\newcommand{\dimoH}{{\it Proof of Theorem \ref{harnack}}\ \ -\ \ }
\newcommand{\ess}{b}

\newcommand{\R}{{\bf R}}

\newcommand{\N}{{\bf N}}

\newfont{\bbten}{bbold12}

\newcommand{\miu}{\leqslant}
\newcommand{\mau}{\geqslant}

\newtheorem{teorema}{Theorem}[section]

\newtheorem{oss}[teorema]{\sc Remark}
\newtheorem{nota}[teorema]{\sc Nota}
\newtheorem{esempio}[teorema]{\sc Example}
\newtheorem{definition}[teorema]{Definition}
\newtheorem{prop}[teorema]{Proposition}
\newtheorem{lemma}[teorema]{Lemma}
\newtheorem{cor}[teorema]{Corollary}

\newcommand{\be}{\begin{equation}}
\newcommand{\ee}{\end{equation}}
\newcommand{\ba}{\begin{array}}
\newcommand{\ea}{\end{array}}
\newcommand{\bea}{\begin{eqnarray}}
\newcommand{\eea}{\end{eqnarray}}
\newcommand{\arst}{\renewcommand{\arraystretch}{1.5}}

\newcommand{\bd}{\begin{definition}}
\newcommand{\ed}{\end{definition}}
\newcommand{\bprop}{\begin{prop}}
\newcommand{\eprop}{\end{prop}}
\newcommand{\bthm}{\begin{teorema}}
\newcommand{\ethm}{\end{teorema}}
\newcommand{\boss}{ \begin{oss}  \rm -\ }
\newcommand{\eoss}{ \end{oss} }
\newcommand{\bex}{ \begin{esempio}  \rm -\ }
\newcommand{\eex}{ \end{esempio} }
\newcommand{\bcor}{\begin{cor}}
\newcommand{\ecor}{\end{cor}}
\newcommand{\besempio}{ \begin{esempio}  \rm -\ }
\newcommand{\eesempio}{ \end{esempio} }
\newcommand{\bnota}{ \begin{nota}  \rm -\ }
\newcommand{\enota}{ \end{nota} }

\newcommand{\dimo}{{\it Proof}\ \ -\ \ }

\newcommand{\finedimo}{\hfill $\square$ \\}
\newcommand{\fine}{ \\}

\newcommand{\om}{\omega}

\newcommand{\V}{\mathcal{V}}

\newcommand{\X}{\mathcal{X}}

\newcommand{\loc}{\rm loc}

\newcommand{\osc}{\mathop{\rm osc}\limits}

\makeglossary
\makeindex

\usepackage{blindtext}

\makeglossary
\makeindex

\usepackage{blindtext}
\title{A Harnack inequality for solutions of elliptic-parabolic equations}

\author{\sc Fabio Paronetto}
\thanks{F. Paronetto - Dipartimento di Matematica ``Tullio Levi Civita'', Universit\`a di Padova, via 
Trieste 63, 35121, Padova, Italy.  \\ e-mail: {\sf fabio.paronetto@unipd.it} - Fax +39 049 8271333, phone +39 049 8271318, ORCID 0000-0002-0017-7805}

\date{\today}

\DeclareMathAlphabet{\mathpzc}{OT1}{pzc}{m}{it}

\begin{document}


\maketitle

\tableofcontents

\begin{abstract}
We want to prove a Harnack type inequality for solutions of strongly degenerate parabolic, or elliptic-parabolic, equations.
To do that, we first define a De Giorgi class of order $p = 2$ that contains the solutions of
evolution equations of the types
$\uprho (x,t) u_t + A u = 0$ and $(\uprho (x,t) u)_t + A u = 0$,
where $\uprho > 0$ almost everywhere and $A$ is a suitable elliptic operator.
For functions belonging to this class we prove an inhomogeneous parabolic Harnack inequality, i.e. a Harnack inequality that
takes into account the mean value of $\uprho$ in different regions of $\Omega \times (0,T)$. \\
As a consequence, thanks to an approximation result and a delicate passage to the limit,
we are able to get a Harnack inequality for solutions, and in these cases only for solutions, of strongly degenerating parabolic equations,
i.e. when $\uprho \geqslant 0$. \\
As a byproduct one obtains H\"older continuity for solutions of a subclass of the first equation
(i.e. $\uprho (x,t) u_t + A u = 0$): in particular the solutions of
this subclass are H\"older continuous in the interface where $\uprho$ changes its sign, from positive to zero.
\end{abstract}



\noindent
Conflict of interest: I have no conflict of interest with anybody. \\

\noindent
Data availability statement: I state that my manuscript has no associated data.

\section{Introduction}

The main goal of this paper is to provide a Harnack type inequality for solutions of equations like
\begin{gather}
\label{equazionegenerale1}
{\displaystyle \uprho (x,t) \frac{\partial u}{\partial t} } - \textrm{div} \big(A (x,t,u,Du) \big) + B (x,t,u,Du) + C(x,t,u) = 0 	\qquad \text{in } \Omega \times (0,T)			\\
\label{equazionegenerale2}
{\displaystyle \frac{\partial}{\partial t} } \big( \uprho (x,t) u \big) - \textrm{div} \big( A (x,t,u,Du) \big) + B (x,t,u,Du) + C(x,t,u) = 0 \qquad \text{in } \Omega \times (0,T)	
\end{gather}
where the coefficient $\uprho \in L^{\infty}$ satisfies 
$$
\uprho \geqslant 0 ,
$$
i.e. it may vanish in regions of positive measure, and where the functions
\begin{align*}
& A : \Omega \times (0,T) \times \R \times \R^n \to \R^n	,	\\
& B : \Omega \times (0,T) \times \R \times \R^n \to \R ,		\\
& C : \Omega \times (0,T) \times \R \to \R
\end{align*}
are Caratheodory functions with $A, B , C \in L^{\infty} (\Omega \times (0,T) \times \R \times \R^n)$ satisfying
\begin{equation}
\label{proprieta`}
\left\{
\begin{array}{l}
\big( A(x,t,u,\xi), \xi \big) \geqslant  \lambda |\xi|^2				\\	[0.5em]
\big| A(x,t,u,\xi) \big| \leqslant  \Lambda |\xi|					\\	[0.5em]
\big| B(x,t,u,\xi) \big| \leqslant M |\xi|	\, ,						\\	[0.5em]
\big| C(x,t,u) \big| \leqslant N |u| \, ,
\end{array}
\right.
\end{equation}
for some given positive constants $\lambda, \Lambda, M, N$
and for every $u \in \R, \xi \in \R^n$ and for a.e. $(x,t) \in \Omega \times (0,T)$. 
Existence results for such equations are given in \cite{fabio4} and \cite{fabio9.1}. \\ [0.3em]
A part of the natural energy at time $t$ of a solution of equations like \eqref{equazionegenerale1} or \eqref{equazionegenerale2} is
\begin{align}
\label{carloamoremio}
\int_{\Omega} u^2 (x,t) \uprho (x,t) \, dx .
\end{align}
When $\uprho > 0$ almost everywhere this is easily handled, but when $\uprho \geqslant 0$
(i.e. $\uprho$ is null in a set of positive measure), the fact that this term may give only a partial information creates problems
and adapting the classical techniques does not work.  \\
These difficulties are evident, for instance, even in proving local boundedness of solutions to equations
\eqref{equazionegenerale1} or \eqref{equazionegenerale2}
with $\uprho \geqslant 0$ (for which we refer to \cite{fabio21}). \\ [0.3em]
For this reason the proof is delicate and technical, and we try to summarise it briefly in the following.
First, write $\Omega \times (0,T)$ as $\mathcal{Q}_+ \cup \mathcal{Q}_0 \cap \Upgamma$ where
$\mathcal{Q}_+$ and $\mathcal{Q}_0$ are open and connected subsets of $\Omega \times (0,T)$ and
$\Upgamma = \partial \mathcal{Q}_+ \cap \partial \mathcal{Q}_0$ satisfying appropriate assumptions.
The function $\uprho$ is assumed to be positive in $\mathcal{Q}_+$ and to vanish in $\mathcal{Q}_0$. \\ [0.3em]
Here we try to give a simplified idea in the simplest case, i.e. when $\uprho$ takes only two values
(in general $\uprho > 0$ almost everywhere in $\mathcal{Q}_+$), and in particular
$$
\uprho :=
\left\{
\begin{array}{ll}
1		&		\text{in } \mathcal{Q}_+	,			\\	[0.2em]
0		&		\text{in } \mathcal{Q}_0	.		
\end{array}
\right.
$$
We approximate a solution $u$ of \eqref{equazionegenerale1} or \eqref{equazionegenerale2} with solutions of standard parabolic equations like
\begin{gather}
\label{equazionegenerale1e}
{\displaystyle \uprho_{\epsilon} (x,t) \frac{\partial u}{\partial t} } - \textrm{div} \big(A (x,t,u,Du) \big) + B (x,t,u,Du) + C(x,t,u) = 0 	\qquad \text{in } \Omega \times (0,T)			\\
\label{equazionegenerale2e}
{\displaystyle \frac{\partial}{\partial t} } \big( \uprho_{\epsilon} (x,t) u \big) - \textrm{div} \big( A (x,t,u,Du) \big) + B (x,t,u,Du) + C(x,t,u) = 0 \qquad \text{in } \Omega \times (0,T)	
\end{gather}
with 
\begin{gather}
\label{rotoy}
\epsilon \in (0,1] \qquad \text{and} \qquad 
\uprho_{\epsilon} :=
\left\{
\begin{array}{ll}
1		&		\text{in } \mathcal{Q}_+	,			\\	[0.2em]
\epsilon	&		\text{in } \mathcal{Q}_0	,	
\end{array}
\right.
\end{gather}
and trying to use first standard techniques for this case, then taking the limit for $\epsilon$ going to $0$. \\
The case for equations \eqref{equazionegenerale1} or \eqref{equazionegenerale2} with a positive coefficient $\uprho$ has already been considered in
\cite{fabio19}, but there a different approach is considered. In the present paper we prove essentially the same thing, but the approach
is designed to prove the Harnack inequality when $\uprho \geqslant 0$. This is essentially the content of Part I.
Part II is devoted to prove the main result when $\uprho \geqslant 0$. \\
The proof relies on four main ingredients (for the first and the second see Theorem \ref{harnack}, for the other see Part II), of which the second and the fourth are not immediate: \\ [0.3em]
$1^{\circ} \, )$  for a function $u_{\epsilon}$ belonging to an appropriate De Giorgi class (see Definition \ref{DG1}),
we show an inhomogeneous parabolic Harnack inequality, where {\em inhomogeneous} is referred to {\em time}, since the ``amount of time one has to wait''
in one of the two regions $\mathcal{Q}_+$ or $\mathcal{Q}_0$ depends on the mean value of $\uprho_{\epsilon}$ in that region.
Consider a point $(x_o, t_o) \in \mathcal{Q}$.
For the toy model with the coefficient \eqref{rotoy}
we derive the existence of a constant $\upeta > 0$ independent of $u_{\epsilon}$ such that
\begin{equation}
\label{harnack_introduzione}
\begin{array}{cc}
{\displaystyle
\sup_{B_{r}^+ (x_o; t_o - r^2 )} u_{\epsilon} (x , t_o - r^2) \leqslant \upeta \, \inf_{B_{r}^+ (x_o; t_o + r^2)} u_{\epsilon} (x , t_o + r^2)}		\\	[0.5em]
{\displaystyle
\sup_{B_{r}^0 (x_o; t_o - \epsilon \, r^2 )} u_{\epsilon} (x , t_o - \epsilon \, r^2) \leqslant \upeta \, \inf_{B_{r}^0 (x_o; t_o + \epsilon \, r^2)} u_{\epsilon} (x , t_o + \epsilon \, r^2)} ,
\end{array}
\end{equation}
where $B_{r}^+ (\bar{x}; \bar{t}) := \big( B_r ( \bar{x}) \times \{ \bar{t} \} \big) \cap \mathcal{Q}_+$,
$B_{r}^0 (\bar{x}; \bar{t}) := \big( B_r ( \bar{x}) \times \{ \bar{t} \} \big) \cap \mathcal{Q}_0$; 		\\ [0.3em]
$2^{\circ} \, )$ paying careful attention to the constants on which $\upeta$ depends on, one realises that
(see Remark \ref{scimmiettabella} and \eqref{bottarga} in the proof of Theorem \ref{harnack0+})
$$
\upeta \text{ does not depend on } \epsilon ;
$$
$3^{\circ} \, )$ 
for $u$ solution of \eqref{equazionegenerale1} or of \eqref{equazionegenerale2} one can find $u_{\epsilon}$ solution of
\eqref{equazionegenerale1e} or of \eqref{equazionegenerale2e}
such that the family $(u_{\epsilon})_{\epsilon > 0}$ pointwise converges to $u$
for ${\epsilon} \to 0$;  \\ [0.3em]
$4^{\circ} \, )$ 
finally an important fact helps to obtain the Harnack inequality for $u$:
although $u_{\epsilon}$ converge to $u$ only pointwise one can, in some sense unexpectedly, take the limit in \eqref{harnack_introduzione}
(using Lemma \ref{lemmino}) and get
\begin{gather*}
\sup_{B_{r}^+ (x_o; t_o - r^2 )} u(x , t_o - r^2) \leqslant \upeta \, \inf_{B_{r}^+ (x_o; t_o + r^2)} u(x , t_o + r^2)	\\
\sup_{B_{r}^0 (x_o; t_o )} u(x , t_o ) \leqslant \upeta \, \inf_{B_{r}^0 (x_o; t_o)} u(x , t_o) .
\end{gather*}
\ \\
To see the precise and more general statements we refer to Theorem \ref{harnack} for the inhomogeneous Harnack inequality for $\uprho > 0$,
to Theorem \ref{harnack0+} for the Harnack inequality for $\uprho \geqslant 0$. \\ [0.3em]
As a consequence we also get some regularity for some solutions.
Precisely, adapting the classical simple argument due to Moser (see \cite{moser61}, but also, e.g., Section 7.9 in \cite{giusti}),
one can show that solutions of \eqref{equazionegenerale1} with $C = 0$ and a bounded boundary datum are H\"older continuous
in $\mathcal{Q}_+ \cup \Upgamma$, while in $\mathcal{Q}_0$
they are locally H\"older continuous with respect to $x$ for almost every $t$ (see Subsection \ref{sottoparagrafo7.1} for more details). \\
In particular the solutions are H\"older continuous in the interface $\Upgamma$. \\ [0.3em]
Equations like \eqref{equazionegenerale1} and \eqref{equazionegenerale2} arise as a natural generalisations of the case $\uprho = \uprho (x)$
to the time dependent case $\uprho = \uprho (x,t)$ and the interest for them comes from the many applications. \\
Applications when $\uprho > 0$ can be found
in some diffusion and  in fluid flow problems.
Specifically, these equations occur in hydrology when dealing with the transport of contaminants in unsaturated or variably saturated media (see \cite{demarsily}), or in density-dependent flow in porous media (see, for instance, \cite{mario}). \\
Regarding the mixed case, many diffusion problems naturally lead to problems that are partially elliptic and partially parabolic. For example, in some composite materials, there can be zones where the evolution is so rapid that the problem is considered "stationary" in those areas and evolutionary in others (i.e., when $\uprho \geqslant 0$).
For examples of equations of this type and their applications, we direct readers to \cite{carrol-show}
(see Chapter 3, in particular Example 1.1, Example 5.14, Example 5.15) and the references therein for numerous applications,
\cite{show1}, \cite{show2} (Section III.3), and \cite{zeni} (Chapter 5) . \\ [0.3em]
Coming to the Harnack inequality, there is a wide literature about that.
We confine to briefly review the highlights of the story of this inequality.
As regards the elliptic case, starting from the very first result due to Harnack for harmonic functions in 1887, 
important steps are due to J. Serrin in 1956 (see \cite{serrin55}) and to J. Moser in 1961 (see \cite{moser61}), where
for the first time weak solutions of (linear) equations in divergence form
\begin{equation}
\label{div_form}
- \textrm{div}(a \cdot D u) = 0 , \hskip30pt a_{ij} = a_{ji} , \quad c |\xi|^2 \leqslant (a \cdot \xi, \xi) \leqslant C |\xi|^2 \quad \textrm{for every } \xi \in \R^n \ (c > 0) \, .
\end{equation}
are considered.
Using the ideas of Moser, which are not strictly connected with linear equations, in 1964 Serrin (see \cite{serrin64})
proved a Harnack inequality also for solutions of a wide class of quasi-linear equations with $p-1$ growth, $p > 1$. \\
Before turning to parabolic equations we have to mention the paper of E. De Giorgi \cite{degiorgi3}, where, in particular,
H\"older continuity of the solutions of linear equations like those in \eqref{div_form} is proved.
The reason why we mention it is that the technique introduced by De Giorgi would have turned out to be useful also to prove local boundedness
and a Harnack inequality for a wide class of functions, just called later De Giorgi classes. \\
As regards parabolic equations, after the very first results due to Hadamard and Pini
(see \cite{hadamard} and \cite{pini}), the next important steps are due to Moser (see \cite{moser64}) and to Aronson and Serrin (see \cite{aron-serrin})
where linear and non-linear but with linear growth are considered. \\
Unlike the elliptic case, for parabolic operators with different growth the problem is more delicate.
A Harnack inequality for the equation $u_t - \textrm{div} \big( |Du|^{p-2} Du \big) = 0$, $p > 2$, was proved, to our knowlwdge, by DiBenedetto and can be found in \cite{dibenedetto2},
but see also \cite{Gia-Ves}, and for more general operators
by DiBenedetto, Gia\-naz\-za, Vespri in 2008 (see \cite{articolo}), after
DiBenedetto in 1986 showed that a Harnack inequality analogous to the linear case cannot hold. \\
For the case $p \in (1,2)$ we mention \cite{Gia-Ves2}. \\ [0.3em]
Just to mention some papers with degenerate coefficients we mention, for the elliptic case, \cite{fa-ke-se}, and for parabolic equations with degenerate coefficients of the type
\begin{gather*}
\uprho (x) \frac{\partial u}{\partial t} - \sum_{i, j = 1}^n \frac{\partial}{\partial x_i} \left( a_{ij} (x,t) \frac{\partial u}{\partial x_j} \right) = 0
\end{gather*}
we mention \cite{gut-wheeden1} and \cite{gut-wheeden2}, where a Harnack inequality is shown for solutions. \\
Coming to more recent papers, in 2006 Gianazza and Vespri used a technique, due also to DiBenedetto,
inspired by the already mentioned paper by De Giorgi and adapted to the parabolic case.
In \cite{gianazza-vespri} they proved a Harnack inequality not only for solutions, but for functions belonging to the natural De Giorgi class. \\
Finally, we recall the already mentioned paper \cite{fabio19}, where it is shown a Harnack inequality for functions belonging to the natural De Giorgi class, where $\uprho > 0$ almost everywhere. \\ [0.3em]
In the present paper we use this technique obtaining the inhomogeneous inequality for functions in the parabolic De Giorgi class, 
while for elliptic-parabolic equations treated in the second part we
confine ourselves to consider solutions, as we already explained. \\ [0.3em]
Finally, we conclude saying that without much additional work one could also consider unbounded $\uprho, \lambda, \Lambda, M, N$, i.e. in $L^1_{\loc}$.
Since the Part I is quite technical and the proof quite long we prefer to confine to consider bounded coefficients just to prevent unnecessary complications.

\section*{Part I - An inhomogeneous parabolic Harnack inequality}
\ \\ [-3em]
\section{Preliminaries, notations and assumptions}

Consider $\Omega \subset \R^n$ an open set, $T > 0$, and a function
\begin{align*}
\uprho \in L^{\infty} \big( \Omega \times (0,T) \big) , \qquad \uprho > 0 \text{ almost everywhere in } \Omega \times (0,T) .
\end{align*}
We consider a partition of $\mathcal{Q} := \Omega \times (0,T)$ in three subsets
\begin{align*}
\mathcal{Q} := \mathcal{Q}_1 \cup \mathcal{Q}_2 \cup \Upgamma
\end{align*}
where $\mathcal{Q}_1$ and $\mathcal{Q}_2$ are two not empty disjoint connected and open subsets of $\mathcal{Q}$,
while $\Upgamma$ is the subset of codimension $1$ defined by
$$
\Upgamma := \partial \mathcal{Q}_1 \cap \partial \mathcal{Q}_2 .
$$
About $\Upgamma$ we make some assumptions. The first one is made just for the sake of simplicity:
we suppose that for each $t \in [0,T]$ there is a partition of $\Omega$ in three not empty connected subsets ($ \Omega_j (t)$ open, $j=1,2$)
$$
\Omega = \Omega_1 (t) \cup \Omega_2 (t) \cup \Upgamma (t),
$$
with $\Gamma (t)$ of codimension $1$, that is
$$
\Upgamma (t) := \partial \Omega_1 (t) \cap \partial \Omega_2 (t) .
$$
Notice that with these requirements one has
$$
\cup_{t \in [0,T]} \Upgamma (t) = \Upgamma
$$
and we avoid a situation like the one shown in the following picture:
\ \\
\begin{center}
\pgfplotsset{every axis/.append style={
extra description/.code={ \node at (0.77,0.42) {${\mathcal Q}_2$}; \node at (0.54, 0.55) {$\Upgamma$}; \node at (0.3,0.6) {${\mathcal Q}_1$};
}}}
\begin{tikzpicture}
\begin{axis} [axis equal, axis lines=middle, xmin=-0.4, xmax=7.8, ymin=0, ymax=4.8, 
xtick={0}, ytick={0}, xticklabels={$$}, yticklabels={$$}, xlabel=$x$, ylabel=$t$, width=7.4cm, height=6cm, title={Figure 1}]
\addplot
[domain=1:7,variable=\t, very thick]
({t},{4.2});
\addplot
[domain=1:7,variable=\t, very thick]
({t},{0});
\addplot
[domain=0:4.2,variable=\t, very thick]
({1},{t});
\addplot
[domain=0:4.2,variable=\t, very thick]
({7},{t});
\addplot
[domain=0:0.7,variable=\t, thick]
({t+1.5+1},{3*t});
\addplot
[domain=0.7:1.4,variable=\t, thick]
({t+3+1},{3*t});
\addplot
[domain=0:1.5,variable=\t, thick]
({t+2.2+1},{1.1+1});
\end{axis}
\end{tikzpicture}
\end{center}
\ \\
Notice that in this way we also have that
$$
\mathcal{Q}_1 = \cup_{t \in [0,T]} \Omega_1 (t) , \quad \mathcal{Q}_2 = \cup_{t \in [0,T]} \Omega_2 (t) .
$$
Other assumptions about $\Upgamma$ are listed below in (H.1)-(H.4). \\ [0.3em]
For every set $E \subset \Omega \times (0,T)$ and ball $B_r ( \bar{x} ) \subset \Omega$ we will denote by
\begin{gather*}
E (\bar{t}) := \{ (x,t) \in E \, | \, t = \bar{t} \}											\\
E^j (\bar{t}) := E (\bar{t}) \, \cap \, \mathcal{Q}_j , \qquad j =1, 2 ,						\\
B_r ( \bar{x} ; \bar{t} ) := B_r ( \bar{x} ) \times \{ \bar{t} \} ,								\\
B_r^j ( \bar{x} ; \bar{t} ) := B_r ( \bar{x} ; \bar{t} ) \, \cap \, \mathcal{Q}_j , \qquad j =1, 2 .
\end{gather*}
We will denote
\begin{align}
\label{iduero}
\uprho_1 (x,t) := 
\left\{
\begin{array}{cl}
\uprho		&	\text{if } (x,t) \in \mathcal{Q}_1	,	\\	[0.3em]
	0		&	\text{if } (x,t) \in \mathcal{Q}_2 ,
\end{array}
\right.
\qquad
\uprho_2 (x,t) := 
\left\{
\begin{array}{cl}
	0		&	\text{if } (x,t) \in \mathcal{Q}_1	,	\\	[0.3em]
\uprho		&	\text{if } (x,t) \in \mathcal{Q}_2 ,
\end{array}
\right.
\end{align}
and alongside
\begin{align}
\label{chi}
\upchi_1 (x,t) := 
\left\{
\begin{array}{cl}
1		&	\text{if } (x,t) \in \mathcal{Q}_1	,	\\	[0.3em]
0		&	\text{if } (x,t) \in \mathcal{Q}_2 ,
\end{array}
\right.
\qquad
\upchi_2 (x,t) := 
\left\{
\begin{array}{cl}
0		&	\text{if } (x,t) \in \mathcal{Q}_1	,	\\	[0.3em]
1		&	\text{if } (x,t) \in \mathcal{Q}_2 .
\end{array}
\right.
\end{align}
We will use the same notations for the functions $\uprho$, $\uprho_j$, $\upchi_j$ and for the measures they induce, for instance we will write
\begin{gather*}
\uprho (E) := \iint_E \uprho (x,t) \, dx dt , \qquad E \subset \Omega \times (0,T),		\\
\uprho (t) (A) := \int_A \uprho (x,t) \, dx , \qquad A \subset \Omega .
\end{gather*}
Notice that $\upchi_1 (E) + \upchi_2 (E) = |E|$ for $E \subset \mathcal{Q}$.
By $Du$ we will always denote the vector of derivatives of a function $u$ with respect to the variables $x_1, \dots x_n$, $(x_1, \dots x_n) \in \Omega$, even if $u$ depends also on $t$.
Once defined, for $A$ open subset of $\R^k$ and $\nu : A \to [0, M]$
\begin{align}
\label{elleduepesato}
L^2 \big( A, \nu\big) := \text{the completion of } C^1_c (A) \text{ w.r.t. the topology induced by } \left( \int_{A} w^2 \nu \, dy \right)^{1/2} ,
\end{align}
we then define
$$
\V := \big\{ u \in L^2 \big( \Omega \times (0,T), \uprho \big) \, \big| \, Du \in L^2 (\Omega \times (0,T)) \big \}
$$
and
$$
\V' \qquad \text{the dual space of } \V .
$$
We also define
$$
\V_{\loc} := \big\{ u \in L^2_{\loc} \big( \Omega \times (0,T), \uprho \big) \, \big| \, Du \in L^2_{\loc} (\Omega \times (0,T)) \big \} .
$$

\noindent
In \cite{fabio4} and \cite{fabio9.1} some existence results are proved (where $\uprho$ may be also zero and negative).
As byproducts we get the existence of solutions for the following equations
($A$, $B$ and $C$ introduced in the previous section)
\begin{gather*}
{\displaystyle \uprho (x,t) \frac{\partial u}{\partial t} } - \textrm{div} \big(A (x,t,u,Du) \big) + B (x,t,u,Du) + C(x,t,u) = 0 ,					\\
{\displaystyle \frac{\partial}{\partial t} } \big( \uprho (x,t) u \big) - \textrm{div} \big( A (x,t,u,Du) \big) + B (x,t,u,Du) + C(x,t,u) = 0 .
\end{gather*}

\noindent
About $\uprho$ we will need some assumptions; before to list them we recall some definitions.

\bd
\label{BeA}
Given $\omega \in L^{\infty} (\R^n)$, 
$\omega > 0$ a.e., $K > 0$ and $q > 2$ we say that
$$
\omega \in B_{2,q}^{1} (K)
$$
if
\begin{equation}
\label{sob-poin-cond}
\frac{r}{\rho} \left( \frac{\omega (B_r(\bar{x}))}{\omega (B_\rho(\bar{x}))} \right)^{1/q} \left( \frac{| B_r (\bar{x})|}{|B_\rho(\bar{x})|} \right)^{-1/2} \leqslant K
\end{equation}
for every pair of concentric balls $B_r$ and $B_\rho$ with $0 < r < \rho$. \\
We say that $\omega$ is {\em doubling} if
for every $t > 0$
there is a positive constant $c_{\omega}(t)$ such that
\begin{equation}
\label{doubling}
\int_{t B} \omega (x)  \, dx \leqslant c_{\omega} (t) \int_B \omega (x) \, dx
\end{equation}
where $B$ is a generic ball $B_{r} (\bar{x})$ and $t B$ denotes the ball $B_{t r} (\bar{x})$. \\
We say that $\omega$ belongs to the class $A_{\infty} (K, \varsigma)$, $\varsigma > 0$, if
\begin{equation}
\label{Ainfinito}
\frac{\omega (S)}{\omega (B)} \leqslant K \left( \frac{|S|}{|B|} \right)^{\varsigma}
\end{equation}
for every ball $B \subset \R^n$ and every measurable set $S \subset B$. \\
We will denote $A_{\infty} (K)$ the set $\cup_{\varsigma > 0} A_{\infty} (K, \varsigma)$ and $A_{\infty} = \cup_{K > 0} A_{\infty} (K)$.
\ed

\boss
\label{rimi}
If $\omega \in A_{\infty}$ then $\omega$ is doubling (see \cite{gc-rdf}).
\eoss

\boss
\label{rimi2}
If $\omega \in A_{\infty} (K)$ then there is $b \geqslant 1$ such that
\begin{align*}
\left(\frac{|S|}{|B|}\right)^b \leqslant K \, \frac{\omega(S)}{\omega(B)}
\end{align*}
for every measurable $S \subset B$, for every $B$ ball of $\R^n$ (see \cite{gc-rdf}).
\eoss

\noindent
{\bf Assumption about the function $\uprho$ - }
We suppose that 
$$
\uprho \in L^{\infty} (\Omega \times (0,T))
$$
and there is a positive constant $L$ such that
\begin{equation}
\tag{\bf H.1}
\begin{array}{l}
[0,T] \ni t \mapsto {\displaystyle \int_{\Omega} v(x) w(x) \uprho (x,t) \, dx } \qquad \text{is absolutely continuous} \hskip30pt \text{and}		\\	[0.5em]
{\displaystyle
\Bigg| \frac{d}{dt}  \int_{\Omega} v(x) w(x) \uprho (x,t) \, dx \Bigg| \leqslant L \, \| v \|_{H^{1} (\Omega)} \| w \|_{H^{1} (\Omega)}
}
\end{array}
\end{equation}
for every $v, w \in H^1 (\Omega)$.  \\ [0.2em]

\noindent
Moreover we will also suppose that
$$
\Upgamma := \partial \mathcal{Q}_1 \cap \partial \mathcal{Q}_2
$$
satisfies the following assumption: 
there are two constants $\bar{R} > 0$ (in this regard, see the example in Figure 3, that shows why we need the existence of such $\bar{R}$)
and $\upkappa \in (0,1/2)$ such that for every $j = 1, 2$
\begin{equation}
\tag{\bf H.2}
\left|
\begin{array}{ll}
 i \, ) & \uprho_j (t_o) \big( B_r (x_o) \big) \geqslant 2 \upkappa \, \uprho (t_o)  \big( B_r (x_o) \big) , 			\\	[0.7em]
 ii \, ) & \upchi_j (t_o) \big( B_r (x_o) \big) \geqslant 2 \upkappa \, \big| B_r (x_o) \big|
\end{array}
\right.
\end{equation}
for every $(x_o,t_o) \in \mathcal{Q}_j \cup \Upgamma$ and every $r \in (0, \bar{R}]$ such that $B_r (x_o) \subset \Omega$. \\
\ \\
Moreover we will suppose the existence of $\updelta > 0$ such that for every $j = 1, 2$
\begin{equation}
\tag{\bf H.3}
\left|
\begin{array}{ll}
 i \, ) & \big| \uprho_j (t) \big( B_r (x_o) \big) - \uprho_j (s) \big( B_r (x_o) \big) \big| \leqslant \upkappa \, {\uprho} (t_o)  \big( B_r (x_o) \big), 		\\	[0.7em]
 ii \, ) & \big| \upchi_j (t) \big( B_r (x_o) \big) - \upchi_j (s) \big( B_r (x_o) \big) \big|  \leqslant \upkappa \, \big| B_r (x_o) \big| , 	
 \end{array}
\right.
\end{equation}
for every $(x_o,t_o) \in \mathcal{Q}_j \cup \Upgamma$, every $r \in (0, \bar{R}]$ such that $B_r (x_o) \subset \Omega$,
every $s, t \in [0,T] \cap [t_o - \updelta, t_o + \updelta]$. \\ [0.3em]
Finally we will suppose
\begin{equation}
\tag{\bf H.4}
\left|
\begin{array}{ll}
 i \, ) & \uprho (t) \in B_{2,q}^{1} (K_1) \quad \text{for some } q > 2 \qquad \text{ for every } t \in (0,T) , 				\\	[0.6em]
 ii \, ) & \uprho (t) \in A_{\infty} (K_2, \varsigma) \quad \quad \text{for some } \varsigma \in (0,1], \qquad \text{ for every } t \in [0,T] .
\end{array}
\right.
\end{equation}

\boss
\label{nota}
We could weaken assumptions (H.2) and (H.3) as follows: we could assume that for every $(t_o, x_o) \in \mathcal{Q}_j \cup \Upgamma$ there are $\bar{R}$, $\upkappa$, $\updelta$
for which (H.2) and (H.3). Every result that follows would be true, but with the constants $\bar{R}$, $\upkappa$, $\updelta$ depending on $t_o, x_o$. \\
For the sake of simplicity we assume $\bar{R}$, $\upkappa$, $\updelta$ uniform in $\mathcal{Q}$.
The same thing we do for the constant ${\sf L}$ in (C.6). Notice that if $\Omega$ is bounded we can drop the dependence of $\bar{R}$, $\upkappa$, $\updelta$ on $x_o$.
\eoss

\noindent
{\bf Comments on the assumptions.} \\ [0.3em]
{\bf (H.1)} - This assumption is needed in \cite{fabio4} and \cite{fabio9.1} to obtain the results about existence and uniqueness of the solutions of
\eqref{equazionegenerale1} and \eqref{equazionegenerale2} with suitable boundary conditions. 
Boundedness could be avoided, but, as already said in the introduction, we prefer to assume $\uprho$ bounded, thus keeping the exposition from becoming overly complex.
For examples with unbounded $\uprho$ we refer to \cite{fabio4} and \cite{fabio9.1}. \\
This assumption is clearly satisfied if
\begin{equation}
\label{doesnotsatisfy}
\uprho,  {\displaystyle \frac{\partial \uprho}{\partial t} } \in L^{\infty} (\Omega \times (0,T))
\end{equation}
since in this case we get
\begin{align*}
\Bigg| \frac{d}{dt}  \int_{\Omega} & v(x) w(x) \uprho (x,t) \, dx \Bigg| = 
		\Bigg| \int_{\Omega} v(x) w(x) \frac{\partial \uprho}{\partial t} (x,t) \, dx \Bigg| \leqslant \| \uprho_t \|_{\infty} \, \| v \|_{L^2 (\Omega)} \| w \|_{L^2 (\Omega)} \, .
\end{align*}
But more interesting are the cases when $t \mapsto \uprho (x,t)$ is not absolutely continuous, it may be only continuous
(see an example in \cite{fabio13}) and even discontinuous for every $x \in \Omega$.
The simplest example is the following: consider ($\Omega = (0,1)$, $T = 1$)
$$
\uprho : (0,1) \times (0,1) \to \R \quad \text{s.t. } \qquad 
\uprho (x,t) := \left\{
\begin{array}{ll}
2			&	\text{ for } x < t				\\	[0.3em]
1			&	\text{ for } x > t
\end{array}
\right. .
$$
Then
\begin{align*}
\frac{d}{dt}  \int_{\Omega} v(x) w(x) \uprho (x,t) \, dx = \frac{d}{dt} \left( 2 \int_0^t v(x) w(x) \, dx +  \int_t^1 v(x) w(x) \, dx \right) = v(x) w(x) .
\end{align*}
Since $H^1(0,1) \subset C^0 ([0,1])$, one can estimate the left hand side in the previous equality by the $H^1$ norm of $v$ and $w$. \\
Analogous examples may be consider in higher dimension, for which we refer to \cite{fabio4} and \cite{fabio9.1}. \\
Also $\uprho \geqslant 0$ may be considered, as we will do in Part II, and indeed the same holds if we consider
$$
\uprho : (0,1) \times (0,1) \to \R \quad \text{s.t. } \qquad 
\uprho (x,t) := \left\{
\begin{array}{ll}
1			&	\text{ for } x < t				\\	[0.3em]
0			&	\text{ for } x > t
\end{array}
\right. .
$$
\\ [0.3em]
{\bf (H.4)} is needed for Theorem \ref{gutierrez-wheeden} to hold. \\
Moreover {\bf (H.4)} - $ii \, )$ implies (C.10) stated below, that
is needed in Lemma \ref{lemma2} and consequently for Theorem \ref{esp_positivita}. \\ [0.3em]
Assumptions {\bf (H.2)} are essentially assumptions on $\Upgamma$:
here we show some pictures to clarify. \\
As observed in (C.7) below, the function $\uprho$ restricted to each of $\mathcal{Q}_1, \mathcal{Q}_2$ is doubling, in particular in the balls centred in a point belonging to $\Upgamma$. \\
Observe that $\Upgamma$ may have also cusps, as shown in Figure 2:
consider, for instance, $\Upgamma$
which is the union of the graphs of $f(x) = x^n$ and $g(x) = - x^n$
for $x > 0$ and $n \in \N$, $n \geqslant 1$. In this case we have that
$$
\uprho_1 \big( {B_{2\rho} (\bar{x}, \bar{t})} \big) \leqslant c \, \uprho_1 \big( {B_{\rho} (\bar{x}, \bar{t})} \big) , \qquad
\uprho_2 \big( {B_{2\rho} (\bar{x}, \bar{t})} \big) \leqslant c \, \uprho_2 \big( {B_{\rho} (\bar{x}, \bar{t})} \big)
$$
for some $c$ depending on $n$. Notice that not every kind of cusp is admitted (see \cite{fabio10}).
\ \\
\begin{center}
\pgfplotsset{every axis/.append style={
extra description/.code={\node at (0.35, 0.4) {$(\bar{x}, \bar{t})$}; \node at (0.64, 0.88) {$\Upgamma$}; \node at (0.77,0.62) {$\uprho_2$};  \node at (0.2,0.6) {$\uprho_1$};
}}}
\begin{tikzpicture}
\begin{axis} [axis equal, xtick={0}, ytick={0}, xticklabels={$$}, yticklabels={$$}, width=8cm, height=6cm, title={Figure 2}]
\addplot coordinates
{(0, 0)};
\addplot
[domain=0:0.9,variable=\t,
samples=40, smooth, thick, black]
({t},{t^3});
\addplot
[domain=0.9:1.1,variable=\t,
samples=40, smooth, thick, dashed]
({t},{t^3});
\addplot
[domain=0:0.9,variable=\t,
samples=40, smooth, thick]
({t},{-t*t*t});
\addplot
[domain=0.9:1.1,variable=\t,
samples=40, smooth, thick, dashed]
({t},{-t*t*t});
\end{axis}
\end{tikzpicture}
\end{center}
\ \\
\noindent
The request of the existence of $\bar{R}$ for which (H.2) holds is due to the fact that we want to admit also interfaces as that shown in Figure 3.
If such $\bar{R}$ did not exist (H.2), in the following example, could not hold for every $r > 0$.
\ \\
\begin{center}
\pgfplotsset{every axis/.append style={
extra description/.code={ \node at (0.77,0.5) {${\mathcal Q}_2$}; \node at (0.54, 0.63) {$\Upgamma$}; \node at (0.2,0.6) {${\mathcal Q}_1$};
}}}
\begin{tikzpicture}
\begin{axis} [axis equal, xtick={2}, ytick={1.2}, xticklabels={$$}, yticklabels={$$}, width=14cm, height=6cm, title={Figure 3}]
\addplot
[domain=-2:8,variable=\t, very thick]
({t},{3.2});
\addplot
[domain=-3.5:-2,variable=\t, very thick, dashed]
({t},{3.2});
\addplot
[domain=8:9.5,variable=\t, very thick, dashed]
({t},{3.2});
\addplot
[domain=-2:8,variable=\t, very thick]
({t},{-1});
\addplot
[domain=-3.5:-2,variable=\t, very thick, dashed]
({t},{-1});
\addplot
[domain=8:9.5,variable=\t, very thick, dashed]
({t},{-1});
\addplot
[domain=0:8,variable=\t, thick]
({t},{1.5});
\addplot
[domain=8:9.5,variable=\t, thick, dashed]
({t},{1.5});
\addplot
[domain=0:8,variable=\t, thick]
({t},{0.7});
\addplot
[domain=0:9.5,variable=\t, thick, dashed]
({t},{0.7});
\addplot
[domain=1.57:4.71,variable=\t, thick]
({0.4*cos(deg(t))},{0.4*sin(deg(t))+1.1});
\end{axis}
\end{tikzpicture}
\end{center}
\ \\
Finally assumptions {\bf (H.3)} is needed to ``transport'' assumptions {\bf (H.2)} for a while (see, e.g, the consequences (C.3) and (C.5)). \\ [0.5em]

\noindent
{\bf Consequences - } Here we show some consequences of (H.1)-(H.4). \\ [-0.5em]
\begin{itemize}
\item[{\bf (C.1)}]
If $\uprho$ satisfies (H.1) then
$$
[0,T] \ni t \mapsto \uprho (t) \big( B_r \big) \qquad \text{is continuous}
$$
for every ball $B_r \subset \Omega$ (see Lemma 2.4 in \cite{fabio19}), and then uniformly continuous. \\
\item[{\bf (C.2)}]
The weights $\{\uprho (t)\}_{t \in [0,T]}$ are {\em doubling}, thanks to Remark \ref{rimi} and assumption (H.4) - $ii\, $).
So there is a constant $c_{\uprho} > 0$ (independent of $t$) such that
$$
\uprho (t) \big( B_{2 r} (x) \big) \leqslant c_{\uprho} \, \uprho (t) \big( B_r (x) \big)
$$
for every $x \in \Omega$ and $r > 0$ such that $B_{2 r} (x) \subset \Omega$. \\
In general we have
$
\uprho (t) \big( B_{R} (x) \big) \leqslant c_{\uprho} \left( \frac{R}{r} \right) \uprho (t) \big( B_r (x) \big)
$
where the constant $c_{\uprho} \left( \frac{R}{r} \right)$ depends on the ratio $R/r$ and we will simply denote it by $c_{\uprho}$ if $R = 2r$. \\
\item[{\bf (C.3)}]
For every $j \in \{1, 2\}$ and every $(x_o, t_o) \in \mathcal{Q}_j \cup \Upgamma$ and every $r \in (0, \bar{R}]$
$$
\begin{array}{r}
\uprho_j (t) \big( B_r (x_o) \big) \geqslant \upkappa \, \uprho (t_o)  \big( B_r (x_o) \big) 		\qquad
\text{for every } t \in [0,T] \text{ such that } |t - t_o| \leqslant \updelta .
\end{array}
$$
Indeed notice that, taking $s = t_o$ in (H.3) - $i \, )$ and using (H.2) - $i \, )$, we have
\begin{align*}
\uprho_j (t) \big( B_r (x_o) \big)	& \geqslant \uprho_j (t_o) \big( B_r (x_o) \big) - \upkappa \, \uprho (t_o)  \big( B_r (x_o) \big) \geqslant \upkappa \, \uprho (t_o)  \big( B_r (x_o) \big)
\end{align*}
for every $t \in [0,T]$ such that $|t - t_o| \leqslant \updelta$, for every $r \in (0, \bar{R}]$. \\
\item[{\bf (C.4)}]
As just done in (C.3), for every $j \in \{1, 2\}$, every $(x_o, t_o) \in \mathcal{Q}_j \cup \Upgamma$,
every $r \in (0, \bar{R}]$, every $t \in [0,T]$ such that $|t - t_o| < \updelta$
$$
\begin{array}{l}
\upchi_j (t) \big( B_r (x_o) \big) \geqslant \upkappa \, \big| B_r (x_o) \big| 		.
\end{array}
$$
\item[{\bf (C.5)}]
For every $j = 1, 2$, $(x_o,t_o) \in \mathcal{Q}_j \cup \Upgamma$ and $r \in (0, \bar{R}]$ with $B_r (x_o) \subset \Omega$
$$
\begin{array}{c}
{\displaystyle
\sup_{t \in (t_o - \updelta, t_o + \updelta) \cap [0,T]} \uprho_j (t) \big( B_r (x_o) \big) \leqslant 
			\frac{1 + \upkappa}{\upkappa} \, \inf_{t \in (t_o - \updelta, t_o + \updelta) \cap [0,T]} \uprho_j (t) \big( B_r (x_o) \big), }				\\	[0.5em]
{\displaystyle
\sup_{t \in (t_o - \updelta, t_o + \updelta) \cap [0,T]} \Big(\uprho_j (t) \big( B_r (x_o) \big) \Big)^{-1} \leqslant 
			\frac{1 + \upkappa}{\upkappa} \, \inf_{t \in (t_o - \updelta, t_o + \updelta) \cap [0,T]} \Big(\uprho_j (t) \big( B_r (x_o) \big) \Big)^{-1} .
}
\end{array}
$$
To prove these inequalities observe that by (H.2) - $i \, )$ and (H.3) - $i \, )$ one has
\begin{gather*}
\uprho_j (s_1) \big( B_r (x_o) \big) < \uprho_j (t_o) \big( B_r (x_o) \big) + \upkappa \, \uprho (t_o)  \big( B_r (x_o) \big) \leqslant (1 + \upkappa) \, \uprho (t_o)  \big( B_r (x_o) \big) ,		\\
\uprho_j (s_2) \big( B_r (x_o) \big) > \uprho_j (t_o) \big( B_r (x_o) \big) - \upkappa \, \uprho (t_o)  \big( B_r (x_o) \big) \geqslant \upkappa \, \uprho (t_o)  \big( B_r (x_o) \big) .
\end{gather*}
for every $s_1, s_2 \in [t_o - \updelta, t_o + \updelta]$. Then in particular
$$
\uprho_j (s_1) \big( B_r (x_o) \big) \leqslant
\frac{1 + \upkappa}{\upkappa} \, \upkappa \, \uprho (t_o)  \big( B_r (x_o) \big) \leqslant \frac{1 + \upkappa}{\upkappa} \, \uprho_j (s_2) \big( B_r (x_o) \big) ,
$$
from which the thesis follows. In an analogous way one proves
$$
\begin{array}{c}
{\displaystyle
\sup_{t \in (t_o - \updelta, t_o + \updelta) \cap [0,T]} \upchi_j (t) \big( B_r (x_o) \big) \leqslant 
			\frac{1 + \upkappa}{\upkappa} \, \inf_{t \in (t_o - \updelta, t_o + \updelta) \cap [0,T]} \upchi_j (t) \big( B_r (x_o) \big), }				\\	[0.5em]
{\displaystyle
\sup_{t \in (t_o - \updelta, t_o + \updelta) \cap [0,T]} \Big(\upchi_j (t) \big( B_r (x_o) \big) \Big)^{-1} \leqslant 
			\frac{1 + \upkappa}{\upkappa} \, \inf_{t \in (t_o - \updelta, t_o + \updelta) \cap [0,T]} \Big(\upchi_j (t) \big( B_r (x_o) \big) \Big)^{-1} .
}
\end{array}
$$
\item[{\bf (C.6)}]
We also have, by (C.1), that for each ball $B = B_r (x_o)$ with $r \in (0, \bar{R}]$ and every $t_o \in (0,T)$ there is a constant ${\sf L}$
$$
\begin{array}{c}
{\displaystyle
\sup_{t \in (0,T)} \uprho (t) \big( B_r (x_o) \big) \leqslant 
			{\sf L} \, \inf_{t \in (t_o - \updelta, t_o + \updelta)} \uprho (t) \big( B_r (x_o) \big), }				\\	[0.5em]
{\displaystyle
\sup_{t \in (0,T)} \Big(\uprho (t) \big( B_r (x_o) \big) \Big)^{-1} \leqslant 
			{\sf L} \, \inf_{t \in (t_o - \updelta, t_o + \updelta)} \Big(\uprho (t) \big( B_r (x_o) \big) \Big)^{-1} .
}
\end{array}
$$
Notice that ${\sf L}$ in principle depends on $x_o$ (and on $\bar{R}$), but, as observed in Remark \ref{nota}, this would not change the results in the following section,
except for the fact that the constant ${\sf L}$ would depend on $x_o$. For this reason, just to simplify the notation, we assume that $\bar{R}$ and ${\sf L}$
are independent of the point $x_o$. \\
\item[{\bf (C.7)}]
The functions $\uprho_j (t)$ satisfy a doubling condition, limited to balls {\em centred} in points $(x_o, t_o)$ belonging to $\mathcal{Q}_j \cup \Upgamma$.
For every $(x_o, t_o) \in \mathcal{Q}_j \cup \Upgamma$ and $r > 0$ with $2r \leqslant \bar{R}$, such that $B_{2r} (x_o) \subset \Omega$ we have that
$$
\uprho_j (t_o)  \big( B_{R} (x_o) \big) \leqslant c_{\uprho_j} \, \uprho_j (t_o) \big( B_r (x_o) \big) , \qquad c_{\uprho_j} = {\displaystyle \frac{c_{\uprho}}{2 \upkappa}} .
$$
Indeed by (H.2) - $i \, )$ and (C.2)
\begin{align*}
\uprho_j (t_o) \big( B_r (x_o) \big) \geqslant 2 \upkappa \, \uprho (t_o)  \big( B_r (x_o) \big)
	\geqslant \frac{2 \upkappa}{c_{\uprho}} \, \uprho (t_o)  \big( B_{2r} (x_o) \big) .
\end{align*}
Notice that similarly
\begin{gather*}
\upchi_j (t_o) \big( B_{R} (x_o) \big) \leqslant 2^n \big| B_{r} (x_o) \big| \leqslant
	\frac{1}{2 \upkappa} 2^n \, \upchi_j (t_o) \big( B_r (x_o) \big) = \frac{2^{n-1}}{\upkappa} \, \upchi_j (t_o) \big( B_r (x_o) \big) .
\end{gather*}
Analogous estimates hold for $R \leqslant \bar{R}$ in the place of $2r$ with $c_{\uprho_j} (R/r) = {\displaystyle \frac{c_{\uprho} (R/r)}{2 \upkappa}}$ in the first
and $\frac{1}{2 \upkappa} \left( \frac{R}{r} \right)^n$ in the second. \\
\item[{\bf (C.8)}]
Consider $r \in (0, \bar{R}]$, $j \in \{ 1, 2 \}$, $(x_o, t_o) \in \mathcal{Q}_j \cup \Upgamma$. Then by (C.3) and (C.5) one has that
\begin{align*}
\uprho_j (t) \big( B_r (x_o) \big) \geqslant \frac{\upkappa}{\sf L} \, \uprho (t)  \big( B_r (x_o) \big) \qquad	
			\text{for every } t \in [t_o - \updelta, t_o + \updelta] .
\end{align*}
\item[{\bf (C.9)}]
Given $(x_o , t_o) \in \mathcal{Q}_j \cup \Upgamma$ for $j \in \{1, 2 \}$, $\tilde\theta , \theta > 0$, $\tilde{r} , r > 0$ such that $B_{\tilde{r}} (x_o) \subset \Omega$, $B_{r} (x_o) \subset \Omega$ and
$\theta \, r^2 \leqslant \updelta$ and $\tilde\theta \, \tilde{r}^2 \leqslant \updelta$ one has that
\begin{gather*}
\uprho_j \big(B_{\tilde{r}} (x_o) \times (t_o , t_o + \tilde\theta \, \tilde{r}^2) \big) \leqslant 
	C_{\uprho_j} \, \uprho_j \big( B_{r} (x_o) \times (t_o , t_o + \theta \, r^2 ) \big) \, ,															\\
\upchi_j \big(B_{\tilde{r}} (x_o) \times (t_o , t_o + \tilde\theta \, \tilde{r}^2 ) \big) \leqslant
	\frac{1}{\upkappa} \, \frac{\tilde{\theta}}{\theta} \left( \frac{\tilde{r}}{r} \right)^{n+2} \upchi_j \big( B_{r} (x_o) \times (t_o , t_o + \theta \, r^2 ) \big)
\end{gather*}
where $C_{\uprho_j} = \frac{1 + \upkappa}{\upkappa} \, c_{\uprho_j} (\frac{\tilde{r}}{r}) \, \frac{\tilde{\theta}}{\theta} \left( \frac{\tilde{r}}{r} \right)^2$
(the constant $c_{\uprho_j} (\tilde{r}/r)$ is $1$ if $\tilde{r} \leqslant r$). \\
%
%
We prove the first: by (C.5)
\begin{align*}
\uprho_j \big( B_{r} (x_o) & \times (t_o , t_o + \theta r^2 ) \big) \geqslant \theta \, r^2 \inf_{t \in (t_o , t_o + \theta r^2 )} \uprho_j  (t) \big( B_{r} (x_o) \big) \geqslant						\\
& \geqslant \frac{\upkappa}{1 + \upkappa} \, \theta \, r^2 \sup_{t \in (t_o , t_o + \tilde\theta \, \tilde{r}^2 )} \uprho_j  (t) \big( B_{r} (x_o) \big) \geqslant 									\\
& \geqslant \frac{\upkappa}{1 + \upkappa} \, \frac{1}{c_{\uprho_j } (\tilde{r}/r)} \, \theta \, r^2 \sup_{t \in (t_o , t_o + \tilde\theta \, \tilde{r}^2 )} \uprho_j  (t) \big( B_{\tilde{r}} (x_o) \big) \geqslant	\\
& \geqslant \frac{\upkappa}{1 + \upkappa} \, \frac{1}{c_{\uprho_j } (\tilde{r}/r)} \, \frac{\theta}{\tilde{\theta}} \left( \frac{r}{\tilde{r}} \right)^2
	\uprho_j  \big( B_{\tilde{r}} (x_o) \times (t_o , t_o + \tilde\theta \, \tilde{r}^2 ) \big) .
\end{align*}
For the other we can estimate as follows:
\begin{align*}
\upchi_j \big( B_{r} (x_o) & \times (t_o , t_o + \theta r^2 ) \big) \geqslant 
	\theta \, r^2 \inf_{t \in (t_o , t_o + \theta r^2 )} \upchi_j (t) \big( B_{r} (x_o) \big) \geqslant								 						\\
& \geqslant \upkappa \, \theta \, r^2 \big| B_{r} (x_o) \big| = 
	\upkappa \, \theta \, \frac{r^2}{\tilde{r}^2 \tilde\theta} \, \tilde\theta \, \tilde{r}^2 \,  \left( \frac{r}{\tilde{r}} \right)^n \,  \big| B_{\tilde{r}} (x_o) \big| \geqslant			\\
& \geqslant \upkappa \, \left( \frac{r}{\tilde{r}} \right)^{n+2} \, \frac{\theta}{\tilde\theta} \, \tilde\theta \, \tilde{r}^2
	\sup_{t \in (t_o , t_o + \tilde\theta \, \tilde{r}^2 )} \upchi_j (t) \big( B_{\tilde{r}} (x_o) \big) \geqslant												\\
& \geqslant \upkappa \, \left( \frac{r}{\tilde{r}} \right)^{n+2} \, \frac{\theta}{\tilde\theta}
	\, \upchi_j \big( B_{\tilde{r}} (x_o) \times (t_o , t_o + \tilde\theta \, \tilde{r}^2 ) \big) .
\end{align*}
\item[{\bf (C.10)}]
Consider an interval $I \subseteq (0,T)$, a ball $B = B_r (\bar{x}) \subset \Omega$ with $r \in (0, \bar{R}]$ and a family of measurable sets $S (t) \subset B$ in such a way that
$S = \cup_{t \in I} S (t) \times \{ t \} \subset B \times I$ turns out to be measurable.
First of all observe that by (C.6)
\begin{align*}
\frac{1}{\uprho (t) (B)} \leqslant {\sf L} \, | I | \, \frac{1}{\uprho (B \times I )} \, .
\end{align*}
Now by (H.4)-$ii \, )$ and Remark \ref{rimi2} we have
\begin{align*}
\left(\frac{|S (t)|}{|B|}\right)^\ess \leqslant K_2 \, \frac{\uprho (t) (S (t))}{\uprho (t) (B)} \qquad \text{for every } t \in I
\end{align*}
with $b \geqslant 1$ independent of $t$.
Now integrating the above inequality
\begin{align*}
\frac{|S|^b}{|B|^b} & =
\frac{1}{|B|^b} \left( \int_I | S (t) | \, dt \right)^{b} \leqslant | I |^{b - 1} \int_I \frac{| S (t) |^b}{|B|^b}\, dt \leqslant | I |^{b - 1} K_2 \int_I \frac{\uprho (t) (S (t))}{\uprho (t) (B)} \, dt \leqslant			\\
	& \leqslant {\sf L} \, K_2 \, | I |^{b} \, \frac{1}{\uprho (B \times I )} \int_I \uprho (t) (S (t)) \, dt =
		{\sf L} \, K_2 \, | I |^b \, \frac{\uprho (S)}{\uprho (B \times I )}
\end{align*}
by which, dividing by $|I|^b$, 
\begin{align*}
\frac{|S|^b}{|B \times I|^b} \leqslant {\sf L} \, K_2 \, \frac{\uprho (S)}{\uprho (B \times I )} .
\end{align*}
Similarly, by (H.4)-$ii \, )$,
\begin{align*}
\frac{\uprho (S)}{\uprho (B \times I )} \leqslant {\sf L} \, K_2 \, \frac{|S|^{\varsigma}}{|B \times I|^{\varsigma}} .
\end{align*}
\end{itemize}

\ \\
\noindent
We now recall two standard results. For the following result see \cite{chanillo-wheeden}.

\bthm
\label{chanillo-wheeden}
Consider $q > 2$, $r > 0$, $x_o \in \R^n$, $\omega \in L^{\infty} (\R^n)$, $\omega > 0$ a.e., 
$\omega \in B_{2,q}^{1} (K_1)$ and doubling, i.e. satisfying \eqref{sob-poin-cond} and \eqref{doubling}.
Then there is a constant $\tilde\upgamma$ depending $($only$)$ on $n, q, K_1, c_{\omega}$ such that
\begin{equation}
\label{disuguaglianzadisobolev}
\bigg[ \frac{1}{\omega(B_r)} {\int_{B_r}} |u(x)|^q \omega (x) dx \bigg]^{1/q} \leqslant \tilde\upgamma \, r \,
	\bigg[ \frac{1}{|B_r|} {\int_{B_r}} |Du(x)|^2 dx \bigg]^{1/2}
\end{equation}
for every $u$ Lipschitz continuous function
defined in $B_{r} = B_{r}(x_0)$, with either support contained in $B_{r}(x_0)$ or with null mean value.
\ethm

\noindent
While the following is a particular case of a result contained in \cite{fabio10} (see Theorem 2.9),
which in turn follows from a classical result (see \cite{gut-wheeden3} and \cite{gut-wheeden2}).

\bthm
\label{gut-whee}
Consider $B_r (x_o)$ a ball of $\R^n$,
$\omega \in B_{2,q}^1(K_1)$ for some $q > 2$ and $\omega \in A_{\infty} (K_2, \varsigma)$.
Then there is $\upsigma_1 \in (1,q)$
and a constant $\tilde\upgamma$ $(\tilde\upgamma$ the same constant of Theorem $\ref{chanillo-wheeden})$
such that for every $A \subset B_r (x_o)$, for every Lipschitz continuous function $u$
defined in $B_r (x_o)$, with either support contained in $B_r (x_o)$ or with null mean value in $B_r (x_o)$,
and for every $\kappa \in (1, \upsigma_1]$
$$
\frac{1}{\upsilon (B_r (x_o))} \int_A |u|^{2\kappa} \upsilon \, dx \leqslant \tilde\upgamma^{2} \, r^2 \, 
	\Bigg( \frac{1}{\omega (B_r (x_o))} \int_A  |u|^{2} \omega \, dx \Bigg)^{\kappa -1} 
	\frac{1}{ |B_r (x_o) |} \int_{B_r (x_o)}  |D u|^2 \, dx
$$
where the inequality holds both with $\upsilon = \omega$ and $\upsilon \equiv 1$.
\ethm

\boss
\label{norreni}
About the proof of the previous theorem we want to observe what follows:
if $\upsilon = \omega$ requiring $\omega$ to be doubling (together to $\omega \in B_{2,q}^{1} (K_1)$) is in fact sufficient to get the thesis.
The assumption $\omega \in A_{\infty} (K_2)$ is needed just to get the thesis with $\upsilon \equiv 1$ (and more in general, with $\upsilon \not= \omega$).
In this regard see Remark 2.3 and Theorem 2.9 in \cite{fabio10}.
\eoss

\noindent
From Theorem \ref{chanillo-wheeden} we derive the following.

\bthm
\label{chanillo-wheeden_bis}
Consider $q > 2$, $r \in (0, \bar{R}]$, and suppose \text{\rm (H.2), (H.3), (H.4) - $i \, )$} to hold.
Then there is a constant $\upbeta$ depending $($only$)$ on $n, q, K_1, c_{\uprho}, {\sf L}, \upkappa$ such that for every point $(x_o, t_o) \in \mathcal{Q}_j \cup \Upgamma$
\begin{equation}
\label{disuguaglianzadisobolev}
\bigg[ \frac{1}{\upeta (t) \big( B_r (x_o) \big)} {\int_{B_r (x_o)}} |u(x)|^q \upeta (x, t) dx \bigg]^{1/q} \leqslant \upbeta \, r \,
	\bigg[ \frac{1}{|B_r (x_o)|} {\int_{B_r (x_o)}} |Du(x)|^2 dx \bigg]^{1/2}
\end{equation}
for every $t \in (0,T)$ with $| t - t_o | < \updelta$,
for every $u$ Lipschitz continuous function
defined in $B_{r}(x_o)$, with either support contained in $B_{r}(x_o)$ or with null mean value
and where the inequality holds for every $\upeta$ satisfying $\uprho_j (t) \leqslant \upeta (t) \leqslant \uprho (t)$, or $\upchi_j (t) \leqslant \upeta (t) \leqslant 1$,
for $t \in (0,T)$ with $| t - t_o | < \updelta$, both for $j = 1$ and $j = 2$.
\ethm
\noindent
\dimo
We start taking $\upeta$ such that $\uprho_j (t) \leqslant \upeta (t) \leqslant \uprho (t)$ for $j=1$ or $j=2$.
Taking $\omega = \uprho (\cdot, t)$ in Theorem \ref{chanillo-wheeden} (here we need (H.4) - $i \, )$ and (C.2)) one gets
$$
\bigg[ \frac{1}{\uprho (t) \big( B_r (x_o) \big)} {\int_{B_r (x_o)}} |u(x)|^q \uprho (x, t) dx \bigg]^{1/q} \leqslant \tilde\upgamma \, r \,
	\bigg[ \frac{1}{|B_r (x_o) |} {\int_{B_r (x_o)}} |Du(x)|^2 dx \bigg]^{1/2} .
$$
Then, by (C.8), for $| t - t_o | < \updelta$
\begin{align*}
\frac{1}{\uprho (t) \big( B_r (x_o) \big)} & {\int_{B_r (x_o)}} |u(x)|^q \uprho (x, t) dx
	\geqslant \frac{1}{\uprho (t) \big( B_r (x_o) \big)} {\int_{B_r (x_o)}} |u(x)|^q \upeta (x, t) dx	\geqslant										\\
&	\geqslant \frac{\upkappa}{\sf L} \, \frac{1}{\uprho_j (t) \big( B_r (x_o) \big)} {\int_{B_r (x_o)}} |u(x)|^q \upeta (x, t) dx \geqslant		\\
&	\geqslant \frac{\upkappa}{\sf L} \,  \frac{1}{\upeta (t) \big( B_r (x_o) \big)} {\int_{B_r (x_o)}} |u(x)|^q \upeta (x, t) dx
\end{align*}
by which we conclude taking $\upbeta = \tilde\upgamma \, {\sf L}^{1/q}/\upkappa^{1/q}$. One concludes in the other cases in a similar way,
taking $\omega \equiv 1$ in Theorem \ref{chanillo-wheeden}.
\finedimo

\noindent
From Theorem \ref{gut-whee} we derive the following theorem.

\bthm
\label{gutierrez-wheeden}
Consider $q > 2$ and suppose $\uprho$ satisfies \text{\rm (H.1), (H.2), (H.3), (H.4) - $i \, )$} and \text{\rm (H.4) - $ii \, )$}.
Then there is $\kappa \in (1,q/2)$ and $\upgamma$ depending $($only$)$ on $n, q, K_1, K_2, \varsigma, {\sf L}, \upkappa$ such that
\begin{align*}
& \int_{a}^{b} \frac{1}{\upeta (t) \big( B_{r} \big)} \int_{B_{r}}  |u|^{2\kappa} (x,t) \upeta (x,t) \, dx dt \leqslant 					\\
& \hskip30pt \leqslant \upgamma^{2} \, r^2 \bigg( \max_{a \leqslant t \leqslant b} \frac{1}{\uprho (t) \big( B_{r} \big)} \int_{B_{r}} u^{2}(x,t) \upeta (x,t) \, dx \bigg)^{\kappa-1}
	\frac{1}{|B_{r}|} \int_{a}^{b}\!\! \int_{B_{r}} |D u|^2 (x,t) \, dx dt
\end{align*}
for every $(x_o, t_o) \in \mathcal{Q}_j \cup \Upgamma$,
for every  $B_{r}(x_o) \subset \Omega$ with $r \in (0, \bar{R}]$, for every Lipschitz continuous function $u$
defined in $B_{r} (x_o) \times [t_o - \updelta, t_o + \updelta]$, $u (\cdot, t)$ with either support contained in $B_{r}(x_o)$ or with null mean value in $B_{r}(x_o)$
and for every $[a,b] \subset [t_o - \updelta, t_o + \updelta]$, where the inequality holds for every $j \in \{ 1,2 \}$, every
$\upeta$ satisfying $\uprho_j (t) \leqslant \upeta (t) \leqslant \uprho (t)$ or $\upchi_j (t) \leqslant \upeta (t) \leqslant 1$ for $t \in (0,T)$ with $| t - t_o | \leqslant \updelta$.
\ethm
\noindent
\dimo
Consider $(x_o, t_o) \in \mathcal{Q}_j \cup \Upgamma$.
In Theorem \ref{gut-whee} consider first $\omega = \uprho ( \cdot , t)$
for $t \in [t_o - \updelta, t_o + \updelta]$, and $\upeta$ such that $\uprho_j \leqslant \upeta \leqslant \uprho$ for $j = 1$ or $j = 2$.
Denote $A = A(t) = \{ x \in B_r (x_o) \, | \, \upeta (x,t) > 0 \}$.
Then, for each $t \in [t_o - \updelta, t_o + \updelta]$, we have
\begin{align*}
\frac{1}{\uprho (t) \big( B_{r} \big)} & \int_{B_r (x_o)} |u|^{2\kappa} (x,t) \upeta (x,t) \, dx \leqslant
\frac{1}{\uprho (t) \big( B_{r} \big)} \int_{A (t)} |u|^{2\kappa} (x,t) \uprho (x,t) \, dx \leqslant 							\\
&	\leqslant \tilde\upgamma^{2} \, r^2 \, 
	\Bigg( \frac{1}{\uprho (t) (B_r (x_o))} \int_{A (t)} |u|^{2} (x,t) \uprho (x,t) \, dx \Bigg)^{\kappa -1} 
	\frac{1}{ |B_r (x_o) |} \int_{B_r (x_o)}  |D u|^2 (x,t) \, dx ,
\end{align*}
By (C.3) and (C.5) we get
\begin{align*}
\upeta (t) \big( B_{r} \big) \geqslant \uprho_j (t) \big( B_{r} \big) \geqslant \upkappa \, \uprho (t_o) \big( B_{r} \big)
	\geqslant \frac{\upkappa}{\sf L}  \, \uprho (t) \big( B_{r} \big)
\end{align*}
and then
\begin{align*}
& \frac{1}{\upeta (t) \big( B_{r} \big)} \int_{B_r (x_o)} |u|^{2\kappa} (x,t) \upeta (x,t) \, dx \leqslant 							\\
& \qquad \leqslant \frac{\sf L}{\upkappa}  \, \tilde\upgamma^{2} \, r^2 \, 
	\Bigg( \frac{1}{\uprho (t) (B_r (x_o))} \int_{B_r (x_o)} |u|^{2} (x,t) \upeta (x,t) \, dx \Bigg)^{\kappa -1} 
	\frac{1}{ |B_r (x_o) |} \int_{B_r (x_o)}  |D u|^2 (x,t) \, dx .
\end{align*}
Integrating in time between $t_o - \updelta$ and $t_o + \updelta$ and taking $\upgamma := \tilde\upgamma \sqrt{{\sf L}/\upkappa}$ we conclude. \\
The other case, $\upchi_j \leqslant \upeta \leqslant 1$, is similar and need (C.4).
\finedimo

\noindent
Another result useful in the sequel is the following lemma whose proof can be easily obtained following the proof of
Lemma 2.12 in \cite{fabio10} and using Theorem \ref{chanillo-wheeden_bis}, and for this reason we do not show it.

\begin{lemma}
\label{lemma2.2}
In the same assumptions of Theorem $\ref{chanillo-wheeden_bis}$, consider $k, l \in \R$ with $k < l$, $p \in (1,2]$, $(x_o, t_o) \in \mathcal{Q}_j \cup \Upgamma$,
$j \in \{1, 2 \}$. Then
$$
(l-k) \, \upeta (t) ( \{ v < k \} ) \, \upeta (t) ( \{ v > l \}) \leqslant
	\, 2 \, \upbeta \, r \, \big( \upeta (t) (B_r (x_o)) \big)^2
	\left( \frac{1}{| B_r (x_o) |} \int_{B_r (x_o) \cap \{ k < v < l \}} |Dv|^p \, dx \right)^{\frac{1}{p}}
$$
for every $v$ Lipschitz continuous function defined in the ball $B_r (x_o)$, for every $t \in (0,T)$ such that
$| t - t_o | < \updelta$ and where the inequality holds
for every $\upeta$ satisfying $\uprho_j (t) \leqslant \upeta (t) \leqslant \uprho (t)$, or $\upchi_j (t) \leqslant \upeta (t) \leqslant 1$,
for $t \in (0,T)$ with $| t - t_o | < \updelta$, both for $j = 1$ and $j = 2$, for $t \in (0,T)$ with $| t - t_o | < \updelta$.
\end{lemma}

\noindent
We conclude stating a standard lemma (see, for instance, Lemma 7.1 in \cite{giusti}) needed later.

\begin{lemma}
\label{giusti}
Let $(y_h)_h$ be a sequence of positive real numbers such that
$$
y_{h+1} \leqslant c \, b^h \, y_h^{1+\alpha}
$$
with $c, \alpha > 0$, $b > 1$.
If $y_0 \miu c^{-1/\alpha} b^{-1/\alpha^2}$ then
$$
\lim_{h\to +\infty} y_h = 0 \, .
$$
\end{lemma}

{\color{red}
Le {\sf L} le ho controllate fino a qui
}

\section{DG classes}

We are going now to introduce a De Giorgi class suited to parabolic equations \eqref{equazionegenerale1} and \eqref{equazionegenerale2} of order 2;
before we need to we consider a suitable class of decay functions, the set
\begin{gather*}
\X_c := \big\{ \zeta \in \text{Lip} (\Omega \times (0,T)) \, \big| \,  \zeta (\cdot, t) \in \text{Lip}_c (\Omega) \text{ for each } t \in (0,T) , \zeta \geqslant 0, \zeta_t \geqslant 0 \big\} .
\end{gather*}
By $v_+$ we will denotes the positive part of a function $v$ and by $v_-$ we will denotes the negative part of a function $v$.

\bd
\label{DG1}
Given $\gamma > 0$ we say that a function
$$
u \in \V_{\loc} \quad \text{such that } (0,T) \ni t \mapsto \int_{\omega} u^2 (x,t) \uprho (x,t) \, dx \in C^0_{\loc} (0,T) \qquad \text{for every } \omega \subset \subset \Omega
$$
belongs to $DG(\Omega, T, \uprho, \gamma)$ if $u$ satisfies
the following inequality for every function $\zeta \in \X_{c}$, for every $t_1, t_2 \in (0,T)$ and for every $k \in \R$
\begin{align}
\label{lavagna_bis}
\int_{\Omega} (u - k)_+^2 & \zeta^2 \uprho (x, t_2) \, dx + \int_{t_1}^{t_2} \!\!\!\! \int_{\Omega} |D (u - k)_+ |^2 \zeta^2 dx dt   \leqslant		\nonumber	\\
& \leqslant \int_{\Omega}  (u - k)_+^2 \zeta^2 \uprho (x, t_1) dx \, +  														\nonumber	\\
& \qquad + \frac{\gamma}{2} \Bigg[   \int_{t_1}^{t_2} \!\!\!\! \int_{\Omega} (u - k)_+^2 | D \zeta |^2 \, dx dt +
		\int_{t_1}^{t_2} \!\!\!\! \int_{\Omega}  (u - k)_+^2 \zeta \zeta_t \, \uprho \, dx dt \, +			 									\\
& \qquad \qquad \ + \int_{t_1}^{t_2} \!\!\!\! \int_{\Omega} (u - k)_+^2 \zeta^2 \, dx dt
		+ k^2 \int_{t_1}^{t_2} \!\!\!\! \int_{\{u(t) > k\}} \big( \zeta^2 + |D \zeta|^2 \big) \, dx dt \Bigg]	.							\nonumber
\end{align}
and if $u$ satisfies an analogous 
inequality for every function $\zeta \in \X_{c}$, for every $t_1, t_2 \in (0,T)$, for every $k \in \R$, with $(u - k )_-$ at the place of $(u - k )_+$
and with
$$
k^2 \int_{t_1}^{t_2} \!\!\!\! \int_{\{u(t) < k\}} \big( \zeta^2 + |D \zeta|^2 \big) \, dx dt
$$
in the place of the last term in inequality \eqref{lavagna_bis}.
\ed

\boss
\label{nonepuedopiu}
Notice that by \eqref{lavagna_bis} one in particular has the two following inequalities. The first is
\begin{align*}
\int_{t_1}^{t_2} \!\!\!\! \int_{\Omega} |D (u - k)_+ |^2 & \, \zeta^2 dx dt \leqslant \int_{\Omega}  (u - k)_+^2 \zeta^2 \uprho (x, t_1) \, dx \, +  		\nonumber	\\
& + \frac{\gamma}{2} \Bigg[   \int_{s_1}^{s_2} \!\!\!\! \int_{\Omega} (u - k)_+^2 | D \zeta |^2 \, dx dt +
		\int_{t_1}^{t_2} \!\!\!\! \int_{\Omega}  (u - k)_+^2 \zeta \zeta_t \, \uprho \, dx dt \, +			 										\\
& \qquad \ + \int_{t_1}^{t_2} \!\!\!\! \int_{\Omega} (u - k)_+^2 \zeta^2 \, dx dt
		+ k^2 \int_{t_1}^{t_2} \!\!\!\! \int_{\{u(t) > k\}} \big( \zeta^2 + |D \zeta|^2 \big) \, dx dt \Bigg]	.								\nonumber
\end{align*}
For the second one consider $t \in [t_1, t_2]$ and write \eqref{lavagna_bis} in the interval $[t_1, t]$: then, taking the maximum for $t \in [t_1, t_2]$ one has
\begin{align*}
\max_{t \in [t_1, t_2]} \int_{\Omega} (u - k)_+^2 & \zeta^2 \uprho (x, t) \, dx   \leqslant \int_{\Omega}  (u - k)_+^2 \zeta^2 \uprho (x, t_1) dx \, +  			\\
& + \frac{\gamma}{2} \Bigg[   \int_{t_1}^{t_2} \!\!\!\! \int_{\Omega} (u - k)_+^2 | D \zeta |^2 \, dx dt +
		\int_{t_1}^{t_2} \!\!\!\! \int_{\Omega}  (u - k)_+^2 \zeta \zeta_t \, \uprho \, dx dt \, +			 									\\
& \qquad \ + \int_{t_1}^{t_2} \!\!\!\! \int_{\Omega} (u - k)_+^2 \zeta^2 \, dx dt
		+ k^2 \int_{t_1}^{t_2} \!\!\!\! \int_{\{u(t) > k\}} \big( \zeta^2 + |D \zeta|^2 \big) \, dx dt \Bigg]	.		
\end{align*}
Summing the two inequalities  one gets that for every $\zeta \in \X_c$ for which $\zeta (\cdot, t_1) \equiv 0$
\begin{align*}
\max_{t \in [t_1, t_2]}  \int_{\Omega} (u - k)_+^2 & \zeta^2 \uprho (x, t) \, dx + \int_{t_1}^{t_2} \!\!\!\! \int_{\Omega} |D (u - k)_+ |^2 \zeta^2 dx dt   \leqslant	\\
& \leqslant \gamma \Bigg[   \int_{t_1}^{t_2} \!\!\!\! \int_{\Omega} (u - k)_+^2 | D \zeta |^2 \, dx dt +
		\int_{t_1}^{t_2} \!\!\!\! \int_{\Omega}  (u - k)_+^2 \zeta \zeta_t \, \uprho \, dx dt \, +			 									\\
& \qquad \qquad \ + \int_{t_1}^{t_2} \!\!\!\! \int_{\Omega} (u - k)_+^2 \zeta^2 \, dx dt
		+ k^2 \int_{t_1}^{t_2} \!\!\!\! \int_{\{u(t) > k\}} \big( \zeta^2 + |D \zeta|^2 \big) \, dx dt \Bigg]	.
\end{align*}
\eoss

\boss
\label{commenti}
{\bf Comments on Definition \ref{DG1}} \\ [0.1em]
- The energy inequality \eqref{lavagna_bis} can be derived as done in \cite{fabio13} starting from both equations
\eqref{equazionegenerale1} and \eqref{equazionegenerale2}.
Notice that in \cite{fabio13} only $k \geqslant 0$ is considered to derive \eqref{lavagna_bis} (and only $k \leqslant 0$ to derive the analogous one for $(u - k)_-$):
indeed these restrictions are due to the fact that in \cite{fabio13} a wider class was considered and can be removed in our case. \\ [0.1em]
- A solution of \eqref{equazionegenerale1} belongs to $\{ v \in \V \, | \, \uprho v' \in \V' \}$ and
a solution of \eqref{equazionegenerale2} belongs to $\{ v \in \V \, | \, (\uprho v)' \in \V' \}$ but in fact these two spaces,
under assumption (H.1), turn out to be the same space  (see Remark 3.1 in \cite{fabio9.1}). \\ [0.1em]
- A function belonging to $DG(\Omega, T, \uprho, \gamma)$ is locally bounded (see \cite{fabio21}).
\eoss

\boss
\label{holder}
The solutions of equation \eqref{equazionegenerale1} with $C \equiv 0$ satisfy
\begin{align}
\int_{\Omega} (u & - k)_+^2 \zeta^2 \uprho (x, t_2) \, dx + \int_{t_1}^{t_2} \!\!\!\! \int_{\Omega} |D (u - k)_+ |^2 \zeta^2 dx dt   \leqslant		\nonumber	\\
\label{lavagna_tris}
& \leqslant \int_{\Omega}  (u - k)_+^2 \zeta^2 \uprho (x, t_1) dx +
		\frac{\gamma}{2} \Bigg[   \int_{t_1}^{t_2} \!\!\!\! \int_{\Omega} (u - k)_+^2 | D \zeta |^2 \, dx dt \, +									\\
& \qquad + \int_{t_1}^{t_2} \!\!\!\! \int_{\Omega}  (u - k)_+^2 \zeta \zeta_t \, \uprho \, dx dt
		+ \int_{t_1}^{t_2} \!\!\!\! \int_{\Omega} (u - k)_+^2 \zeta^2 \, dx dt \Bigg]											\nonumber
\end{align}
and the analogous one for $(u - k)_-$.
Notice that if $u$ satisfies these inequalities also $u + c$ satisfies them, where $c$ is an arbitrary constant
(in particular this holds for solutions of \eqref{equazionegenerale1}, but not of \eqref{equazionegenerale2}, with $C \equiv 0$). \\
This is the key point to use the argument due to Moser, by which one can derive the local H\"older continuity from the Harnack inequality
for functions satisfying \eqref{lavagna_tris} and the analogous estimate with $(u - k)_-$.
\eoss

\noindent
Given $E \subset \Omega \times (0,T)$, for each $t \in (0,T)$ we define
\begin{align}
\label{wlemedie}
E(t) := E \cap \big( \Omega \times \{ t \} \big) .
\end{align}
Then for $r > 0$ such that $r^2 \in (0, \updelta)$
such that $B_r (x_o) \times (t_o - \updelta^2, t_o + \updelta^2) \subset \Omega \times (0,T)$ we define, for $j \in \{ 1, 2 \}$, the quantities
\begin{equation}
\label{funzioneacca}
\begin{array}{c}
{\sf h}_{\ast}^j (x_o, t_o, r) := 0 \qquad \text{ if } \big( B_r (x_o) \times (t_o - r^2 , t_o) \big) \cap \mathcal{Q}_j  = \emptyset	,	\\	[0.7em]
{\sf h}_{\ast}^j (x_o, t_o, r) := 
{\displaystyle \frac{\uprho \big( \big( B_r (x_o) \times (t_o - r^2 , t_o) \big) \cap \mathcal{Q}_j \big)}{\big| \big( B_r (x_o) \times (t_o - r^2 , t_o) \big) \cap \mathcal{Q}_j \big|}}		
	\qquad \text{otherwise,} 																			\\	[1.5em]
{\sf h}^{\ast}_j (x_o, t_o, r) := 0 \qquad \text{ if } \big( B_r (x_o) \times (t_o, t_o + r^2) \big) \cap \mathcal{Q}_j  = \emptyset ,	\\	[0.7em]
{\sf h}^{\ast}_j (x_o, t_o, r) := 
{\displaystyle \frac{\uprho \big( \big( B_r (x_o) \times (t_o, t_o + r^2) \big) \cap \mathcal{Q}_j \big)}{\big| \big( B_r (x_o) \times (t_o, t_o + r^2) \big) \cap \mathcal{Q}_j \big| }}
	\qquad \text{otherwise,} 																			\\	[1.5em]
{\sf h}_{\ast} (x_o, t_o, r) = \max\{ {\sf h}_{\ast}^1 (x_o, t_o, r), {\sf h}_{\ast}^2 (x_o, t_o, r) \} ,		\\	[1em]
{\sf h}^{\ast} (x_o, t_o, r) = \max\{ {\sf h}^{\ast}_1 (x_o, t_o, r), {\sf h}^{\ast}_2 (x_o, t_o, r) \} .
\end{array}
\end{equation}
Notice that
\begin{equation}
\label{accaasterisco}
{\sf h}^{\ast}_j \big( x_o, t_o - r^2 , r \big) = {\sf h}_{\ast}^j ( x_o, t_o, r ) ,			\qquad \quad
{\sf h}^{\ast}_j ( x_o, t_o, r ) = {\sf h}_{\ast}^j \big( x_o, t_o + r^2 , r \big) .
\end{equation}
Notice that for $0 < r < \tilde{r}$ with $\tilde{r}^2 < \updelta$, by (H.2), one has
\begin{equation}
\label{valsusa}
\begin{array}{c}
{\displaystyle 
{\sf h}_{\ast}^j (x_o, t_o, r) \leqslant \frac{1}{2 \upkappa} \left( \frac{\tilde{r}}{r} \right)^{n+2} {\sf h}_{\ast}^j (x_o, t_o, \tilde{r}) ,		\qquad
{\sf h}^{\ast}_j (x_o, t_o, r) \leqslant \frac{1}{2 \upkappa} \left( \frac{\tilde{r}}{r} \right)^{n+2} {\sf h}^{\ast}_j (x_o, t_o, \tilde{r})  .				}
\end{array}
\end{equation}
Once fixed $(x_o,t_o) \in \Omega \times (0,T)$, for a generic $R > 0$ and $R^2 \, {\sf h}_{\ast} (x_o, t_o, R) < \updelta$ and
$R^2 \, {\sf h}^{\ast} (x_o, t_o, R) < \updelta$ we will denote
\begin{equation}
\label{querre+}
\begin{array}{l}
Q_{R} (x_o, t_o) := B_R (x_o) \times [ t_o - R^2 \, {\sf h}_{\ast} (x_o, t_o, R) , t_o ]	,	\\	[0.5em]
Q^{R} (x_o, t_o) := B_R (x_o) \times [  t_o, t_o + R^2 \, {\sf h}^{\ast} (x_o, t_o, R) ] .
\end{array}
\end{equation}
We also introduce
\begin{equation}
\label{querreteta_+}
\begin{array}{c}
B_r^j (x_o; t) = \big( B_r ( x_o) \times \{ t \} \big) \cap \mathcal{Q}_j	,									\\	[0.5em]
Q_{r, \theta}^j (x_o, t_o) = \bigcup_{t \in [t_o - \theta \, r^2 \, {\sf h}_{\ast}^j ( x_o, t_o, R ), t_o]} B_r^j ( x_o ; t) ,		\\	[0.5em]
Q^{r, \theta}_j (x_o, t_o) = \bigcup_{t \in [t_o , t_o + \theta \, r^2 \, {\sf h}^{\ast}_j ( x_o, t_o, R )]} B_r^j ( x_o ; t) ,		\\	[0.5em]
Q_{r, \theta} (x_o, t_o) := B_r (x_o) \times [t_o - \theta \, r^2 \, {\sf h}_{\ast} ( x_o, t_o, R ), t_o ] ,					\\	[0.5em]
Q^{r, \theta} (x_o, t_o) := B_r (x_o) \times [t_o, t_o + \theta \, r^2 \, {\sf h}^{\ast} ( x_o, t_o, R )] .
\end{array}
\end{equation}
Now consider $B_r^j (x_o; t)$, $Q_{r, \theta}^j (x_o, t_o)$ and $Q^{r, \theta}_j (x_o, t_o)$ for $j = 1,2$.
Then one can define $\varsigma := (\varsigma_1, \varsigma_2) \in \R^2$ and the sets
\begin{equation}
\label{svaroschi}
\begin{array}{lc}
\big( B_r^j (x_o; t) \big)^{\varsigma_1} := B_r^j (x_o; t)  \ \cup \ 
		\big( \big\{ y \in B_r (x_o) \setminus B_r^j (x_o; t) \, \big| \, \text{dist} \, \big( y, B_r^j (x_o; t) \big) \leqslant \varsigma_1 \big\} \times \{ t \} \big) \, ,			\\	[0.5em]
\big( Q_{r, \theta}^j (x_o, t_o) \big)^{\varsigma} := 
	\big( \bigcup_{t \in [t_o - \theta \, r^2 \, {\sf h}_{\ast}^j ( x_o, t_o, R ) - \varsigma_2, t_o]} \big( B_r^j (x_o; t) \big)^{\varsigma_1} \big) ,						\\	[0.5em]
\big( Q^{r, \theta}_j (x_o, t_o) \big)^{\varsigma} :=
	\big( \bigcup_{t \in [t_o, t_o + \theta \, r^2 \, {\sf h}^{\ast}_j ( x_o, t_o, R ) + \varsigma_2]} \big( B_r^j (x_o; t) \big)^{\varsigma_1} \big) .
\end{array}
\end{equation}


\noindent
To lighten the notation we will omit to write $(x_o, t_o)$ if not strictly necessary and simply write
$$
Q_{R}, \quad Q^{R}, \quad Q_{r, \theta}^j , \quad Q^{r, \theta}_j
$$
if it is clear what point we are referring to.

\section{Expansion of positivity}

Functions belonging to the class defined in Definition \ref{DG1} are locally bounded. 
This can be shown adapting the proof for the standard evolution equations
(i.e. when $\uprho \equiv 1$), but we refer also to \cite{fabio21}, where the case when $\uprho \geqslant 0$ is treated, which clearly holds
also when $\uprho > 0$.
Therefore from now on we will assume that functions in the class $DG(\Omega, T, \uprho, \gamma)$ are locally bounded. \\ [0.3em]
%
{\bf Preliminary computations - } Here we show some estimates useful 
to prove Proposition \ref{carlettobello} and Proposition \ref{anninabella}. \\ [0.3em]
In the following computations we will need directly or indirectly, i.e. using some of their consequences, all assumptions (H.1)-(H.4).  \\ [0.3em]
Consider $j \in \{1, 2 \}$, $(\bar{x}, \bar{t}) \in \mathcal{Q}_j \cup \Upgamma$, $R \in (0 , \bar{R}]$, let ${\sf h}_{\ast}^j = {\sf h}_{\ast}^j ( \bar{x}, \bar{t}, R)$.
Suppose that $Q_{R} = Q_{R} (\bar{x}, \bar{t}) \subset B_R (\bar{x}) \times [\bar{t} - \updelta, \bar{t}] \subset \Omega \times (0,T)$.
Consider $r, \tilde{r}, \theta, \tilde{\theta}, k, \tilde{k}, \varsigma = (\varsigma_{1}, \varsigma_{2}), \tilde\varsigma = (\tilde\varsigma_{1}, \tilde\varsigma_{2})$ satisfying
\begin{gather*}
0 < r < \tilde{r} \leqslant  R \, , \qquad 0 < \theta \leqslant \tilde\theta \qquad 
	\text{s.t. }	Q_{\tilde{r}, \tilde\theta} (\bar{x}, \bar{t}) \subset B_R (\bar{x}) \times [\bar{t} - \updelta, \bar{t}] \subset \Omega \times (0,T) \, ,		\\
0 \leqslant k < \tilde{k}, \qquad 
\varsigma_1 \geqslant 0 , \quad \varsigma_2 \geqslant 0, \qquad
\tilde\varsigma_1 = \varsigma_1 + \tilde{r} - r ,
\quad \tilde\varsigma_2 = \varsigma_2 + (\tilde\theta \tilde{r}^2 - \theta r^2) {\sf h}_{\ast}^j \, .
\end{gather*}
Notice that the following estimate holds:
\begin{align}
\label{pippo}
\tilde\theta \tilde{r}^2 - \theta r^2 \geqslant \tilde\theta (\tilde{r}^2 - r^2) \geqslant \tilde\theta (\tilde{r} - r)^2 \, .
\end{align}
Consider now a function $\zeta$ satisfying
\begin{gather*}
\zeta \equiv 1 \hskip10pt \text{in } \big( Q_{r, \theta}^{j} (\bar{x}, \bar{t}) \big)^{\varsigma}	\, , \qquad 
	\zeta \equiv 0 \hskip10pt \text{ outside of } \big( Q_{\tilde{r}, \tilde\theta}^{j} (\bar{x}, \bar{t}) \big)^{\tilde\varsigma}
	\quad \text{for } t \leqslant \bar{t}	\, , 																								\\
0 \leqslant \zeta \leqslant 1 \, , \hskip15pt {\displaystyle |D \zeta| \leqslant \frac{1}{\tilde{r} - r} } \, ,														\\
0 \leqslant  \zeta_t \leqslant \frac{1}{(\tilde\theta \tilde{r}^2 - \theta r^2) {\sf h}^{\ast}_j} \leqslant \frac{1}{\tilde\theta (\tilde{r} - r)^2 {\sf h}_{\ast}^j}
	\quad (\text{by } \eqref{pippo})
\end{gather*}
and finally define
\begin{gather*}
A (l) = \big\{ (x,t) \in \big( Q_{r, \theta}^{j} (\bar{x}, \bar{t}) \big)^{\varsigma} \, | \, u (x,t) > l \big\}   ,						\\
\tilde{A} (l) = \big\{ (x,t) \in \big( Q_{\tilde{r}, \tilde\theta}^{j} (\bar{x}, \bar{t}) \big)^{\tilde\varsigma} \, | \, u (x,t) > l \big\}   .
\end{gather*}
In the following $\upeta$ denotes a function defined in $Q_{R} (\bar{x}, \bar{t})$ 
\begin{align}
\label{eta1}
\text{either as} \hskip80pt & \upeta (x,t) := 
\left\{
\begin{array}{cl}
\uprho		&	\text{if } (x,t) \in \big( Q_{\tilde{r}, \tilde\theta}^{j} (\bar{x}, \bar{t}) \big)^{\tilde\varsigma}	,		\\	[0.3em]
	0		&	\text{if } (x,t) \in Q_{R} (\bar{x}, \bar{t}) \setminus \big( Q_{\tilde{r}, \tilde\theta}^{j} (\bar{x}, \bar{t}) \big)^{\tilde\varsigma}
\end{array}
\right.																						\\
																\ 				\nonumber	\\
\label{eta2}
\text{or as} \hskip80pt & \upeta (x,t) := 
\left\{
\begin{array}{cl}
	1		&	\text{if } (x,t) \in \big( Q_{\tilde{r}, \tilde\theta}^{j} (\bar{x}, \bar{t}) \big)^{\tilde\varsigma}	,		\\	[0.3em]
	0		&	\text{if } (x,t) \in Q_{R} (\bar{x}, \bar{t}) \setminus \big( Q_{\tilde{r}, \tilde\theta}^{j} (\bar{x}, \bar{t}) \big)^{\tilde\varsigma} .
\end{array}
\right.
\end{align}
Notice that for $t \in (0,T)$ with $| t - \bar{t} | \leqslant \updelta$ and $j = 1$ or $j = 2$ $\upeta$ satisfies
$\uprho_j (t) \leqslant \upeta (t) \leqslant \uprho (t)$ in the first case and  $\upchi_j (t) \leqslant \upeta (t) \leqslant 1$ in the second. \\
Observe that by (in the order) (C.1), (C.8), (C.5) and (C.6) one has that for every $(a,b) \subset (\bar{t} - \updelta, \bar{t} + \updelta)$ and for every $\varrho \leqslant R$ one has
(for some $s \in [a,b]$)
\begin{equation}
\begin{array}{l}
\label{succofresco_0}
{\displaystyle \sup_{t \in (a,b)} \upeta (t) \big( B_{\varrho} (\bar{x}) \big) \leqslant 
	\sup_{t \in (a,b)} \uprho (t) \big( B_{\varrho} (\bar{x}) \big) = \uprho (s) \big( B_{\varrho} (\bar{x}) \big) \leqslant
	\frac{\sf L}{\upkappa} \, \uprho_j (s) \big( B_{\varrho} (\bar{x}) \big) \leqslant									}	\\	[1.2em]
\hskip40pt 
{\displaystyle \leqslant \frac{(1 + \upkappa) {\sf L}}{\upkappa^2} \inf_{t \in (a,b)} \uprho_j (t) \big( B_{\varrho}(\bar{x}) \big)
	\leqslant \frac{(1 + \upkappa) {\sf L}}{\upkappa^2} \inf_{t \in (a,b)} \upeta (t) \big( B_{\varrho} (\bar{x}) \big)			}
\end{array}
\end{equation}
and similarly
\begin{align}
\label{succofresco}
\sup_{t \in (a,b)} & \frac{1}{\upeta (t) \big( B_{\varrho} (\bar{x}) \big)}
	\leqslant \frac{(1 + \upkappa) {\sf L}}{\upkappa^2} \inf_{t \in (a,b)} \frac{1}{\upeta (t) \big( B_{\varrho} (\bar{x}) \big)} \, .
\end{align}
Now consider ($\kappa$ is the value in Theorem \ref{gutierrez-wheeden})
\begin{equation}
\label{dentista}
\begin{array}{l}
{\displaystyle
\iint_{A (k)} \frac{(u - \tilde{k})_+^2}{\upeta (t) \big( B_{\tilde{r}} (\bar{x}) \big)} \, \upeta \, dx dt \leqslant	
	\iint_{\tilde{A} (k)} \frac{(u - \tilde{k})_+^2 \zeta^2}{\upeta (t) \big( B_{\tilde{r}} (\bar{x}) \big)} \, \upeta \, dx dt \leqslant		}			\\	[1.4em]
\qquad {\displaystyle
\leqslant \left( \iint_{\tilde{A} (k)} \frac{1}{\upeta (t) \big( B_{\tilde{r}} (\bar{x}) \big)} \, \upeta \, dx dt \right)^{\frac{\kappa - 1}{\kappa}}
\left( \iint_{\big( Q_{\tilde{r}, \tilde\theta}^{j} (\bar{x}, \bar{t}) \big)^{\tilde\varsigma}} 
	\frac{(u - \tilde{k})_+^{2 \kappa} \zeta^{2 \kappa}}{\upeta (t) \big( B_{\tilde{r}} (\bar{x}) \big)}
	\, \upeta \, dx dt \right)^{\frac{1}{\kappa}}	.			
}
\end{array}
\end{equation}
The second factor on the right hand side can be estimated, using in the following order Theorem \ref{gutierrez-wheeden}, (C.1), the inequality
$\alpha^{\frac{\kappa-1}{\kappa}} \beta^{\frac{1}{\kappa}} \leqslant \alpha + \beta$ ($\alpha, \beta > 0$),
the energy inequality defining De Giorgi classes (see in particular Remark \ref{nonepuedopiu}), as
\begin{align*}
& \left( \iint_{\big( Q_{\tilde{r}, \tilde\theta}^{j} (\bar{x}, \bar{t}) \big)^{\tilde\varsigma}} 
	\frac{1}{\upeta (t) \big( B_{\tilde{r}} (\bar{x}) \big)} \, (u - \tilde{k})_+^{2 \kappa} \zeta^{2 \kappa} \upeta \, dx dt \right)^{\frac{1}{\kappa}}	\leqslant					\\
& \qquad \leqslant \upgamma^{\frac{2}{\kappa}} \, {\tilde{r}}^{\frac{2}{\kappa}}
	\bigg( \max_{t \in [\bar{t} - \tilde\theta \, \tilde{r}^2 {\sf h}_{\ast}^j , \bar{t}]} \frac{1}{\uprho (t) \big( B_{\tilde{r}} \big)} 
	\int_{B_{\tilde{r}}} (u - \tilde{k})_+^{2} \zeta^{2} \uprho (x,t) \, dx \bigg)^{\frac{\kappa-1}{\kappa}} \cdot													\\
& \hskip80pt \cdot \left( \frac{1}{|B_{\tilde{r}}|} 
	\iint_{\big( Q_{\tilde{r}, \tilde\theta}^{j} (\bar{x}, \bar{t}) \big)^{\tilde\varsigma}}  |D ((u - \tilde{k})_+ \zeta)|^2 (x,t) \, dx dt \right)^{\frac{1}{\kappa}}	\leqslant			\\
& \qquad \leqslant \upgamma^{\frac{2}{\kappa}} \, {\tilde{r}}^{\frac{2}{\kappa}}
	\left( \max_{\bar{t} - \tilde\theta \, \tilde{r}^2 {\sf h}_{\ast}^j  \leqslant t \leqslant \bar{t}} \frac{1}{\uprho (t) \big( B_{\tilde{r}} \big)} \right)^{\frac{\kappa-1}{\kappa}}
	\left( \frac{1}{|B_{\tilde{r}}|}  \right)^{\frac{1}{\kappa}} \cdot																				\\
& \hskip40pt \cdot \Bigg[ \max_{t \in [\bar{t} - \tilde\theta \, \tilde{r}^2 {\sf h}_{\ast}^j , \bar{t}]} \int_{B_{\tilde{r}}} (u - \tilde{k})_+^{2} \zeta^2 \uprho (x,t) \, dx +
	2 \iint_{\big( Q_{\tilde{r}, \tilde\theta}^{j} (\bar{x}, \bar{t}) \big)^{\tilde\varsigma}} |D (u - \tilde{k})_+|^2 \zeta^2 \, dx dt +										\\
& \hskip80pt + 2 \iint_{\big( Q_{\tilde{r}, \tilde\theta}^{j} (\bar{x}, \bar{t}) \big)^{\tilde\varsigma}} |D \zeta |^2 (u - \tilde{k})_+|^2 \, dx dt \Bigg] 							\\
& \qquad \leqslant \upgamma^{\frac{2}{\kappa}} \, \tilde{r}^{\frac{2}{\kappa}}
	\left( \max_{\bar{t} -\tilde\theta \, \tilde{r}^2 {\sf h}_{\ast}^j  \leqslant t \leqslant \bar{t}} \frac{1}{\uprho (t) \big( B_{\tilde{r}} \big)} \right)^{\frac{\kappa-1}{\kappa}}
	\left( \frac{1}{|B_{\tilde{r}}|}  \right)^{\frac{1}{\kappa}} \cdot																					\\
& \hskip40pt  \cdot (2 \gamma + 2) \Bigg[ \iint_{\big( Q_{\tilde{r}, \tilde\theta}^{j} (\bar{x}, \bar{t}) \big)^{\tilde\varsigma}}
	(u - \tilde{k})_+^2 \Big[ | D \zeta |^2 + \zeta \zeta_t \, \uprho + \zeta^2 \Big] \, dx dt \, +																\\
& \hskip100pt + \tilde{k}^2 \iint_{\big( Q_{\tilde{r}, \tilde\theta}^{j} (\bar{x}, \bar{t}) \big)^{\tilde\varsigma} \cap \{u > \tilde{k}\}} \big( \zeta^2 + |D \zeta|^2 \big) \, dx dt \Bigg]	.
\end{align*}
Summing up, by \eqref{succofresco_0}, \eqref{dentista}, the last estimates, (C.5) and dividing by $\tilde\theta \, \tilde{r}^2 {\sf h}_{\ast}^j$ we get
\begin{align}
& \frac{1}{\upeta \big( \big( Q_{\tilde{r}, \tilde\theta}^{j} (\bar{x}, \bar{t}) \big)^{\tilde\varsigma} \big)} \iint_{A (\tilde{k})}(u - \tilde{k})_+^2  \upeta\, dx dt \leqslant
	\upgamma^{\frac{2}{\kappa}} \, (2 \gamma + 2) \, \frac{(1 + \upkappa) {\sf L}}{\upkappa^2}
	\left( \frac{1 + \upkappa}{\upkappa} \right)^{\frac{\kappa - 1}{\kappa}}  \cdot																	\nonumber	\\
& \qquad \cdot {\tilde{r}}^{\frac{2}{\kappa}} \big( \tilde\theta \, {\tilde{r}}^2 {\sf h}_{\ast}^j \big)^{\frac{\kappa - 1}{\kappa}}
	\left( \frac{\upeta \big( \tilde{A} (\tilde{k}) \big)}{\upeta \big( \big( Q_{\tilde{r}, \tilde\theta}^{j} (\bar{x}, \bar{t}) \big)^{\tilde\varsigma} \big)} \right)^{\frac{\kappa - 1}{\kappa}}
	\left( \frac{1}{\uprho_j \big( Q_{\tilde{r}, \tilde\theta}^{j} (\bar{x}, \bar{t}) \big)} \right)^{\frac{\kappa - 1}{\kappa}}
	\left( \frac{1}{\big| Q_{\tilde{r}, \tilde\theta}^{j} (\bar{x}, \bar{t})\big| } \right)^{\frac{1}{\kappa}} \cdot													\nonumber	\\ 
\label{partitella}
& \qquad \cdot \Bigg[ \left( \frac{1}{(\tilde{r} - r)^2} + 1 \right)
	\iint_{\big( Q_{\tilde{r}, \tilde\theta}^{j} (\bar{x}, \bar{t}) \big)^{\tilde\varsigma}} (u - \tilde{k})_+^2 \, dx dt \, +															\\
& \hskip50pt  + \frac{1}{\tilde\theta (\tilde{r} - r)^2 {\sf h}_{\ast}^j} \iint_{\big( Q_{\tilde{r}, \tilde\theta}^{j} (\bar{x}, \bar{t}) \big)^{\tilde\varsigma}} (u - \tilde{k})_+^2 \uprho \, dx dt +
	\tilde{k}^2 \left( \frac{1}{(\tilde{r} - r)^2} + 1 \right) \big| \tilde{A} (\tilde{k}) \big| \Bigg]	.															\nonumber
\end{align}
Now we need to estimate
$\Big( \uprho_j \big( Q_{\tilde{r}, \tilde\theta}^{j} (\bar{x}, \bar{t}) \big) \Big)^{-\frac{\kappa - 1}{\kappa}}
\Big( \big| Q_{\tilde{r}, \tilde\theta}^{j} (\bar{x}, \bar{t})\big| \Big)^{-\frac{1}{\kappa}}$
in a suitable way.
Notice that for $[a,b], [A,B] \subset (\bar{t} - \updelta, \bar{t} + \updelta)$, $\bar{t} \in [a, b]$, $\bar{t} \in [A, B]$, $\varrho_1 , \varrho_2 > 0$ and by (C.5) and (C.7) one has
\begin{align*}
& \int_A^B \uprho_j (t) (B_{\varrho_2} (\bar{x})) \, dt \leqslant (B-A) \sup_{t \in (A, B)} \uprho_j (t) (B_{\varrho_2} (\bar{x})) \leqslant 
	\frac{1 + \upkappa}{\upkappa} \, (B - A) \inf_{t \in (A, B)} \uprho_j (t) (B_{\varrho_2} (\bar{x})) \leqslant										\\	
& \hskip30pt \leqslant \frac{1 + \upkappa}{\upkappa} \, (B - A) \, c_{\uprho_j} (\varrho_2 / \varrho_1) \, \uprho_j ( \bar{t}) (B_{\varrho_1} (\bar{x})) \leqslant
	\frac{1 + \upkappa}{\upkappa} \, (B - A) \, c_{\uprho_j} (\varrho_2 / \varrho_1) \, \sup_{t \in (a, b)} \uprho_j (t) (B_{\varrho_1} (\bar{x}) \leqslant			\\
& \hskip30pt \leqslant \left(\frac{1 + \upkappa}{\upkappa} \right)^2 \frac{B - A}{b-a} \, c_{\uprho_j} (\varrho_2 / \varrho_1) \, 
	(b - a) \inf_{t \in (a, b)} \uprho_j (t) (B_{\varrho_1} (\bar{x}) \leqslant																	\\
& \hskip30pt \leqslant \left(\frac{1 + \upkappa}{\upkappa} \right)^2 \frac{B - A}{b-a} \, c_{\uprho_j} (\varrho_2 / \varrho_1) \int_a^b \uprho_j (t) (B_{\varrho_1} (\bar{x})) \, dt
\end{align*}
and an analogous estimate for $\upchi_j$
\begin{align*}
& \int_A^B \upchi_j (t) (B_{\varrho_2} (\bar{x})) \, dt \leqslant 
	\left(\frac{1 + \upkappa}{\upkappa} \right)^2 \frac{B - A}{b-a} \, c_{\upchi_j} (\varrho_2 / \varrho_1) \int_a^b \upchi_j (t) (B_{\varrho_1} (\bar{x})) \, dt \, .
\end{align*}
In a similar way, if $\uprho_j \leqslant \upeta \leqslant \uprho$, one can derive
\begin{align}
\label{piscinaieri}
& \int_A^B \upeta (t) (B_{\varrho_2} (\bar{x})) \, dt \leqslant 
	\frac{(1 + \upkappa ) {\sf L}^2}{\upkappa^2} \, \frac{B - A}{b-a} \, c_{\uprho} (\varrho_2 / \varrho_1) \int_a^b \upeta (t) (B_{\varrho_1} (\bar{x})) \, dt ,
\end{align}
while if $\upchi_j \leqslant \upeta \leqslant 1$
\begin{align}
\label{piscinaoggi}
& \int_A^B \upeta (t) (B_{\varrho_2} (\bar{x})) \, dt \leqslant 
	\frac{\sf L}{\upkappa} \, \frac{B - A}{b-a} \, (\varrho_2 / \varrho_1)^n \int_a^b \upeta (t) (B_{\varrho_1} (\bar{x})) \, dt ,
\end{align}
Therefore one gets
\begin{align}
& \uprho_j \big( \big( B_R (\bar{x}) \times (\bar{t} - R^2 , \bar{t}) \big) \cap \mathcal{Q}_j \big) \big) \leqslant
	\left(\frac{1 + \upkappa}{\upkappa} \right)^2 \frac{R^2}{\tilde\theta \, \tilde{r}^2 {\sf h}_{\ast}^j} \, c_{\uprho_j} (R / \tilde{r}) \, \uprho_j \big( Q_{\tilde{r}, \tilde\theta}^{j} (\bar{x}, \bar{t}) \big) ,	\nonumber	\\
\label{euno}
& \uprho_j \big( Q_{\tilde{r}, \tilde\theta}^{j} (\bar{x}, \bar{t}) \big) \leqslant
	\left(\frac{1 + \upkappa}{\upkappa} \right)^2 \frac{\tilde\theta \, \tilde{r}^2 {\sf h}_{\ast}^j}{R^2} \, \uprho_j \big( \big( B_R (\bar{x}) \times (\bar{t} - R^2 , \bar{t}) \big) \cap \mathcal{Q}_j \big) \big) ,		\\
& \upchi_j \big( \big( B_R (\bar{x}) \times (\bar{t} - R^2 , \bar{t}) \big) \cap \mathcal{Q}_j \big) \big) \leqslant
	\left(\frac{1 + \upkappa}{\upkappa} \right)^2 \frac{R^2}{\tilde\theta \, \tilde{r}^2 {\sf h}_{\ast}^j} \, c_{\upchi_j} (R / \tilde{r}) \, \upchi_j \big( Q_{\tilde{r}, \tilde\theta}^{j} (\bar{x}, \bar{t}) \big) .	\nonumber
\end{align}
From these we derive
\begin{align}
\label{edue}
\frac{\uprho_j \big( Q_{\tilde{r}, \tilde\theta}^{j} (\bar{x}, \bar{t}) \big)}{\big| Q_{\tilde{r}, \tilde\theta}^{j} (\bar{x}, \bar{t})\big| } \leqslant
	\left(\frac{1 + \upkappa}{\upkappa} \right)^4 \, c_{\upchi_j} (R / \tilde{r}) \, {\sf h}_{\ast}^j \, .
\end{align}
Moreover by (C.6) and (C.8) we also get, for $\upeta$ satisfying $\uprho_j \leqslant \upeta \leqslant \uprho$ and for every $t, s \in [\bar{t} - \updelta, \bar{t} + \updelta]$,
$$
\upeta (t) \big( B_{\tilde{r}} (\bar{x}) \big) \leqslant \uprho (t) \big( B_{\tilde{r}} (\bar{x}) \big) \leqslant {\sf L} \, \uprho (\bar{t}) \big( B_{\tilde{r}} (\bar{x}) \big) \leqslant
	\frac{\sf L}{\upkappa} \, \uprho_j (s) \big( B_{\tilde{r}} (\bar{x}) \big) 
$$
while for $\upeta$ satisfying $\upchi_j \leqslant \upeta \leqslant 1$, by (C.4) and for every $t, s \in [\bar{t} - \updelta, \bar{t} + \updelta]$ we get
$$
\upeta (t) \big( B_{\tilde{r}} (\bar{x}) \big) \leqslant \big| B_{\tilde{r}} (\bar{x}) \big| \leqslant \frac{1}{\upkappa} \, \upchi_j (t) \big( B_{\tilde{r}} (\bar{x}) \big) .
$$
In particular we derive the three useful estimates:
\begin{equation}
\label{etre}
\begin{array}{l}
{\displaystyle
\uprho \big( \big( Q_{\tilde{r}, \tilde\theta}^{j} (\bar{x}, \bar{t}) \big)^{\tilde\varsigma} \big) \leqslant
	\frac{\sf L}{\upkappa} \, \uprho_j \big( \big( Q_{\tilde{r}, \tilde\theta}^{j} (\bar{x}, \bar{t}) \big) \big) ,	 } 													\\	[0.7em]
{\displaystyle
\uprho \big( \big( Q_{\tilde{r}, \tilde\theta}^{j} (\bar{x}, \bar{t}) \big)^{\tilde\varsigma} \big) \leqslant
	\frac{\sf L}{\upkappa} \, \frac{\tilde\theta \, \tilde{r}^2 {\sf h}_{\ast}^j}{R^2} \, \uprho_j \big( B_R (\bar{x}) \times (\bar{t} - R^2 , \bar{t}) \big) \cap \mathcal{Q}_j \big) ,	 } 	\\	[0.7em]
{\displaystyle
\big| \big( Q_{\tilde{r}, \tilde\theta}^{j} (\bar{x}, \bar{t}) \big)^{\tilde\varsigma} \big| \leqslant 
	\frac{1}{\upkappa}  \, \frac{\tilde\theta \, \tilde{r}^2 {\sf h}_{\ast}^j}{R^2} \, \upchi_j \big( B_R (\bar{x}) \times (\bar{t} - R^2 , \bar{t}) \big) \cap \mathcal{Q}_j \big) . }
\end{array}
\end{equation}
Writing $\Big( \uprho_j \big( Q_{\tilde{r}, \tilde\theta}^{j} (\bar{x}, \bar{t}) \big) \Big)^{-\frac{\kappa - 1}{\kappa}}
\Big( \big| Q_{\tilde{r}, \tilde\theta}^{j} (\bar{x}, \bar{t})\big| \Big)^{-\frac{1}{\kappa}}$ as
$\Big( \uprho_j \big( Q_{\tilde{r}, \tilde\theta}^{j} (\bar{x}, \bar{t}) \big) \Big)^{-1}
\Big( \uprho_j \big( Q_{\tilde{r}, \tilde\theta}^{j} (\bar{x}, \bar{t}) \big) / \big| Q_{\tilde{r}, \tilde\theta}^{j} (\bar{x}, \bar{t})\big| \Big)^{\frac{1}{\kappa}}$,
starting from \eqref{partitella} and using \eqref{euno}, \eqref{edue}, \eqref{etre} we get
\begin{align*}
& \frac{1}{\upeta \big( \big( Q_{\tilde{r}, \tilde\theta}^{j} (\bar{x}, \bar{t}) \big)^{\tilde\varsigma} \big)} \iint_{A (\tilde{k})}(u - \tilde{k})_+^2  \upeta\, dx dt \leqslant
	\upgamma^{\frac{2}{\kappa}} \, (2 \gamma + 2) \, \frac{(1 + \upkappa) {\sf L}}{\upkappa^2}
	\left( \frac{1 + \upkappa}{\upkappa} \right)^{\frac{\kappa - 1}{\kappa}}  \cdot																			\nonumber	\\
& \hskip25pt  \cdot {\tilde{r}}^{\frac{2}{\kappa}} \big( \tilde\theta \, {\tilde{r}}^2 {\sf h}_{\ast}^j \big)^{\frac{\kappa - 1}{\kappa}}
	\left( \frac{\upeta \big( \tilde{A} (\tilde{k}) \big)}{\upeta \big( \big( Q_{\tilde{r}, \tilde\theta}^{j} (\bar{x}, \bar{t}) \big)^{\tilde\varsigma} \big)} \right)^{\frac{\kappa - 1}{\kappa}}
	\left( \left(\frac{1 + \upkappa}{\upkappa} \right)^4 \, c_{\upchi_j} (R / \tilde{r}) \, {\sf h}_{\ast}^j \right)^{\frac{1}{\kappa}} \cdot											\nonumber	\\ 
& \hskip25pt  \cdot \Bigg[ \frac{(1 + \upkappa)^2}{\upkappa^3}  \, \frac{c_{\uprho_j} ( R/\tilde{r})}{{\sf h}_{\ast}^j} \, 
\left( \frac{1}{(\tilde{r} - r)^2} + 1 \right) \frac{1}{\big| \big( Q_{\tilde{r}, \tilde\theta}^{j} (\bar{x}, \bar{t}) \big)^{\tilde\varsigma} \big|}
	\iint_{\big( Q_{\tilde{r}, \tilde\theta}^{j} (\bar{x}, \bar{t}) \big)^{\tilde\varsigma}} (u - \tilde{k})_+^2 \, dx dt \, +																	\\
& \hskip50pt  + \frac{\sf L}{\upkappa} \, \frac{1}{\tilde\theta (\tilde{r} - r)^2 {\sf h}_{\ast}^j} 
	\frac{1}{\uprho \big( \big( Q_{\tilde{r}, \tilde\theta}^{j} (\bar{x}, \bar{t}) \big)^{\tilde\varsigma} \big)}
	\iint_{\big( Q_{\tilde{r}, \tilde\theta}^{j} (\bar{x}, \bar{t}) \big)^{\tilde\varsigma}} (u - \tilde{k})_+^2 \uprho \, dx dt \, +												\nonumber	\\ 
& \hskip50pt  + \frac{(1 + \upkappa)^2}{\upkappa^3}  \, \frac{c_{\uprho_j} ( R/\tilde{r})}{{\sf h}_{\ast}^j} \, 
	\tilde{k}^2 \, \left( \frac{1}{(\tilde{r} - r)^2} + 1 \right) \frac{\big| \tilde{A} (\tilde{k}) \big|}{\big| \big( Q_{\tilde{r}, \tilde\theta}^{j} (\bar{x}, \bar{t}) \big)^{\tilde\varsigma} \big|}  \Bigg]	\nonumber
\end{align*}
by which finally
\begin{align}
& \frac{1}{\upeta \big( \big( Q_{\tilde{r}, \tilde\theta}^{j} (\bar{x}, \bar{t}) \big)^{\tilde\varsigma} \big)} \iint_{A (\tilde{k})}(u - \tilde{k})_+^2  \upeta\, dx dt \leqslant
	{\sf C}_1 \, {\tilde{r}}^2 \, \tilde\theta^{\frac{\kappa - 1}{\kappa}}
	\left( \frac{\upeta \big( \tilde{A} (\tilde{k}) \big)}{\upeta \big( \big( Q_{\tilde{r}, \tilde\theta}^{j} (\bar{x}, \bar{t}) \big)^{\tilde\varsigma} \big)} \right)^{\frac{\kappa - 1}{\kappa}} \cdot	\nonumber	\\ 
\label{corsetta}
& \hskip25pt  \cdot \Bigg[ \frac{(1 + \upkappa)^2}{\upkappa^3}  \, c_{\uprho_j} ( R/\tilde{r}) \, 
\left( \frac{1}{(\tilde{r} - r)^2} + 1 \right) \frac{1}{\big| \big( Q_{\tilde{r}, \tilde\theta}^{j} (\bar{x}, \bar{t}) \big)^{\tilde\varsigma} \big|}
	\iint_{\big( Q_{\tilde{r}, \tilde\theta}^{j} (\bar{x}, \bar{t}) \big)^{\tilde\varsigma}} (u - \tilde{k})_+^2 \, dx dt \, +																	\\
& \hskip50pt  + \frac{\sf L}{\upkappa} \, \frac{1}{\tilde\theta (\tilde{r} - r)^2} 
	\frac{1}{\uprho \big( \big( Q_{\tilde{r}, \tilde\theta}^{j} (\bar{x}, \bar{t}) \big)^{\tilde\varsigma} \big)}
	\iint_{\big( Q_{\tilde{r}, \tilde\theta}^{j} (\bar{x}, \bar{t}) \big)^{\tilde\varsigma}} (u - \tilde{k})_+^2 \uprho \, dx dt \, +												\nonumber	\\ 
& \hskip50pt  + \frac{(1 + \upkappa)^2}{\upkappa^3}  \, c_{\uprho_j} ( R/\tilde{r}) \, 
	\tilde{k}^2 \, \left( \frac{1}{(\tilde{r} - r)^2} + 1 \right) \frac{\big| \tilde{A} (\tilde{k}) \big|}{\big| \big( Q_{\tilde{r}, \tilde\theta}^{j} (\bar{x}, \bar{t}) \big)^{\tilde\varsigma} \big|} \Bigg]	\nonumber
\end{align}
where
$$
{\sf C}_1 = \upgamma^{\frac{2}{\kappa}} \, (2 \gamma + 2) \, \frac{(1 + \upkappa) {\sf L}}{\upkappa^2}
	\left( \frac{1 + \upkappa}{\upkappa} \right)^{\frac{\kappa - 1}{\kappa}} 	\left( \left(\frac{1 + \upkappa}{\upkappa} \right)^4 \, c_{\upchi_j} (R / \tilde{r}) \right)^{\frac{1}{\kappa}} \, .
$$
Moreover notice that
\begin{align}
\label{carlofrigna}
(\tilde{k} - k)^2 \, & \upeta \big( \tilde{A} (\tilde{k}) \big) \leqslant
	\iint_{\tilde{A} (\tilde{k})} (u - k)_+^2 \upeta \, dx dt \leqslant \iint_{\big( Q_{\tilde{r}, \tilde\theta}^{j} (\bar{x}, \bar{t}) \big)^{\tilde\varsigma}} (u - k)_+^2 \upeta \, dx dt .
\end{align}

\bprop
\label{carlettobello}
Suppose $(\text{\rm H}.1)$-$(\text{\rm H}.4)$ hold. 
Consider $u \in DG(\Omega, T, \uprho, \gamma)$; $j = 1$ or $j = 2$;
$s, \tilde{s}, R \in (0, \bar{R}]$ with $s < \tilde{s} \leqslant R$; $\vartheta, \tilde{\vartheta}$ with $0 < \vartheta \leqslant \tilde{\vartheta}$; $(\bar{x}, \bar{t}) \in \mathcal{Q}_j \cup \Upgamma$;
suppose moreover $Q_{\tilde{s}, \tilde{\vartheta}} (\bar{x}, \bar{t}) \subset \Omega \times (0,T)$, $Q_{R} (\bar{x}, \bar{t}) \subset \Omega \times (0,T)$ and
$\tilde{\vartheta} \, \tilde{s}^2 {\sf h}_{\ast} ( \bar{x}, \bar{t}, R) < \updelta$. 
Consider $\varsigma_0 = \big( \tilde{s} - s , \big( \tilde\vartheta \, \tilde{s}^2 - \vartheta \, s^2 \big) {\sf h}_{\ast}^j ( \bar{x}, \bar{t}, R) \big)$,
$\varsigma^0 = \big( \tilde{s} - s , \big( \tilde\vartheta \, \tilde{s}^2 - \vartheta \, s^2 \big) {\sf h}^{\ast}_j ( \bar{x}, \bar{t}, R) \big)$. 
Then for every choice of $a, \sigma \in (0,1)$, every
$u \in DG (\Omega, T, \uprho, \gamma)$ and every $\mu, \omega$ satisfying the conditions below there is $\nu$ depending (only) on
$\upgamma, \gamma, \kappa, (1-a), \tilde{s}, \tilde{\vartheta}, \tilde{s} - s, \vartheta, c_{\uprho_j} (R / \tilde{s}), c_{\upchi_j} (R / \tilde{s}), \upkappa , {\sf L}^{-1},
\left(\frac{\mu}{\sigma \omega} - a \right)^{-1}$ such that
under the conditions
\begin{gather*}
Q_{\tilde{s}, \tilde{\vartheta}} (\bar{x}, \bar{t}) \subset \Omega \times (0,T) , \qquad Q_{R} (\bar{x}, \bar{t}) \subset \Omega \times (0,T) ,			\\
\mu \geqslant \sup_{\big( Q_{\tilde{s}, \tilde\vartheta}^{j} (\bar{x}, \bar{t}) \big)^{\varsigma_0}} u , \qquad
	\omega \geqslant \osc_{\big( Q_{\tilde{s}, \tilde\vartheta}^{j} (\bar{x}, \bar{t}) \big)^{\varsigma_0}} u
\end{gather*}
if
$$
\frac{\big|\{ (x,t) \in \big( Q_{\tilde{s}, \tilde\vartheta}^{j} (\bar{x}, \bar{t}) \big)^{\varsigma_0} \, | \, u(x,t) > \mu - \sigma \omega \} \big|}
	{\big| \big( Q_{\tilde{s}, \tilde\vartheta}^{j} (\bar{x}, \bar{t}) \big)^{\varsigma_0} \big|} +
\frac{\uprho \big( \{ (x,t) \in \big( Q_{\tilde{s}, \tilde\vartheta}^{j} (\bar{x}, \bar{t}) \big)^{\varsigma_0} \, | \, u(x,t) > \mu - \sigma \omega \} \big)}
	{\uprho \big( \big( Q_{\tilde{s}, \tilde\vartheta}^{j} (\bar{x}, \bar{t}) \big)^{\varsigma_0} \big)} \leqslant \nu
$$
then
$$
u(x,t) \leqslant \mu - a \, \sigma \, \om \hskip30pt \textit{for a.e. } (x,t) \in Q_{s, \vartheta}^{j} (\bar{x}, \bar{t}) \, .
$$
\eprop

\boss
This is the only result stated for each function belonging to
$DG(\Omega, T, \uprho, \gamma)$ while in all the following we consider non-negative functions. \\
The constant $\nu$ may depend on $\left(\frac{\mu}{\sigma \omega} - a \right)^{-1}$. Having $\mu/\omega$ very little is not a problem since, as one can see
by the proof, $\nu$ precisely depends on $\big( 1  + \left(\frac{\mu}{\sigma \omega} - a \right)^2 \big)^{-1}$ which is always less or equal to 1. \\
The quantity $\left(\frac{\mu}{\sigma \omega} - a \right)^{-1}$ may affect the value of $\nu$ exactly when one considers non-negative functions
and one takes $\mu = \sup u$, $\omega = \osc u$: in this case the value of $\nu$ decreases as the ratio $\mu/\omega$ increases.
Notice that anyway one can always consider, if $u \geqslant 0$, $\mu = \sup u$, $\omega = \mu$.
\eoss

\noindent
\dimo
Consider for $h \in \N$
\begin{gather*}
r_h := s + \frac{\tilde{s} - s}{2^{h}} \searrow s , \hskip30pt \theta_{h} = \theta_h := \vartheta + \frac{\tilde\vartheta - \vartheta}{2^{h}}	\searrow \vartheta , 		\\
\sigma_h = a \, \sigma + \frac{1 - a}{2^h} \, \sigma \searrow a\sigma, \hskip30pt  k_h = \mu - \sigma_h \omega \nearrow \mu - a \sigma \omega	,			\\
\varsigma_h := ( r_h - s, \theta_h r_h^2 - \vartheta s^2) , 																				\\
Q_h^j := \big( Q_{r_h, \theta_h}^{j} (\bar{x}, \bar{t}) \big)^{\varsigma_h} , \hskip30pt A_h^j = \big\{ (x,t) \in Q_h^j \, | \, u (x,t) > k_h \big\} \, .
\end{gather*}
Using
$$
(u - k_h)_+ \leqslant \sigma_h \omega ,
$$
using \eqref{carlofrigna} in \eqref{corsetta} with
\begin{align*}
& r_{h+1}			& \text{in the place of}  &  \hskip10pt	r	\, ,			&	r_{h}	 			& \hskip10pt	\text{in the place of}  		&  \hskip10pt	\tilde{r} \, ,			\nonumber	\\
& \theta_{h+1}		& \text{in the place of}  & \hskip10pt \theta	\, ,			&	\theta_{h}			& \hskip10pt	\text{in the place of}		& \tilde\theta \, ,								\\
& k_{h+1}  		& \text{in the place of}  &  \hskip10pt	\tilde{k}	\, ,		&	k_{h}				& \hskip10pt	\text{in the place of}		& k \, ,						\nonumber	\\
& \varsigma_{h+1}  	& \text{in the place of}  &  \hskip10pt	\tilde{\varsigma} \, ,	&	\varsigma_{h}		& \hskip10pt	\text{in the place of}		& \varsigma \, ,			
\end{align*}	
we get
\begin{align*}
& (k_{h+1} - k_h)^2 \frac{\upeta \big( A^j_{h+1} \big)}{\upeta \big( Q^j_{h+1} \big)} \leqslant
	{\sf C}_1 \, {r_h}^2 \, \theta_h^{\frac{\kappa - 1}{\kappa}}
	\left( \frac{\upeta \big( A^j_{h} \big)}{\upeta \big( Q^j_{h} \big)} \right)^{\frac{\kappa - 1}{\kappa}} \frac{1}{(r_{h-1} - r_h)^2} \, \cdot			\nonumber	\\ 
& \hskip25pt  \cdot \Bigg[ \frac{(1 + \upkappa)^2}{\upkappa^3}  \, c_{\uprho_j} ( R/r_h) \, 
\left( (r_{h-1} - r_{h})^2 + 1 \right) \big( (\sigma_{h} \omega)^2  + k_h^2 \big) \frac{\big| A_{h}^j \big|}{\big| Q_{h}^j \big|} \, +								\\
& \hskip50pt  + \frac{\sf L}{\upkappa} \, \frac{(\sigma_{h} \omega)^2}{\theta_{h}} 
	\frac{\uprho \big( A_{h}^j \big)}{\uprho \big( Q_{h}^j \big)}  \Bigg]															\nonumber
\end{align*}
where $\upeta$ is either $\uprho$ or $1$. By this we get
\begin{align*}
& \frac{\upeta \big( A^j_{h+1} \big)}{\upeta \big( Q^j_{h+1} \big)} \leqslant
	{\sf C}_1 \, \frac{2^{2h+2}}{((1-a) \, \sigma \, \omega)^2} \, {\tilde{s}}^2 \, \tilde\vartheta^{\frac{\kappa - 1}{\kappa}}
	\left( \frac{\upeta \big( A^j_{h} \big)}{\upeta \big( Q^j_{h} \big)} \right)^{\frac{\kappa - 1}{\kappa}} \frac{2^{2h}}{(\tilde{s} - s)^2} \, \cdot		\nonumber	\\ 
& \hskip25pt  \cdot \Bigg[ \frac{(1 + \upkappa)^2}{\upkappa^3}  \, c_{\uprho_j} \left( \frac{R}{s} \right) \, 
\left( (\tilde{s} - s)^2 + 1 \right) \big( (\sigma \omega)^2  + (\mu - a \sigma \omega)^2 \big) \frac{\big| A_{h}^j \big|}{\big| Q_{h}^j \big|} \, +					\\
& \hskip50pt  + \frac{\sf L}{\upkappa} \, \frac{(\sigma \omega)^2}{\vartheta} 
	\frac{\uprho \big( A_{h}^j \big)}{\uprho \big( Q_{h}^j \big)}  \Bigg]															\nonumber
\end{align*}
that may be rewritten as
\begin{align*}
& \frac{\upeta \big( A^j_{h+1} \big)}{\upeta \big( Q^j_{h+1} \big)} \leqslant
	\frac{{\sf C}_2}{2} \, 16^h \, 
	\left( \frac{\upeta \big( A^j_{h} \big)}{\upeta \big( Q^j_{h} \big)} \right)^{\frac{\kappa - 1}{\kappa}}  \, 
	\Bigg[ \frac{\big| A_{h}^j \big|}{\big| Q_{h}^j \big|} \, +  \frac{\uprho \big( A_{h}^j \big)}{\uprho \big( Q_{h}^j \big)}  \Bigg]														
\end{align*}
where
$$
{\sf C}_2 :=  \frac{8 \, {\sf C}_1}{(1-a)^2} \, \frac{{\tilde{s}}^2 \, \tilde\vartheta^{\frac{\kappa - 1}{\kappa}}}{(\tilde{s} - s)^2}
	\max \left\{ 
	\frac{(1 + \upkappa)^2}{\upkappa^3}  \, c_{\uprho_j} \left( \frac{R}{s} \right) \, 
		\left( (\tilde{s} - s)^2 + 1 \right) \Big( 1  + \left(\frac{\mu}{\sigma \omega} - a \right)^2 \Big) , \frac{\sf L}{\upkappa} \, \frac{1}{\vartheta} \, .
	\right\}
$$
Taking first $\upeta = \uprho$, then $\upeta = 1$ and summing the two obtained inequalities one derives
$$
y_{h+1} \leqslant {\sf C}_2 \, 16^h \, y_h^{1 + \frac{\kappa - 1}{\kappa}} \, ,
	\hskip40pt y_h := \frac{\big| A_{h}^j \big|}{\big| Q_{h}^j \big|} \, +  \frac{\uprho \big( A_{h}^j \big)}{\uprho \big( Q_{h}^j \big)}		.		
$$
Using then Lemma \ref{giusti} we get that if
$$
\frac{\big| A_{0}^j \big|}{\big| Q_{0}^j \big|} \, +  \frac{\uprho \big( A_{0}^j \big)}{\uprho \big( Q_{0}^j \big)} \leqslant 16^{-\left(\frac{\kappa}{\kappa - 1}\right)^2} {\sf C}_2^{- \frac{\kappa}{\kappa - 1}}
\qquad \text{ then } \qquad \lim_{h\to +\infty} y_h = 0 \, .
$$
We conclude choosing $\nu = 16^{-\left(\frac{\kappa}{\kappa - 1}\right)^2} {\sf C}_2^{-\frac{\kappa}{\kappa - 1}}$.
\finedimo

\bprop
\label{anninabella}
Consider $u \in DG(\Omega, T, \uprho, \gamma)$; $j = 1$ or $j = 2$;
$s, \tilde{s}, R \in (0, \bar{R}]$ with $s < \tilde{s} \leqslant R$; $\hat\vartheta, \tilde{\vartheta}$ with $0 < \hat\vartheta \leqslant \tilde{\vartheta} s^2 / \tilde{s}^2$;
$(\bar{x}, \bar{t}) \in \mathcal{Q}_j \cup \Upgamma$ and suppose
$Q_{\tilde{s}, \tilde{\vartheta}} (\bar{x}, \bar{t}) \subset \Omega \times (0,T)$, $Q_{R} (\bar{x}, \bar{t}) \subset \Omega \times (0,T)$ and
$\tilde{\vartheta} \, \tilde{s}^2 {\sf h}_{\ast}^j ( \bar{x}, \bar{t}, R) < \updelta$ and $\tilde{\vartheta} \, \tilde{s}^2 {\sf h}^{\ast}_j ( \bar{x}, \bar{t}, R) < \updelta$. 
Consider $\varsigma_0 = \big( \tilde{s} - s , \big( \tilde\vartheta - \hat\vartheta \big) \tilde{s}^2 {\sf h}_{\ast}^j ( \bar{x}, \bar{t}, R) \big)$,
$\varsigma^0 = \big( \tilde{s} - s , \big( \tilde\vartheta - \hat\vartheta \big) \tilde{s}^2 {\sf h}^{\ast}_j ( \bar{x}, \bar{t}, R) \big)$. 
Then for every choice of $a, \sigma \in (0,1)$, every
$u \in DG (\Omega, T, \uprho, \gamma)$ and every $\mu, \omega$ satisfying the conditions below there is $\nu$ depending $($only$)$ on
$\upgamma, \gamma, \kappa, (1-a), \tilde{s}, \tilde{\vartheta}, \tilde{s} - s, \hat\vartheta, c_{\uprho_j} (R / \tilde{s}), c_{\upchi_j} (R / \tilde{s}), \upkappa , {\sf L}^{-1}$ and
on $\left(\frac{\mu_+}{\sigma \omega} - a \right)^{-1}$ in points $i \, )$ and $ii \, )$, and on $\left(\frac{\mu_-}{\sigma \omega} + a \right)^{-1}$  in points $iii \, )$ and $iv \, )$ such that
\begin{itemize}
\item[$i\, )$]
under the conditions
\begin{gather*}
Q_{\tilde{s}, \tilde{\vartheta}} (\bar{x}, \bar{t}) \subset \Omega \times (0,T) , \qquad Q_{R} (\bar{x}, \bar{t}) \subset \Omega \times (0,T) ,			\\
\mu_+ \geqslant \sup_{\big( Q_{\tilde{s}, \tilde\vartheta}^{j} (\bar{x}, \bar{t}) \big)^{\varsigma_0}} u , \qquad
	\omega \geqslant \osc_{\big( Q_{\tilde{s}, \tilde\vartheta}^{j} (\bar{x}, \bar{t}) \big)^{\varsigma_0}} u
\end{gather*}
if
$$
\frac{\big|\{ (x,t) \in \big( Q_{\tilde{s}, \tilde\vartheta}^{j} (\bar{x}, \bar{t}) \big)^{\varsigma_0} \, | \, u(x,t) > \mu - \sigma \omega \} \big|}
	{\big| \big( Q_{\tilde{s}, \tilde\vartheta}^{j} (\bar{x}, \bar{t}) \big)^{\varsigma_0} \big|} +
\frac{\uprho \big( \{ (x,t) \in \big( Q_{\tilde{s}, \tilde\vartheta}^{j} (\bar{x}, \bar{t}) \big)^{\varsigma_0} \, | \, u(x,t) > \mu - \sigma \omega \} \big)}
	{\uprho \big( \big( Q_{\tilde{s}, \tilde\vartheta}^{j} (\bar{x}, \bar{t}) \big)^{\varsigma_0} \big)} \leqslant \nu
$$
then
$$
u(x,t) \leqslant \mu_+ - a \, \sigma \, \om \hskip30pt \textit{for a.e. } (x,t) \in
	\Big(B_{s} (\bar{x}) \times [ \bar{t} - \hat\vartheta \tilde{s}^2 {\sf h}^{\ast}_j ( \bar{x}, \bar{t}, R) , \bar{t} ] \Big) \cap \mathcal{Q}_j ;
$$
\item[$ii\, )$]
under the conditions
\begin{gather*}
Q^{\tilde{s}, \tilde{\vartheta}} (\bar{x}, \bar{t}) \subset \Omega \times (0,T) , \qquad Q^{R} (\bar{x}, \bar{t}) \subset \Omega \times (0,T) ,			\\
\mu_+ \geqslant \sup_{\big( Q^{\tilde{s}, \tilde\vartheta}_{j} (\bar{x}, \bar{t}) \big)^{\varsigma^0}} u , \qquad
	\omega \geqslant \osc_{\big( Q^{\tilde{s}, \tilde\vartheta}_{j} (\bar{x}, \bar{t}) \big)^{\varsigma^0}} u
\end{gather*}
if
$$
\frac{\big|\{ (x,t) \in \big( Q^{\tilde{s}, \tilde\vartheta}_{j} (\bar{x}, \bar{t}) \big)^{\varsigma^0} \, | \, u(x,t) > \mu - \sigma \omega \} \big|}
	{\big| \big( Q^{\tilde{s}, \tilde\vartheta}_{j} (\bar{x}, \bar{t}) \big)^{\varsigma^0} \big|} +
\frac{\uprho \big( \{ (x,t) \in \big( Q^{\tilde{s}, \tilde\vartheta}_{j} (\bar{x}, \bar{t}) \big)^{\varsigma^0} \, | \, u(x,t) > \mu - \sigma \omega \} \big)}
	{\uprho \big( \big( Q^{\tilde{s}, \tilde\vartheta}_{j} (\bar{x}, \bar{t}) \big)^{\varsigma^0} \big)} \leqslant \nu
$$
then
\begin{align*}
& u(x,t) \leqslant \mu_+ - a \, \sigma \, \om \hskip30pt 	\\
& \textit{for a.e. } (x,t) \in \Big(B_{s} (\bar{x}) \times [ \bar{t} + \big( \tilde\vartheta - \hat\vartheta \big) \tilde{s}^2 {\sf h}^{\ast}_j ( \bar{x}, \bar{t}, R) ,
	\bar{t} + \tilde\vartheta \tilde{s}^2 {\sf h}^{\ast}_j ( \bar{x}, \bar{t}, R) ] \Big) \cap \mathcal{Q}_j \, ;
\end{align*}
\item[$iii\, )$]
under the conditions
\begin{gather*}
Q_{\tilde{s}, \tilde{\vartheta}} (\bar{x}, \bar{t}) \subset \Omega \times (0,T) , \qquad Q_{R} (\bar{x}, \bar{t}) \subset \Omega \times (0,T) ,			\\
\mu_- \leqslant \inf_{\big( Q_{\tilde{s}, \tilde\vartheta}^{j} (\bar{x}, \bar{t}) \big)^{\varsigma_0}} u , \qquad
	\omega \geqslant \osc_{\big( Q_{\tilde{s}, \tilde\vartheta}^{j} (\bar{x}, \bar{t}) \big)^{\varsigma_0}} u
\end{gather*}
if
$$
\frac{\big|\{ (x,t) \in \big( Q_{\tilde{s}, \tilde\vartheta}^{j} (\bar{x}, \bar{t}) \big)^{\varsigma_0} \, | \, u(x,t) < \mu + \sigma \omega \} \big|}
	{\big| \big( Q_{\tilde{s}, \tilde\vartheta}^{j} (\bar{x}, \bar{t}) \big)^{\varsigma_0} \big|} +
\frac{\uprho \big( \{ (x,t) \in \big( Q_{\tilde{s}, \tilde\vartheta}^{j} (\bar{x}, \bar{t}) \big)^{\varsigma_0} \, | \, u(x,t) < \mu + \sigma \omega \} \big)}
	{\uprho \big( \big( Q_{\tilde{s}, \tilde\vartheta}^{j} (\bar{x}, \bar{t}) \big)^{\varsigma_0} \big)} \leqslant \nu
$$
then
$$
u(x,t) \geqslant \mu_- + a \, \sigma \, \om \hskip30pt \textit{for a.e. } (x,t) \in
	\Big(B_{s} (\bar{x}) \times [ \bar{t} - \hat\vartheta \tilde{s}^2 {\sf h}^{\ast}_j ( \bar{x}, \bar{t}, R) , \bar{t} ] \Big) \cap \mathcal{Q}_j ;
$$
\item[$iv\, )$]
under the conditions
\begin{gather*}
Q^{\tilde{s}, \tilde{\vartheta}} (\bar{x}, \bar{t}) \subset \Omega \times (0,T) , \qquad Q^{R} (\bar{x}, \bar{t}) \subset \Omega \times (0,T) ,			\\
\mu_- \leqslant \inf_{\big( Q^{\tilde{s}, \tilde\vartheta}_{j} (\bar{x}, \bar{t}) \big)^{\varsigma^0}} u , \qquad
	\omega \geqslant \osc_{\big( Q^{\tilde{s}, \tilde\vartheta}_{j} (\bar{x}, \bar{t}) \big)^{\varsigma^0}} u
\end{gather*}
if
$$
\frac{\big|\{ (x,t) \in \big( Q^{\tilde{s}, \tilde\vartheta}_{j} (\bar{x}, \bar{t}) \big)^{\varsigma^0} \, | \, u(x,t) < \mu + \sigma \omega \} \big|}
	{\big| \big( Q^{\tilde{s}, \tilde\vartheta}_{j} (\bar{x}, \bar{t}) \big)^{\varsigma^0} \big|} +
\frac{\uprho \big( \{ (x,t) \in \big( Q^{\tilde{s}, \tilde\vartheta}_{j} (\bar{x}, \bar{t}) \big)^{\varsigma^0} \, | \, u(x,t) < \mu + \sigma \omega \} \big)}
	{\uprho \big( \big( Q^{\tilde{s}, \tilde\vartheta}_{j} (\bar{x}, \bar{t}) \big)^{\varsigma^0} \big)} \leqslant \nu
$$
then
\begin{align*}
& u(x,t) \geqslant \mu_- + a \, \sigma \, \om \hskip30pt 	\\
& \textit{for a.e. } (x,t) \in \Big(B_{s} (\bar{x}) \times [ \bar{t} + \big( \tilde\vartheta - \hat\vartheta \big) \tilde{s}^2 {\sf h}^{\ast}_j ( \bar{x}, \bar{t}, R) ,
	\bar{t} + \tilde\vartheta \tilde{s}^2 {\sf h}^{\ast}_j ( \bar{x}, \bar{t}, R) ] \Big) \cap \mathcal{Q}_j \, .
\end{align*}
\end{itemize}
\eprop
\noindent
\dimo
To prove the first point it is sufficient to consider $\hat \vartheta$ such that 
$\vartheta s^2 = \hat \vartheta \tilde{s}^2$ where $\vartheta, s, \tilde{s}$ are like in Proposition \ref{carlettobello}, i.e.
$\hat \vartheta = \vartheta s^2 / \tilde{s}^2$. \\
The other points can be proved similarly to the proof of Proposition \ref{carlettobello}.
\finedimo

\ \\
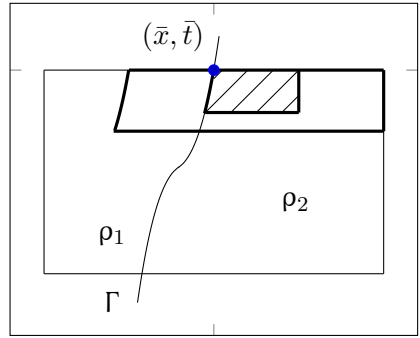
\begin{center}
\pgfplotsset{every axis/.append style={
extra description/.code={ \node at (0.4,0.9) {$(\bar{x}, \bar{t})$}; \node at (0.25,0.1) {$\Upgamma$}; \node at (0.7,0.4) {$\uprho_2$};  \node at (0.25,0.3) {$\uprho_1$};
}}}
\begin{tikzpicture}
\begin{axis} [axis equal, xtick={2}, ytick={1.2}, xticklabels={$$}, yticklabels={$$}, width=7cm, height=6cm, title={Figure 4.a}]
\addplot coordinates
{(2, 1.2)};
\addplot
[domain=0:2.15,variable=\t, very thick]
({t},{-1.2});
\addplot
[domain=2.15:4,variable=\t, smooth]
({t},{-1.2});
\addplot
[domain=-1.2:1.2,variable=\t, very thick]
({0}, {t});
\addplot
[domain=-1.2:1.2,variable=\t, smooth]
({4}, {t});
\addplot
[domain=0:3,variable=\t, very thick]
({t},{1.2});
\addplot
[domain=3:4,variable=\t, smooth]
({t},{1.2});
\addplot
[domain=1.2:1.6,variable=\t,smooth]
({0.6*t^(1/3)+2-0.64}, {t});
\addplot
[domain=0.05:1.2,variable=\t, very thick]
({0.6*t^(1/3)+1.36}, {t});
\addplot
[domain=0.05:0.6,variable=\t, very thick]
({-0.6*t^(1/3)+1.36 + 0.45}, {-t + 0.11});
\addplot
[domain=0.6:1.65,variable=\t, smooth]
({-0.6*t^(1/3)+1.36 + 0.45}, {-t + 0.11});
\addplot
[domain=0.05:1.2,variable=\t, very thick]
({0.6*t^(1/3)+1.36+1}, {t});
\addplot
[domain=0.05:1.31,variable=\t, very thick]
({-0.6*t^(1/3)+1.36 + 0.45+1}, {-t + 0.11});
\addplot
[domain=-0.5:1.2,variable=\t, very thick]
({1}, {t});
\addplot
[domain=1:1.3,variable=\t, very thick]
({t},{-0.5});
\addplot
[domain=1:1.3,variable=\t, smooth]
({t},{t - 0.1});
\addplot
[domain=1:1.6,variable=\t, smooth]
({t},{t - 0.4});
\addplot
[domain=1:1.9,variable=\t, smooth]
({t},{t - 0.7});
\addplot
[domain=1:1.95,variable=\t, smooth]
({t},{t - 1});
\addplot
[domain=1:1.85,variable=\t, smooth]
({t},{t - 1.3});
\addplot
[domain=1.1:1.4,variable=\t, smooth]
({t},{t - 1.6});
\end{axis}
\end{tikzpicture}
\qquad
\begin{tikzpicture}
\begin{axis} [axis equal, xtick={2}, ytick={1.2}, xticklabels={$$}, yticklabels={$$}, width=7cm, height=6cm, title={Figure 4.b}]  
\addplot coordinates
{(2, 1.2)};
\addplot
[domain=0:4,variable=\t, smooth]
({t},{-1.2});
\addplot
[domain=-1.2:1.2,variable=\t, smooth]
({0}, {t});
\addplot
[domain=-1.2:0.7,variable=\t, smooth]
({4}, {t});
\addplot
[domain=0.48:1.2,variable=\t, very thick]
({4}, {t});
\addplot
[domain=0.7:1.2,variable=\t, very thick]
({3}, {t});
\addplot
[domain=1.9:3,variable=\t, very thick]
({t},{0.7});
\addplot
[domain=1:4,variable=\t, very thick]
({t},{1.2});
\addplot
[domain=0:1,variable=\t, smooth]
({t},{1.2});
\addplot
[domain=0.82:4,variable=\t, very thick]
({t},{0.48});
\addplot
[domain=1.2:1.6,variable=\t, smooth]
({0.6*t^(1/3)+2-0.64}, {t});
\addplot
[domain=0.05:1.2,variable=\t, smooth]
({0.6*t^(1/3)+1.36}, {t});
\addplot
[domain=0.69:1.2,variable=\t, very thick]
({0.6*t^(1/3)+1.36}, {t});
\addplot
[domain=0.05:1.31,variable=\t, smooth]
({-0.6*t^(1/3)+1.36 + 0.45}, {-t + 0.11});
\addplot
[domain=1.31:1.65,variable=\t, smooth]
({-0.6*t^(1/3)+1.36 + 0.45}, {-t + 0.11});
\addplot
[domain=0.49:1.2,variable=\t, very thick]
({0.6*t^(1/3)+1.36-1}, {t});
\addplot
[domain=2:2.1,variable=\t, smooth]
({t},{t - 0.9});
\addplot
[domain=1.9:2.4,variable=\t, smooth]
({t},{t - 1.2});
\addplot
[domain=2.2:2.7,variable=\t, smooth]
({t},{t - 1.5});
\addplot
[domain=2.5:3,variable=\t, smooth]
({t},{t - 1.8});
\addplot
[domain=2.8:3,variable=\t, smooth]
({t},{t - 2.1});
\end{axis}
\end{tikzpicture}
\end{center}
\ \\

\noindent
In Figure 4 there is an example of the sets involved in points $i \,)$ and $iii \,)$, while in Figure 5, an example of the sets involved in points $ii \,)$ and $iv \,)$.
In all the examples we suppose
$$
{\sf h}_{\ast}^1 ( \bar{x}, \bar{t}, R) > {\sf h}_{\ast}^2 ( \bar{x}, \bar{t}, R) 
$$
and the dashed set is $Q_{s, \vartheta}^j (\bar{x}, \bar{t}))$ in Figure 4, $Q^{s, \vartheta}_j (\bar{x}, \bar{t}))$ in Figure 5,
the bigger one with the boundary in bold is $\big( Q_{\tilde{s}, \tilde\vartheta}^{j} (\bar{x}, \bar{t}) \big)^{\varsigma_0}$ in Figure 4,
$\big( Q^{\tilde{s}, \tilde\vartheta}_{j} (\bar{x}, \bar{t}) \big)^{\varsigma^0}$ in Figure 5.

\ \\
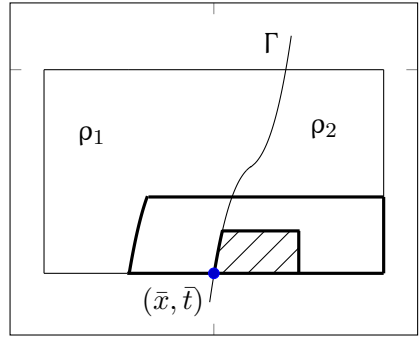
\begin{center}
\pgfplotsset{every axis/.append style={
extra description/.code={\node at (0.4, 0.1) {$(\bar{x}, \bar{t})$}; \node at (0.64, 0.88) {$\Upgamma$}; \node at (0.77,0.62) {$\uprho_2$};  \node at (0.2,0.6) {$\uprho_1$};
}}}
\begin{tikzpicture}
\begin{axis} [axis equal, xtick={2}, ytick={1.2}, xticklabels={$$}, yticklabels={$$}, width=7cm, height=6cm, title={Figure 5.a}]
\addplot coordinates
{(2, -1.2)};
\addplot
[domain=0:3,variable=\t, very thick]
({t},{-1.2});
\addplot
[domain=3:4,variable=\t, smooth]
({t},{-1.2});
\addplot
[domain=-1.2:1.2,variable=\t, very thick]
({0}, {t});
\addplot
[domain=-1.2:1.2,variable=\t, smooth]
({4}, {t});
\addplot
[domain=0:3.85,variable=\t, very thick]
({t},{1.2});
\addplot
[domain=2.85:4,variable=\t, smooth]
({t},{1.2});
\addplot
[domain=1.2:1.6,variable=\t, smooth]
({0.6*t^(1/3)+2-0.64 + 0.85}, {t});
\addplot
[domain=0.5:1.2,variable=\t, smooth]
({0.6*t^(1/3)+1.36 + 0.85}, {t});
\addplot
[domain=0.05:0.5,variable=\t, very thick]
({0.6*t^(1/3)+1.36 + 0.85}, {t});
\addplot
[domain=0.05:0.6,variable=\t, very thick]
({-0.6*t^(1/3)+1.36 + 0.45 + 0.85}, {-t + 0.11});
\addplot
[domain=0.6:1.3,variable=\t, very thick]
({-0.6*t^(1/3)+1.36 + 0.45 + 0.85}, {-t + 0.11});
\addplot
[domain=1.3:1.65,variable=\t, smooth]
({-0.6*t^(1/3)+1.36 + 0.45 + 0.85}, {-t + 0.11});
\addplot
[domain=0.05:1.2,variable=\t, very thick]
({0.6*t^(1/3)+1.36+1.85}, {t});
\addplot
[domain=0.05:1.31,variable=\t, very thick]
({-0.6*t^(1/3)+1.36 + 0.45+1.85}, {-t + 0.11});
\addplot
[domain=-1.2:0.5,variable=\t, very thick]
({1}, {t});
\addplot
[domain=1:2.7,variable=\t, very thick]
({t},{0.5});
\addplot
[domain=1:1.2,variable=\t, smooth]
({t},{t - 0.7});
\addplot
[domain=1:1.5,variable=\t, smooth]
({t},{t - 1});
\addplot
[domain=1:1.8,variable=\t, smooth]
({t},{t - 1.3});
\addplot
[domain=1:2.1,variable=\t, smooth]
({t},{t - 1.6});
\addplot
[domain=1:2.4,variable=\t, smooth]
({t},{t - 1.9});
\addplot
[domain=1:2.7,variable=\t, smooth]
({t},{t - 2.2});
\addplot
[domain=1.3:2.2,variable=\t, smooth]
({t},{t - 2.5});
\addplot
[domain=1.6:2.1,variable=\t, smooth]
({t},{t - 2.8});
\addplot
[domain=1.9:2,variable=\t, smooth]
({t},{t - 3.1});
\end{axis}
\end{tikzpicture}
\qquad
\begin{tikzpicture}
\begin{axis} [axis equal, xtick={2}, ytick={1.2}, xticklabels={$$}, yticklabels={$$}, width=7cm, height=6cm, title={Figure 5.b}]  
\addplot coordinates
{(2, -1.2)};
\addplot
[domain=0:0.98,variable=\t, smooth]
({t},{-1.2});
\addplot
[domain=0.98:4,variable=\t, very thick]
({t},{-1.2});
\addplot
[domain=-1.2:1.2,variable=\t, smooth]
({0}, {t});
\addplot
[domain=-1.2:-0.3,variable=\t, very thick]
({4}, {t});
\addplot
[domain=-0.3:1.2,variable=\t, smooth]
({4}, {t});
\addplot
[domain=1:4,variable=\t, smooth]
({t},{1.2});
\addplot
[domain=0:1,variable=\t, smooth]
({t},{1.2});
\addplot
[domain=1.2:1.6,variable=\t, smooth]
({0.6*t^(1/3)+2-0.64 + 0.85}, {t});
\addplot
[domain=0.05:1.2,variable=\t, smooth]
({0.6*t^(1/3)+1.36 + 0.85}, {t});
\addplot
[domain=0.05:0.1,variable=\t, smooth]
({-0.6*t^(1/3)+1.36 + 0.45 + 0.85}, {-t + 0.11});
\addplot
[domain=0.1:0.4,variable=\t, smooth]
({-0.6*t^(1/3)+1.36 + 0.45 + 0.85}, {-t + 0.11});
\addplot
[domain=0.4:0.92,variable=\t, smooth]
({-0.6*t^(1/3)+1.36 + 0.45 + 0.85}, {-t + 0.11});
\addplot
[domain=0.8:1.3,variable=\t, very thick]
({-0.6*t^(1/3)+1.36 + 0.45 + 0.85}, {-t + 0.11});
\addplot
[domain=1.3:1.65,variable=\t, smooth]
({-0.6*t^(1/3)+1.36 + 0.45 + 0.85}, {-t + 0.11});
\addplot
[domain=0.4:1.31,variable=\t, very thick]
({-0.6*t^(1/3)+1.36 + 0.45 - 0.15}, {-t + 0.11});
\addplot
[domain=-1.2:-0.7,variable=\t, very thick]
({3}, {t});
\addplot
[domain=2.1:3.01,variable=\t, very thick]
({t},{-0.7});
\addplot
[domain=1.23:4,variable=\t, very thick]
({t},{-0.3});
\addplot
[domain=2.05:2.3,variable=\t, smooth]
({t},{t - 3});
\addplot
[domain=2.1:2.6,variable=\t, smooth]
({t},{t - 3.3});
\addplot
[domain=2.4:2.9,variable=\t, smooth]
({t},{t - 3.6});
\addplot
[domain=2.7:3,variable=\t, smooth]
({t},{t - 3.9});
\end{axis}
\end{tikzpicture}
\end{center}
\ \\

\boss
\label{Aacca_ro}
Fix $(\tilde{x}, \tilde{t}) \in \mathcal{Q}_j \cup \Upgamma$ and $h > 0$ and denote
\begin{gather*}
D_{r} (h ; \tilde{t}) (\tilde{x}) := \{ x \in B_{r}(\tilde{x}) \, | \, u(x,\tilde{t}) < h \} ,		\\
D_{r}^j (h ; \tilde{t}) (\tilde{x}) := \{ x \in B_{r}^j (\tilde{x}; \tilde{t}) \, | \, u(x,\tilde{t}) < h \} .
\end{gather*}
First notice that a function $\upeta$ such that $\uprho_j \leqslant \upeta \leqslant \uprho$
satisfies a doubling condition, limited to balls {\em centred} in points $(\tilde{x}, \tilde{t}) \in \mathcal{Q}_j \cup \Upgamma$, with the doubling constant of $\uprho_j$
(this follows by (H.2) - $i \, )$ and (C.7)), i.e.
$$
\upeta (\tilde{t})  \big( B_{2r} (\tilde{x}) \big) \leqslant c_{\uprho_j} \, \upeta (\tilde{t})  \big( B_{r} (\tilde{x}) \big) .
$$
Indeed, using in the order (C.2), (H.2), (C.7), we have
\begin{align*}
\upeta (\tilde{t})  \big( B_{2r} (\tilde{x}) \big) \leqslant \uprho (\tilde{t}) \big( B_{2r} (\tilde{x}) \big) \leqslant c_{\uprho} \, \uprho (\tilde{t}) \big( B_{r} (\tilde{x}) \big)
	\leqslant \frac{c_{\uprho}}{2 \upkappa} \, \uprho_j (\tilde{t}) \big( B_{r} (\tilde{x}) \big)
	\leqslant \frac{c_{\uprho}}{2 \upkappa} \, \upeta (\tilde{t}) \big( B_{r} (\tilde{x}) \big)
\end{align*}
and then we conclude since $c_{\uprho_j} = c_{\uprho} / 2 \upkappa$.
Clearly also the following inequality is true
$$
\upeta (\tilde{t})  \big( B_{2r}^j (\tilde{x}; \tilde{t}) \big) \leqslant c_{\uprho_j} \, \upeta (\tilde{t})  \big( B_{r} (\tilde{x}; \tilde{t}) \big)  \, , \qquad j = 1, 2 ,
$$
which is equivalent to
$\uprho_j (\tilde{t})  \big( B_{2r} (\tilde{x}) \big) \leqslant c_{\uprho_j} \, \uprho_j (\tilde{t})  \big( B_{r} (\tilde{x}) \big)$.
Observe moreover that the condition $u(x,\tilde{t}) \mau h$ for every $x \in B_{r}^j (\tilde{x}; \tilde{t})$ implies
$D_{4 r}(h; \tilde{t}) (\tilde{x}) \subset B_{4 r}(\tilde{x}) \setminus B_{r}^j (\tilde{x}; \tilde{t})$ and then, by the previous inequalities and by (C.7),
\begin{align*}
\uprho_j (\tilde{t}) \big( D_{4 r} (h; \tilde{t}) (\tilde{x}) \big) & \leqslant \uprho_j (\tilde{t}) (B_{4r} (\tilde{x})) - c_{\uprho_j}^{-2}  \uprho_j (\tilde{t}) (B_{4r}^j (\tilde{x}; \tilde{t})) =		\\
& = \big( 1 - c_{\uprho_j}^{-2} \big) \uprho_j (\tilde{t}) (B_{4r}^j (\tilde{x}; \tilde{t})) .
\end{align*}
\eoss

\boss
Clearly the best doubling constant (limited to balls centred in $\mathcal{Q}_1 \cup \Upgamma$) of $\uprho_1$
may be different from the the best doubling constant (limited to balls centred in $\mathcal{Q}_2 \cup \Upgamma$) of $\uprho_2$.
By $c_{\uprho_1}$ and $c_{\uprho_2}$ we denote in fact the same constant, i.e. $c_{\uprho}/ 2 \upkappa$ (see (C.7)),
which may be not the best doubling constant.
\eoss

\begin{lemma}
\label{lemma1}
Consider $j \in \{ 1, 2 \}$, $(\bar{y}, \bar{s}) \in \mathcal{Q}_j \cup \Upgamma$ for which, given $R, r , \theta > 0$, $R \leqslant \bar{R}$,
$B_{4r} (\bar{y}) \subseteq B_R (\bar{y})$, $Q^R (\bar{y}, \bar{s}) = B_R (\bar{y}) \times (\bar{s}, \bar{s} + {\sf h}^{\ast}) \subset \Omega \times (0,T)$ 
$({\sf h}^{\ast} = {\sf h}^{\ast} ( \bar{y}, \bar{s}, R))$,
$Q^{4 r , \theta}(\bar{y}, \bar{s}) = B_{4r} (\bar{y}) \times (\bar{s} , \bar{s} + 16 \, \theta \, r^2 \, {\sf h}^{\ast} ) \subset \Omega \times (0,T)$,
$\theta (4 r)^2 {\sf h}^{\ast}_j ( \bar{y}, \bar{s}, R) < \updelta$.
Then there exist $\vartheta_j \in (0, \theta)$ and $\eta \in (0,1)$ such that,
given $h > 0$ and $u \geqslant 0$ in $DG(\Omega, T, \uprho, \gamma)$ for which
$$
u(x, \bar{s} ) \geqslant h \hskip20pt \text{a.e. in } B_{r}^j (\bar{y}; \bar{s})
$$
then
\begin{gather*}
\uprho_j (t) \big( D_{4r} (\eta h ; t) (\bar{y}) \big) < \left( 1 - \frac{1}{2 c_{\uprho_j}^2} \right) \, \uprho_j (t) \left( B_{4r}(\bar{y}) \right)
		\qquad \forall \ t \in [\bar{s}, \bar{s} + \vartheta_j (4 r)^2 {\sf h}^{\ast}_j ( \bar{y}, \bar{s}, R)] .
\end{gather*}
\end{lemma}

\boss
\label{garofano}
Clearly if $(\bar{y}, \bar{s}) \in \Upgamma$ or $B_r (\bar{y}) \cap \Upgamma \not= \emptyset$ we derive both the estimates
\begin{gather*}
\uprho_1 (t) \big( D_{4r} (\eta h ; t) (\bar{y}) \big) < \left( 1 - \frac{1}{2 c_{\uprho_1}^2} \right) \, \uprho_1 (t) \left( B_{4r}(\bar{y}) \right)
		\qquad \forall \ t \in [\bar{s}, \bar{s} + \vartheta_1 (4 r)^2 {\sf h}^{\ast}_1 ( \bar{y}, \bar{s}, R)] ,										\\
\uprho_2 (t) \big( D_{4r} (\eta h ; t) (\bar{y}) \big) < \left( 1 - \frac{1}{2 c_{\uprho_2}^2} \right) \, \uprho_2 (t) \left( B_{4r}(\bar{y}) \right)
		\qquad \forall \ t \in [\bar{s}, \bar{s} + \vartheta_2 (4 r)^2 {\sf h}^{\ast}_2 ( \bar{y}, \bar{s}, R)]
\end{gather*}
for some $\vartheta_1, \vartheta_2$.
Taking $\bar\vartheta = \min \{ \vartheta_1, \vartheta_2 \}$ we can choose a value independent of $j$.
\eoss
\noindent
\dimo
Fix $j = 1$ or $j = 2$.
Consider
$$
\varsigma = ( 4 r \sigma, (4 r \sigma)^2 {\sf h}^{\ast}_j)
$$
and $\upeta_j^{\varsigma}$ defined as follows:
$$
\upeta_j^{\varsigma} = \uprho \quad \text{in } \big( Q^{4r, \theta}_j (\bar{y}; \bar{s}) \big)^{\varsigma} , \qquad
		\upeta_j^{\varsigma} = 0 \quad \text{in } \mathcal{Q} \setminus \big( Q^{4r, \theta}_j (\bar{y}; \bar{s}) \big)^{\varsigma}.
$$
In the following we denote
$$
\big( D^{r, \theta}_j \big)^{\varsigma} (k) = \big\{ (x,t) \in \big( Q^{4r, \theta}_j (\bar{y}; \bar{s}) \big)^{\varsigma} \, | \, u (x,t) < k \big\} \, .
$$
Now we consider a function $\zeta \in \X_c$ satisfying
\begin{gather*}
\zeta \equiv 1 \hskip10pt \text{in } Q^{4r (1 - \sigma), \theta / (1 - \sigma)^2}_j (\bar{y}; \bar{s}) \, , \qquad 
		\zeta \equiv 0 \hskip10pt \text{ outside of } \big( Q^{4r, \theta}_j (\bar{y}; \bar{s}) \big)^{\varsigma} \, , 	\\
0 \leqslant \zeta \leqslant 1 \, , \hskip15pt {\displaystyle |D \zeta| \leqslant \frac{1}{4 r \sigma} }  \, , \hskip15pt {\displaystyle 0 \leqslant \zeta_t  \leqslant \frac{1}{(4 r \sigma)^2 {\sf h}^{\ast}_j} } \, .
\end{gather*}
Notice that
\begin{align*}
Q^{4r (1 - \sigma), \theta / (1 - \sigma)^2}_j (\bar{y}; \bar{s}) =
	\bigcup_{t \in [\bar{s} , \bar{s} + 16 \, \theta \, r^2 \, {\sf h}^{\ast}_j ( \bar{y}, \bar{s}, R )]} B_{4r(1 - \sigma)}^j (\bar{y}; \bar{s})  \, .				\\
\end{align*}
Consider $\vartheta > 0$ to be fixed below. 
We apply now the energy estimate \eqref{lavagna_bis} with this choice of $\zeta$ to the function $(u - h)_-$ and get
\begin{align*}
& \sup_{t \in (\bar{s}, \bar{s} + \vartheta (4 r)^2  {\sf h}^{\ast}_j)} \int_{B_{4r (1 - \sigma)} (\bar{y})} (u - h)_-^2 \uprho_j (x, t) \, dx   \leqslant
	\int_{B_{4r} (\bar{y})}  (u - h)_-^2 \upeta_j^{\varsigma} (x, \bar{s}) dx \, +  	\nonumber	\\
& \qquad + \frac{\gamma}{2} \left( \frac{1}{(4 r \sigma)^2} + 1 \right) \Bigg[  \iint_{\big( Q^{4r, \theta}_j (\bar{y}; \bar{s}) \big)^{\varsigma}} (u - h)_-^2 \, dx dt
		+ h^2 \big| \big( D^{4r, \vartheta}_j \big)^{\varsigma} (h) \big| \Bigg]  \, +  	\nonumber	\\
& \qquad + \frac{\gamma}{2} \frac{1}{(4 r \sigma)^2 {\sf h}^{\ast}_j} \iint_{\big( Q^{4r, \theta}_j (\bar{y}; \bar{s}) \big)^{\varsigma}} (u - h)_-^2 \upeta_j^{\varsigma} \, dx dt  .		
\end{align*}
Now we use Remark \ref{Aacca_ro} and
the fact that $u(x, \bar{s} ) \geqslant 0$ implies $(u - h)_- (x, \bar{s}) \leqslant h$ to estimate the first term on the right hand side. We have
\begin{align*}
& \int_{B_{4r} (\bar{y})}  (u - h)_-^2 \upeta_j^{\varsigma} (x, \bar{s}) dx = 
	\int_{B_{4r} (\bar{y})}  (u - h)_-^2 \uprho_j (x, \bar{s}) dx + \int_{B_{4r} (\bar{y})}  (u - h)_-^2 (\upeta_j^{\varsigma} - \uprho_j) (x, \bar{s}) dx \leqslant			\\
& \qquad \qquad \leqslant h^2 \big( 1 - c_{\uprho_j}^{-2} \big) \uprho_j (\bar{s}) (B_{4r}^j (\bar{y}; \bar{s}))
	+ h^2 (\upeta_j^{\varsigma} (\bar{s}) - \uprho_j (\bar{s}) ) (B_{4r}^j (\bar{y}; \bar{s}))
\end{align*}
by which
\begin{align}
& \sup_{t \in (\bar{s}, \bar{s} + \vartheta (4 r)^2  {\sf h}^{\ast}_j)} \int_{B_{4r (1 - \sigma)}(\bar{y})} (u - h)_-^2 \uprho_j (x, t) \, dx    \leqslant		\nonumber	\\
\label{anninalamoremio}
& \qquad \qquad \leqslant h^2 \big( 1 - c_{\uprho_j}^{-2} \big) \uprho_j (\bar{s}) (B_{4r}^j (\bar{y}; \bar{s}))
	+ h^2 (\upeta_j^{\varsigma} (\bar{s}) - \uprho_j (\bar{s}) ) (B_{4r}^j (\bar{y}; \bar{s})) \, +													\\
& \qquad \qquad \qquad + \gamma \left( \frac{1}{(4 r \sigma)^2} + 1 \right) h^2 \big| \big( D^{4r, \vartheta}_j \big)^{\varsigma} (h) \big|		
	+ \frac{\gamma}{2} \frac{h^2}{(4 r \sigma)^2 {\sf h}^{\ast}_j} \, \uprho \big( \big( D^{4r, \vartheta}_j \big)^{\varsigma} (h) \big) .			\nonumber
\end{align}
Writing, for some $\eta \in (0,1)$ to be fixed and for $t \in [\bar{s}, \bar{s} + \vartheta (4 r)^2  {\sf h}^{\ast}_j]$,
$$
D_{4r} (\eta h ; t) (\bar{y}) = D_{4r (1 - \sigma)} (\eta h ; t) (\bar{y}) \cup
	\left( \{x \in B_{4r} (\bar{y}) \setminus B_{4r (1-\sigma)}(\bar{y}) \, | \,  u(x,t) < \eta \} \right)
$$
we derive that
\begin{align*}
\uprho_j (t) \big( D_{4r} (\eta h ; t) (\bar{y}) \big) \leqslant \uprho_j (t) \big( D_{4r (1 - \sigma)} (\eta h ; t) (\bar{y}) \big) + \uprho_j (t) \big(B_{4r} (\bar{y}) \setminus B_{4r (1-\sigma)}(\bar{y}) \big).
\end{align*}
Then, on the other hand, for each $t \in [\bar{s}, \bar{s} + \vartheta (4 r)^2  {\sf h}^{\ast}_j]$ we get
\begin{align*}
& \int_{B_{4r (1 - \sigma)}(\bar{y})} (u - h)_-^2 \uprho_j (x, t) \, dx \geqslant \int_{D_{4r (1 - \sigma)} (\eta h ; t) (\bar{y})} (u - h)_-^2 \uprho_j (x, t) \, dx \geqslant				\\
& \hskip50pt \geqslant h^2 (1 - \eta)^2 \uprho_j (t) \big( D_{4r (1 - \sigma)} (\eta h ; t) (\bar{y}) \big) .
\end{align*}
Finally, by this last inequality and \eqref{anninalamoremio}, we obtain
\begin{align}
\label{stimasigma} 
& \uprho_j (t) \big( D_{4r} (\eta h ; t) (\bar{y}) \big) \leqslant \uprho_j (t) \big( D_{4r (1 - \sigma)} (\eta h ; t) (\bar{y}) \big) + 
	\uprho_j (t) \big(B_{4r} (\bar{y}) \setminus B_{4r (1-\sigma)}(\bar{y}) \big)	\leqslant												\nonumber	\\
& \qquad \leqslant \frac{1}{h^2 (1 - \eta)^2} \int_{B_{4r (1 - \sigma)}(\bar{y})} (u - h)_-^2 \uprho_j (x, t) \, dx + 
			\uprho_j (t) \big(B_{4r} (\bar{y}) \setminus B_{4r (1-\sigma)}(\bar{y}) \big) 	\leqslant										\nonumber	\\
& \qquad \leqslant \frac{1}{(1 - \eta)^2} \bigg[ \big( 1 - c_{\uprho_j}^{-2} \big) \uprho_j (\bar{s}) (B_{4r}^j (\bar{y}; \bar{s}))
			+ (\upeta_j^{\varsigma} (\bar{s}) - \uprho_j (\bar{s}) ) (B_{4r}^j (\bar{y}; \bar{s})) \, +													\\
& \qquad \qquad + \gamma \left( \frac{1}{(4 r \sigma)^2} + 1 \right) \big| \big( D^{4r, \vartheta}_j \big)^{\varsigma} (h) \big|
			+ \frac{\gamma}{2} \frac{1}{(4 r \sigma)^2 {\sf h}^{\ast}_j} \, \uprho \big( \big( D^{4r, \vartheta}_j \big)^{\varsigma} (h) \big) \bigg] +	\nonumber	\\
& \qquad \qquad + \uprho_j (t) \big(B_{4r} (\bar{y}) \setminus B_{4r (1-\sigma)}(\bar{y}) \big) 											\nonumber
\end{align}
for every $t \in [\bar{s}, \bar{s} + \vartheta (4 r)^2 {\sf h}^{\ast}_j]$. \\
Now we argue by contradiction: if the thesis were not true
then for every $\vartheta, \eta \in (0,1)$ we would have $\tau \in [\bar{s}, \bar{s} + \vartheta (4 r)^2  {\sf h}^{\ast}_j]$ for which
$$
\uprho_j (\tau) \big( D_{4r} (\eta h ; \tau) (\bar{y}) \big) \geqslant \left( 1 - \frac{1}{2 c_{\uprho_j}^2} \right) \uprho_j (\tau) \left( B_{4r}(\bar{y}) \right).
$$
By this and \eqref{stimasigma} (with $t = \tau$) we would get
\begin{align*}
& \left( 1 - \frac{1}{2 c_{\uprho_j}^2} \right) \uprho_j (\tau) \left( B_{4r}(\bar{y}) \right) \leqslant \uprho_j (\tau) \big(B_{4r} (\bar{y}) \setminus B_{4r (1-\sigma)}(\bar{y}) \big) \, +	\\
& \qquad \qquad + \frac{1}{(1 - \eta)^2} \bigg[ \big( 1 - c_{\uprho_j}^{-2} \big) \uprho_j (\bar{s}) (B_{4r}^j (\bar{y}; \bar{s}))
			+ (\upeta_j^{\varsigma} (\bar{s}) - \uprho_j (\bar{s}) ) (B_{4r}^j (\bar{y}; \bar{s})) \, +															\\
& \qquad \qquad + \gamma \left( \frac{1}{(4 r \sigma)^2} + 1 \right) \big| \big( D^{4r, \vartheta}_j \big)^{\varsigma} (h) \big|
			+ \frac{\gamma}{2} \frac{1}{(4 r \sigma)^2 {\sf h}^{\ast}_j} \, \uprho \big( \big( D^{4r, \vartheta}_j \big)^{\varsigma} (h) \big) \bigg] .
\end{align*}
Then we could choose $\vartheta$ depending on $\sigma$ and going to zero when $\sigma \to 0^+$ in such a way that
$$
\lim_{\sigma \to 0^+} \frac{\big| \big( D^{4r, \vartheta}_j \big)^{\varsigma} (h) \big|}{\sigma^2} = 0
\qquad \text{and} \qquad
\lim_{\sigma \to 0^+} \frac{\uprho \big( \big( D^{4r, \vartheta}_j \big)^{\varsigma} (h) \big)}{\sigma^2} = 0 \, .
$$
Moreover notice that 
$$
\lim_{\sigma \to 0^+} (\upeta_j^{\varsigma} (\bar{s}) - \uprho_j (\bar{s}) ) (B_{4r}^j (\bar{y}; \bar{s})) = 0 .
$$
By the uniform continuity (assumption (H.2)-$\, i \,)$) of the function $[0,T] \ni t \mapsto \uprho_j (t) (B_{4r} (\bar{y}))$
we have that for every $\epsilon > 0$ there is $\delta > 0$ such that
\begin{gather*}
\text{if } \quad \vartheta_j (4 r)^2 {\sf h}^{\ast}_j ( \bar{y}, \bar{s}, R) < \delta	\qquad \text{then}		\\
\big| \uprho_j (t_2)(B_{4r} (\bar{y})) - \uprho_j (t_1)(B_{4r} (\bar{y})) \big| < \epsilon
		\qquad \text{for every } t_1, t_2 \in [\bar{s}, \bar{s} + \vartheta_j (4 r)^2 {\sf h}^{\ast}_j ( \bar{y}, \bar{s}, R)] .
\end{gather*}
Then, adding and subtracting on the left hand side the term $\uprho_j (\bar{s})(B_{4r} (\bar{y}))$, letting $\sigma$ go to zero we would get
\begin{align*}
\left( 1 - \frac{1}{2 c_{\uprho_j}^2} \right) \, \uprho_j (\bar{s})(B_{4r} (\bar{y})) \leqslant 
	\frac{1}{(1 - \eta)^2} \left( 1 - \frac{1}{c_{\uprho_j}^2} \right) \, \uprho_j (\bar{s})(B_{4r} (\bar{y})) + \left( 1 - \frac{1}{2 c_{\uprho_j}^2} \right) \, \epsilon .
\end{align*}
Since this inequality holds independently of the values of $\eta$ and $\epsilon$ we would have a contradiction.
Indeed we could choose $\eta$ and $\epsilon$ satisfying
\begin{gather*}
\frac{1}{(1 - \eta)^2} \left( 1 - \frac{1}{c_{\uprho_j}^2} \right) < \alpha \, \left( 1 - \frac{1}{2 c_{\uprho_j}^2} \right) , 			\\
\left( 1 - \frac{1}{2 c_{\uprho_j}^2} \right) \frac{\epsilon}{\uprho_j (\bar{s}) \left( B_{4r}(\bar{y}) \right)} < (1 - \alpha) \left( 1 - \frac{1}{2 c_{\uprho_j}^2} \right)
\end{gather*}
for some $\alpha \in (0,1)$.
These imply to choose $\eta$ and $\epsilon$ satisfying
\begin{gather*}
(1 - \eta)^2 > \frac{1}{\alpha} \left( 1 - \frac{1}{c_{\uprho_j}^2} \right) \left( 1 - \frac{1}{2 c_{\uprho_j}^2} \right)^{-1}  , 	
	\hskip30pt		\epsilon < (1 - \alpha) \, \uprho_j (\bar{s}) \left( B_{4r}(\bar{y}) \right) .
\end{gather*}
Once chosen $\eta$ and $\alpha$ satisfying the first inequality one can choose $\epsilon$ and consequently $\vartheta_j$. \\
Then there exists $\vartheta_j$ such that the thesis is true for each $t \in [\bar{s}, \bar{s} + \vartheta_j (4 r)^2 {\sf h}^{\ast}_j ( \bar{y}, \bar{s}, R)]$. \\
Notice that, since $\alpha$ is free and $\left( 1 - c_{\uprho_j}^{-2} \right) \left( 1 - 2^{-1} c_{\uprho_j}^{-2} \right)^{-1} < 1$,
$\eta$ does not depend on $c_{\uprho_j}$, while $\vartheta_j$ does.
Notice moreover that, since
\begin{gather*}
{\displaystyle
\inf \, \big\{ \uprho_j (s) \left( B_{4r}(y) \right) \, \big| \, (y,s) \in \mathcal{Q}_j \cup \Upgamma , B_{4r}(y) \subset \Omega \big\} > 0 ,		}		\\
{\displaystyle
\sup \, \big\{ \uprho_j (s) \left( B_{4r}(y) \right) \, \big| \, (y,s) \in \mathcal{Q}_j \cup \Upgamma, B_{4r}(y) \subset \Omega \big\} < + \infty } ,
\end{gather*}
one can choose the value of $\vartheta_j$ to be independent of $\uprho_j (\bar{s}) \left( B_{4r}(\bar{y}) \right)$. 
\finedimo

\begin{lemma}
\label{lemma2}
Consider $j \in \{ 1, 2 \}$, $r, \theta > 0$, $(\bar{y}, \bar{s}) \in \mathcal{Q}_j \cup \Upgamma$
for which  $Q^{R} (\bar{y}, \bar{s}) = B_{R} (\bar{y}) \times (\bar{s}, \bar{s} + R^2 {\sf h}^{\ast}_j ) \subset \Omega \times (0,T)$ with
$R = 5r$ where ${\sf h}^{\ast}_j = {\sf h}^{\ast}_j (\bar{y}, \bar{s}, 5r)$ and ${\sf h}_{\ast}^j = {\sf h}_{\ast}^j (\bar{y}, \bar{s}, 5r)$,
and $Q^{4 r , \theta}(\bar{y}, \bar{s}) = B_{4r} (\bar{y}) \times (\bar{s} , \bar{s} + \theta (4r)^2 {\sf h}^{\ast}_j ) \subset \Omega \times (0,T)$.
Suppose moreover $5r \leqslant \bar{R}$ and $\theta (4r)^2 {\sf h}^{\ast}_j < \updelta$.
Consider $u \in DG(\Omega, T, \uprho, \gamma)$, $u \geqslant 0$, $h > 0$, such that
$$
u(x, \bar{s} ) \geqslant h \hskip20pt \text{a.e. in } B_{r}^j (\bar{y}; \bar{s}) .
$$
Consider $\vartheta_j \in (0, \theta)$ as in Lemma $\ref{lemma1}$ and suppose moreover that
$B_{5r} (\bar{y}) \times (\bar{s} - \vartheta_j (4r)^2 {\sf h}_{\ast}^j, \bar{s} + \vartheta_j (4r)^2 {\sf h}^{\ast}_j) \subset \Omega \times (0,T)$. 
Finally consider $\varsigma^0 = (\varsigma^0_1, \varsigma^0_2) \in \R^2$ with $\varsigma^0_1 < 4r$ and $\varsigma^0_2 < \vartheta_j (4r)^2 {\sf h}^{\ast}_j$. \\
Then there is $b \geqslant 1$ and for every $\epsilon > 0$ there exists
$\eta_o = \eta_o (\epsilon, \upkappa, \gamma^{-1}, \upbeta^{-1}, c_{\uprho_j}^{-1}, \vartheta_j, r, \frac{{\sf h}^{\ast}_j}{{\sf h}_{\ast}^j + {\sf h}^{\ast}_j} ) \in (0,1)$
such that
\begin{gather*}
\uprho \big( \big\{ u < \eta_o h \big \} \cap \big( Q^{4r, \vartheta_j}_j (\bar{x}, \bar{t}) \big)^{\varsigma^0} \big) <
	\epsilon \, \uprho \big(  \big( Q^{4r, \vartheta_j}_j (\bar{x}, \bar{t}) \big)^{\varsigma^0} \big) \, ,				\\
\big| \big\{ u < \eta_o h \big \} \cap \big( Q^{4r, \vartheta_j}_j (\bar{x}, \bar{t}) \big)^{\varsigma^0} \big| < {\sf a} \, \epsilon^{\frac{1}{b}}
	\big| \big( Q^{4r, \vartheta_j}_j (\bar{x}, \bar{t}) \big)^{\varsigma^0} \big| \, ,
\end{gather*}
where ${\sf a} = \big( K_ 2 \, {\sf L} \big)^{1/b} / 2 \upkappa$.
\end{lemma}

\boss
The dependence of $\eta_o$ on the parameters expressed in the previous statement means that the more one of those parameters is small,
the more $\eta_o$ is small. \\
Notice in particular the dependence of $\eta_o$ on $\frac{{\sf h}^{\ast}_j}{{\sf h}_{\ast}^j + {\sf h}^{\ast}_j}$ and on $\vartheta_j$,
that in turn depends on $\uprho_j (\bar{s}) \left( B_{4r}(\bar{y}) \right)$. \\
In particular $\eta_o$ decreases as $\uprho_j (\bar{s}) \left( B_{4r}(\bar{y}) \right)$ decreases or as
${\sf h}_{\ast}^j$ decreases with respect to ${\sf h}^{\ast}_j$.
\eoss

\noindent
\dimo
Consider $\eta > 0$, $m \in \N$ and, for $u$ and $h$ as in the assumptions, consider the function $(u - \eta h 2^{-m})_-$.
Fix $j \in \{1, 2 \}$ and let $\vartheta_j$ the value established in Lemma \ref{lemma1}.				\\
Denote for $s \in (\bar{s}, \bar{s} + \theta (4r)^2 {\sf h}^{\ast})$
\begin{gather*}
A_{4r} (h ; s) (y) := \{ x \in B_{4r} (y)  \, | \, u(x,s) \geqslant h \},	\\
D_{4r} (h ; s) (y) := \{ x \in B_{4r} (y)  \, | \, u(x,s) < h \} .
\end{gather*}
Consider $\varsigma^0 = (\varsigma^0_1, \varsigma^0_2) \in \R^2$ with $\varsigma^0_1 < 4r$ and $\varsigma^0_2 < \vartheta_j (4r)^2 {\sf h}^{\ast}(\bar{y}, \bar{s}, 5r)$
and $\upeta_j^{\varsigma^0}$ defined as 
$$
\upeta_j^{\varsigma^0} = \uprho \quad \text{in } \big( Q^{4r, \theta}_j (\bar{y}; \bar{s}) \big)^{\varsigma} , \qquad
		\upeta_j^{\varsigma^0} = 0 \quad \text{in } \mathcal{Q} \setminus \big( Q^{4r, \theta}_j (\bar{y}; \bar{s}) \big)^{\varsigma}.
$$
where $Q^{{4r}, \vartheta_j}_{j} (\bar{y}, \bar{s}) = B_{4r} (\bar{y}) \times (\bar{s}, \bar{s} + \vartheta_j (4r)^2 {\sf h}^{\ast}_j )$ with
${\sf h}^{\ast}_j = {\sf h}^{\ast}_j (\bar{y}, \bar{s}, 5r)$.
Taking $k = \eta h /2^m$, $l = \eta h /2^{m-1}$, $p \in (1,2)$ and
$\upeta = \upeta_j^{\varsigma^0}$ in Lemma \ref{lemma2.2} we get
for each $\sigma \in (\bar{s}, \bar{s} + \theta (4r)^2 {\sf h}^{\ast})$
\begin{equation}
\label{servesubito}
\arst
\begin{array}{l}
{\displaystyle \int_{B_{4 r} (\bar{y})} \Big( u - \frac{\eta h }{2^{m}} \Big)_- \upeta_j^{\varsigma^0} (x,\sigma) \, dx \leqslant
	\frac{\eta h }{2^{m}} \, \upeta_j^{\varsigma^0} (\sigma) \big( D_{4r} (\eta h\, 2^{-m} ; \sigma) (\bar{y}) \big) \leqslant	}			 \\
\hskip20pt
{\displaystyle
\leqslant \frac{2 \, \upbeta \, r \, \big( \upeta_j^{\varsigma^0} (\sigma) \big( B_{4r} (\bar{y}) \big) \big)^2}{\upeta_j^{\varsigma^0} (\sigma) \big( A_{4r} (\eta h\, 2^{-m+1} ; \sigma) (\bar{y}) \big)}
\left( \frac{1}{|B_{4r}|} \int_{D_{4r} (\eta h\, 2^{-m+1} ; \sigma) (\bar{y}) \setminus D_{4r} (\eta h\, 2^{-m} ; \sigma) (\bar{y}) } \!\!\!\!\!\!\!\!\!\!\!\!\!\!\!\!\!\!\!\!\! |Du (x, \sigma)|^p \, dx \right)^{\frac{1}{p}} }.
\end{array}
\end{equation}
Now, since (at least for $m \geqslant 2$)
$A_{4r} (\eta h\, 2^{-m+1} ; \sigma) (\bar{y}) \supset B_{4 r}(\bar{y})  \setminus D_{4r} (\eta h ; \sigma) (\bar{y})$,
if we consider $\eta \in (0,1)$ as in Lemma \ref{lemma1} we get
\begin{align*}
\upeta_j^{\varsigma^0} (\sigma) \big( A_{4r} (\eta h\, 2^{-m+1} ; \sigma) (\bar{y}) \big) \geqslant
\uprho_j (\sigma) \big( A_{4r} (\eta h\, 2^{-m+1} ; \sigma) (\bar{y}) \big) > \frac{1}{2 c_{\uprho_j}^2} \, \uprho_j (\sigma) \left( B_{4r}(\bar{y}) \right)
\end{align*}
for every $\sigma \in [\bar{s}, \bar{s} + \vartheta_j (4r)^2 {\sf h}^{\ast}_j] \subset [\bar{s}, \bar{s} + \theta (4r)^2 {\sf h}^{\ast}]$.
Then by that and \eqref{servesubito} and (C.8) we derive
\begin{equation}
\label{servedopo}
\arst
\begin{array}{l}
{\displaystyle \
\int_{B_{4 r} (\bar{y})} \Big( u - \frac{\eta h }{2^{m}} \Big)_- \upeta_j^{\varsigma^0} (x,\sigma) \, dx \leqslant	
	\frac{\sf L}{\upkappa} \, 4 \, \upbeta \, c_{\uprho_j}^2 \, r \, \upeta_j^{\varsigma^0} (\sigma) \big( B_{4r} (\bar{y}) \big) \, 	\cdot }		\\
\hskip80pt 
{\displaystyle
\cdot \left( \frac{1}{|B_{4r}|} \int_{D_{4r} (\eta h\, 2^{-m+1} ; \sigma) (\bar{y}) \setminus D_{4r} (\eta h\, 2^{-m} ; \sigma) (\bar{y}) } \!\!\!\!\!\!\!\!\!\!\!\!\!\!\!\!\!\!\!\!\! |Du (x, \sigma)|^p \, dx \right)^{\frac{1}{p}} .
}
\end{array}
\end{equation}
Now we define the sequences
\begin{align*}
a_m^j := \int_{\bar{s}}^{\bar{s} + \vartheta_j (4r)^2 {\sf h}^{\ast}_j} \!\!\!\! \upeta_j^{\varsigma^0} (\sigma) \big( D_{4r} (\eta h\, 2^{-m} ; \sigma) (\bar{y}) \big) d \sigma & =
	\upeta_j^{\varsigma^0} \Big( Q^{4 r , \vartheta_j \, }(\bar{y}, \bar{s}) \cap \Big\{ u < \frac{\eta h}{2^{m}} \Big \} \Big) =										\\
\ & = \uprho \Big(  \big( Q^{4r, \vartheta_j}_j (\bar{x}, \bar{t}) \big)^{\varsigma^0} \cap \Big\{ u < \frac{\eta h}{2^{m}} \Big \} \Big) \, ,									\\
b_m := \int_{\bar{s}}^{\bar{s} + \vartheta_j (4r)^2 {\sf h}^{\ast}_j} \!\! \big|  D_{4r} (\eta h\, 2^{-m} ; \sigma) (\bar{y}) \big| d \sigma \, .
\end{align*}
Since $(u - \alpha)_- \geqslant \alpha / 2$ in $\{ u < \alpha / 2 \}$, integrating the left hand side in \eqref{servedopo} we get
\begin{align}
\label{primopasso}
\int_{\bar{s}}^{\bar{s} + \vartheta_j (4r)^2 {\sf h}^{\ast}_j} \!\!\!\! \int_{B_{4 r} (\bar{y})} \Big( u - \frac{\eta h }{2^{m}} \Big)_- \upeta_j^{\varsigma^0} \, dx d \sigma
	\geqslant \frac{\eta h}{2^{m+1}} \, a_{m+1}^j \, .
\end{align}
While integrating the right hand side and using first (H.2), then (C.5) and (C.8)
\begin{align}
& \int_{\bar{s}}^{\bar{s} + \vartheta_j (4r)^2 {\sf h}^{\ast}_j} \!\! \upeta_j^{\varsigma^0} (\sigma) \big( B_{4r} \big)
\left( \frac{1}{|B_{4r}|} \int_{D_{4r} (\eta h\, 2^{-m+1} ; \sigma) (\bar{y}) \setminus D_{4r} (\eta h\, 2^{-m} ; \sigma) (\bar{y}) } \!\!\!\!\!\!\!\!\!\!\!\!\!\!\!\!\!\!\!\!\! 
	|Du (x, \sigma)|^p \, dx \right)^{\frac{1}{p}} d \sigma \leqslant																				\nonumber	\\
\label{secondopasso}
& \hskip20pt \leqslant \max_{\bar{s} \leqslant \sigma \leqslant \bar{s} + \vartheta_j (4r)^2 {\sf h}^{\ast}_j} \upeta_j^{\varsigma^0} (\sigma) \big( B_{4r} \big) 
	\frac{\big( \vartheta_j (4r)^2 {\sf h}^{\ast}_j \big)^{\frac{p - 1}{p}} }{|B_{4r}|^{\frac{1}{p}}}  \left( \int_{\bar{s}}^{\bar{s} + \vartheta_j \, (4r)^2 {\sf h}^{\ast}_j} \!\!\!\!
	\int_{D_{4r} (\eta h\, 2^{-m+1} ; \sigma) (\bar{y}) \setminus D_{4r} (\eta h\, 2^{-m} ; \sigma) (\bar{y}) } \!\!\!\!\!\!\!\!\!\!\!\!\!\!\!\!\!\!\!\!\! |Du|^p \, dx d \sigma \right)^{\frac{1}{p}}  \leqslant	\\
& \hskip20pt \leqslant \frac{{\sf L}^2}{\upkappa}  \, 
	\frac{\uprho_j \big( B_{4r} (\bar{y}) \times (\bar{s}, \bar{s} + \vartheta_j (4r)^2 {\sf h}^{\ast}_j) \big)}{\big( \vartheta_j (4r)^2 {\sf h}^{\ast}_j |B_{4r}|\big)^{\frac{1}{p}}}
	\left( \iint_{Q^{4 r , \vartheta_j}} \Big|D \Big( u - \frac{\eta h }{2^{m}} \Big)_-  \Big|^2 \, dx d \sigma \right)^{\frac{1}{2}}  (b_{m-1} - b_m)^{\frac{2 - p}{2p}} .			\nonumber
\end{align}
Now we estimate $\iint_{Q^{4 r , \vartheta_j}} | D ( u - \frac{\eta h }{2^{m}} )_- |^2 \, dx d \sigma$ using the energy estimate \eqref{lavagna_bis}
in the set $B_{5r} (\bar{y}) \times (\bar{s} - \vartheta_j (4r)^2 {\sf h}_{\ast}^j, \bar{s} + \vartheta_j (4r)^2 {\sf h}^{\ast}_j)$.
By the energy inequality \eqref{lavagna_bis} we get
\begin{align}
& \int_{\bar{s}}^{\bar{s} + \vartheta_j (4r)^2 {\sf h}^{\ast}_j} \!\!\!\! \int_{B_{4 r}}
	{\displaystyle	\Big|D\Big(u - \frac{\eta h}{ 2^m}\Big)_-\Big|^2(x,t) dx dt} \leqslant 															\nonumber	\\
& \qquad {\displaystyle \leqslant \frac{\gamma}{2} \, {{\eta^2h^2}\over{2^{2m}}} \left[ 2 \left( 1 + \frac{1}{r^2} \right) \vartheta_j (4r)^2 ( {\sf h}_{\ast}^j + {\sf h}^{\ast}_j ) \, |B_{5 r}| +
	\frac{1}{\vartheta_j (4r)^2 {\sf h}_{\ast}^j} \,
	\uprho \big( B_{5r} (\bar{y}) \times (\bar{s} - \vartheta_j (4r)^2 {\sf h}_{\ast}^j, \bar{s} + \vartheta_j (4r)^2 {\sf h}^{\ast}_j) \big) \right] } \leqslant				\nonumber	\\
& \qquad {\displaystyle \leqslant \frac{\gamma}{2} \, \frac{\eta^2h^2}{2^{2m}} \, |B_{5 r}| \left[ 32 \left( 1 + r^2 \right) \vartheta_j ( {\sf h}_{\ast}^j + {\sf h}^{\ast}_j ) +
	\frac{25}{16} \frac{1}{\vartheta_j} 
	\frac{\uprho \big( B_{5r} (\bar{y}) \times (\bar{s} - \vartheta_j (4r)^2 {\sf h}_{\ast}, \bar{s} + \vartheta (4r)^2 {\sf h}^{\ast}) \big)}
		{\uprho \big( B_{5r} (\bar{y}) \times (\bar{s} - (5r)^2, \bar{s}) \big)} \right] }	.															\nonumber
\end{align}
Notice that by (C.5)
\begin{align*}
\frac{\uprho \big( B_{5r} (\bar{y}) \times (\bar{s} - \vartheta_j (4r)^2 {\sf h}_{\ast}^j, \bar{s} + \vartheta_j (4r)^2 {\sf h}^{\ast}_j) \big)}
		{\uprho \big( B_{5r} (\bar{y}) \times (\bar{s} - (5r)^2, \bar{s}) \big)} \leqslant 
		{\sf L} \, \frac{16}{25} \, \vartheta_j \, ( {\sf h}_{\ast}^j + {\sf h}^{\ast}_j ) ,
\end{align*}
From that and \eqref{servedopo}, \eqref{primopasso} and \eqref{secondopasso} we get
\begin{align*}
& a_{m+1}^j < 8 \, \frac{{\sf L}^3}{\upkappa^2} \, \upbeta \, c_{\uprho_j}^2 \sqrt{\frac{\gamma}{2}} \, \sqrt{|B_{5 r}|}
	\, \frac{\uprho_j \big( B_{4r} (\bar{y}) \times (\bar{s}, \bar{s} + \vartheta_j (4r)^2 {\sf h}^{\ast}_j) \big)}{\big( \vartheta_j (4r)^2 {\sf h}^{\ast}_j |B_{4r}|\big)^{\frac{1}{p}}} \, \cdot		\\
& \qquad \quad \cdot \big[ 32 \left( 1 + r^2 \right) \vartheta_j ( {\sf h}_{\ast}^j + {\sf h}^{\ast}_j ) +
	{\sf L} \, ( {\sf h}_{\ast}^j + {\sf h}^{\ast}_j ) \big]^{\frac{1}{2}} (b_{m-1} - b_m)^{\frac{2 - p}{2p}} .
\end{align*}
Now, summing up to a generic $m_o \in \N$ and since $\{ a_m^j\}_m$ is decreasing we get
\begin{align*}
m_o \, ( a_{m_o +1}^j)^{\frac{2p}{2 - p}} \leqslant \sum_{m = 1}^{m_o} (a_{m+1}^j)^{\frac{2p}{2 - p}} < {\sf Q} \, (b_0 - b_{m_o}) < {\sf Q} \, b_0
\end{align*}
where ${\sf Q}$ denotes
\begin{align*}
& \left[ 8 \, \frac{{\sf L}^3}{\upkappa^2} \, \upbeta \, c_{\uprho_j}^2 \sqrt{\frac{\gamma}{2}} \, \sqrt{|B_{5 r}|}
	\, \frac{\uprho_j \big( B_{4r} (\bar{y}) \times (\bar{s}, \bar{s} + 
	\vartheta_j (4r)^2 {\sf h}^{\ast}_j) \big)}{\big( \vartheta_j (4r)^2 {\sf h}^{\ast}_j |B_{4r}|\big)^{\frac{1}{p}}} \right]^{\frac{2p}{2 - p}} \cdot				\\
& \hskip80pt \cdot \big[ 32 \left( 1 + r^2 \right) \vartheta_j ( {\sf h}_{\ast}^j + {\sf h}^{\ast}_j ) +
	{\sf L} \, ( {\sf h}_{\ast}^j + {\sf h}^{\ast}_j ) \big]^{\frac{p}{2 - p}} .
\end{align*}
Since
\begin{gather*}
b_0 = \Big| \Big( B_{4 r} (\bar{y}) \times (\bar{s}, \bar{s} + \vartheta_j (4r)^2 {\sf h}^{\ast}_j) \Big) \cap \big\{ u < \eta h \big \} \Big| 
	< \vartheta_j (4r)^2 {\sf h}^{\ast}_j \big| B_{4 r}(\bar{y}) \big|	,		\\
\frac{|B_{5 r}|^{\frac{1}{2}}}{|B_{4 r}|^{\frac{1}{p}}} \leqslant \frac{|B_{8 r}|^{\frac{1}{2}}}{|B_{4 r}|^{\frac{1}{p}}} \leqslant \frac{\sqrt{2^n} |B_{4 r}|^{\frac{1}{2}}}{|B_{4 r}|^{\frac{1}{p}}}
	= \frac{\sqrt{2^n}}{|B_{4 r}|^{\frac{2-p}{2p}}}
\end{gather*}
and then
\begin{gather*}
\frac{|B_{5 r}|^{\frac{1}{2}}}{\big( \vartheta_j (4r)^2 {\sf h}^{\ast}_j |B_{4r}|\big)^{\frac{1}{p}}}
	\, b_0^{\frac{2 - p}{2p}} \leqslant \sqrt{2^n} \frac{\big( \vartheta_j (4r)^2 {\sf h}^{\ast}_j  \big)^{\frac{2 - p}{2p}}}{{\big( \vartheta_j (4r)^2 {\sf h}^{\ast}_j \big)^{\frac{1}{p}}}} =
	\frac{\sqrt{2^n}}{{\big( \vartheta_j (4r)^2 {\sf h}^{\ast}_j \big)^{\frac{1}{2}}}}
\end{gather*}
we finally get
\begin{align*}
a_{m_o+1}^j	& < {\sf C} \, \left( \frac{1}{m_o} \right)^{\frac{2 - p}{2p}} \uprho_j \big( B_{4r} (\bar{y}) \times (\bar{s}, \bar{s} + \vartheta_j (4r)^2 {\sf h}^{\ast}_j) \big) \leqslant			\\
			& \leqslant  {\sf C} \, \left( \frac{1}{m_o} \right)^{\frac{2 - p}{2p}} \upeta_j^{\varsigma^0} \big( B_{4r} (\bar{y}) \times (\bar{s}, \bar{s} + \vartheta_j (4r)^2 {\sf h}^{\ast}_j) \big)
\end{align*}
where
\begin{align*}
{\sf C} := 
8 \, \frac{{\sf L}^3}{\upkappa^2} \, \upbeta \, c_{\uprho_j}^2 \sqrt{\frac{2^n \, \gamma}{2}} \, 
	\sqrt{ \frac{32 \left( 1 + r^2 \right) \vartheta_j + {\sf L}}{\vartheta_j (4r)^2} \frac{{\sf h}_{\ast}^j + {\sf h}^{\ast}_j}{{\sf h}^{\ast}_j } } \ .
\end{align*}
Now for each $\epsilon > 0$ we can find $m_o$ for which ${\sf C}  \, \left( \frac{1}{m_o} \right)^{\frac{2 - p}{2p}} \leqslant \epsilon$, i.e.
$m_o \mau ( {\sf C} / \epsilon )^{\frac{2p}{2 - p}}$; 
therefore the first point of the thesis is proved
taking $\eta_o := \eta h /2^{m_o + 1}$. \\
To get the second point first define
\begin{align*}
\upchi_j^{\varsigma^0} = 1 \quad \text{in } \big( Q^{4r, \vartheta_j}_j (\bar{x}, \bar{t}) \big)^{\varsigma^0}, \qquad
\upchi_j^{\varsigma^0} = 0 \quad \text{in } Q^{4r, \vartheta_j} (\bar{x}, \bar{t}) \setminus \big( Q^{4r, \vartheta_j}_j (\bar{x}, \bar{t}) \big)^{\varsigma^0} .
\end{align*}
Then observe that, by (H.2), (C.10) and Remark \ref{rimi2}, for $S \subset B_{4r} (\bar{x})$:
\begin{align*}
& \left( \frac{|S|}{\upchi_j^{\varsigma^0} (\sigma) \big( B_{4r} (\bar{x}) \big)} \right)^b 
	\leqslant \left( \frac{1}{2 \upkappa} \frac{|S|}{\upchi_j^{\varsigma^0} (\sigma) \big( B_{4r} (\bar{x}) \big)} \right)^b \leqslant				\\
& \hskip50pt	
	\leqslant \left( \frac{1}{2 \upkappa} \right)^b K_ 2 \frac{ \uprho (\sigma) (S)}{\uprho (\sigma) \big( B_{4r} (\bar{x}) \big)}
	\leqslant \left( \frac{1}{2 \upkappa} \right)^b K_ 2 \frac{ \uprho (\sigma) (S)}{\upeta_j^{\varsigma^0} (\sigma) \big( B_{4r} (\bar{x}) \big)}
\end{align*}
by which and (C.10)
\begin{align*}
& \frac{\big| \big\{ u < \eta_o h \big \} \cap \big( Q^{4r, \vartheta_j}_j (\bar{x}, \bar{t}) \big)^{\varsigma^0} \big|^b}
	{\big| \big( Q^{4r, \vartheta_j}_j (\bar{x}, \bar{t}) \big)^{\varsigma^0} \big|^b} \leqslant 												\\
& \hskip30pt \leqslant \left( \frac{1}{2 \upkappa} \right)^b K_ 2 \, {\sf L} \, \frac{\uprho \Big( \big\{ u < \eta_o h \big \} \cap \big( Q^{4r, \vartheta_j}_j (\bar{x}, \bar{t}) \big)^{\varsigma^0} \Big)}
{\uprho \big( \big( Q^{4r, \vartheta_j}_j (\bar{x}, \bar{t}) \big)^{\varsigma^0} \big)} < \left( \frac{1}{2 \upkappa} \right)^b K_ 2 \, {\sf L} \, \epsilon
\end{align*}
and finally
$$
\displaylines{
\hfill
\big| \big\{ u < \eta_o h \big \} \cap \big( Q^{4r, \vartheta_j}_j (\bar{x}, \bar{t}) \big)^{\varsigma^0} \big| < \frac{1}{2 \upkappa} \, \big( K_ 2 \, {\sf L} \, \epsilon \big)^{1/b}
	\big| \big( Q^{4r, \vartheta_j}_j (\bar{x}, \bar{t}) \big)^{\varsigma^0} \big| \, .
\hfill\llap{$\square$}}
$$
\fine

\bthm
\label{esp_positivita}
Consider $j \in \{ 1, 2 \}$, $(x_o, t_o) \in \mathcal{Q}_j \cup \Upgamma$, $r > 0$ such that $B_{5 r} (x_o) \times (t_o - (5r)^2 {\sf h}_{\ast}^j, t_o + (5r)^2 {\sf h}^{\ast}_j) \subset \Omega \times (0,T)$ where
$$
{\sf h}_{\ast}^j = {\sf h}_{\ast}^j (x_o, t_o, 5r) , \qquad {\sf h}^{\ast}_j = {\sf h}^{\ast} (x_o, t_o, 5r)  .
$$
Suppose $5r \leqslant \bar{R}$ and $(4 r)^2 {\sf h}^{\ast}_j ( x_o, y_o, 5r) < \updelta$.
Let $\vartheta_j \in (0, 1)$ be the value determined in Lemma $\ref{lemma1}$ corresponding to $\theta = 1$.
Then for every $\hat\vartheta \in (0, \vartheta_j / 4)$ there is $\lambda_j \in (0,1)$ depending $(\!$only$)$ on
$$
\upgamma, \gamma, \kappa, r, \vartheta_j, \hat\vartheta, c_{\uprho_j}, c_{\upchi_j} , \upkappa, {\sf L}, {\sf L}^{-1}, \upbeta, \frac{{\sf h}^{\ast}_j}{{\sf h}_{\ast}^j + {\sf h}^{\ast}_j}
$$
such that for every $h > 0$ and $u \geqslant 0$ in $DG(\Omega, T, \uprho, \gamma)$ if
$$
u(\cdot , t_o) \geqslant h \qquad \text{ a.e. in }  B_{r}^j (x_o; t_o) .
$$
then
$$
u \geqslant \lambda_j h \qquad \text{a.e. in } \Big( B_{2r} (x_o) \times \big[t_o + \hat\vartheta (4r)^2 {\sf h}^{\ast}_j (x_o, t_o, 5r) , t_o + \vartheta_j (4r)^2 {\sf h}^{\ast}_j (x_o, t_o, 5r) \big]
						\Big) \cap \mathcal{Q}_j .
$$
\ethm
\boss
As observed after Lemma \ref{lemma2} and being $\lambda_j$ depending on the same quantities on which $\eta_o$ depends on,
we have that the more $\uprho_j (t_o) \left( B_{4r}(x_o) \right)$ is little the more $\vartheta_j$ and $\lambda_j$ are little.
\eoss
\boss
\label{garofano2}
Notice that the constants $c_{\uprho_j}$ $c_{\upchi_j}$ may be chosen not depending on $j$, as in fact we did (see (C.7)).
Moreover also the dependence of $\vartheta_j$ on $j$ may be dropped, as observed in Remark \ref{garofano}.
Then the only dependence of $\lambda_j$ on $j$ relies on $\frac{{\sf h}^{\ast}_j}{{\sf h}_{\ast}^j + {\sf h}^{\ast}_j}$.
Then taking
$$
\min \left\{ \frac{{\sf h}^{\ast}_1}{{\sf h}_{\ast}^1 + {\sf h}^{\ast}_1} , \frac{{\sf h}^{\ast}_2}{{\sf h}_{\ast}^2 + {\sf h}^{\ast}_2} \right\} ,
$$
or take $\lambda := \min \{ \lambda_1, \lambda_2 \}$, one can drop the dependence on $j$ and have the same value of $\lambda$ in $\mathcal{Q}_1$ and in $\mathcal{Q}_2$.
\eoss
\boss
Notice that by (C.3) one has that for each $t \in [t_o - \updelta, t_o + \updelta]$ (and for every $j \in \{ 1, 2 \}$ and $r \in (0, \bar{R}]$)
$\uprho_j (t) \big( B_r (x_o) \big) > 0$, i.e. $(\mathcal{Q}_j \cup \Upgamma) \cap \big( B_r (x_o) \times \{t\} \big) \not= \emptyset$ for every $t \in [t_o - \updelta, t_o + \updelta]$.
\eoss
\noindent
\dimo
In Proposition \ref{anninabella}
consider $\bar{x} = x_o$, $\bar{t} = t_o$, $\mu_- = 0$, $\sigma \omega = c$ ($c$ arbitrary, to be chosen),
$R = 5r$, $\tilde{s} = 4r$, $s = 2r$, $\tilde\vartheta = \vartheta_j$
where $\vartheta_j$ is the value determined in Lemma $\ref{lemma1}$ corresponding to $\theta = 1$
and $\hat\vartheta \in (0, \vartheta_j/4)$.
Consider $\varsigma^0 = (2r, (\vartheta_j - \hat\vartheta ) (4r)^2 {\sf h}^{\ast}_j (x_o, t_o, 5r) )$. \\
By Proposition \ref{anninabella}, point $\, iv \,)$, we have the existence of $\nu$ such that if
$$
\frac{\big|\{ (x,t) \in \big( Q^{4r, \vartheta_j}_{j} (\bar{x}, \bar{t}) \big)^{\varsigma^0} \, | \, u(x,t) < c \} \big|}
	{\big| \big( Q^{4r, \vartheta_j}_{j} (\bar{x}, \bar{t}) \big)^{\varsigma^0} \big|} +
\frac{\uprho \big( \{ (x,t) \in \big( Q^{4r, \vartheta_j}_{j} (\bar{x}, \bar{t}) \big)^{\varsigma^0} \, | \, u(x,t) < c \} \big)}
	{\uprho \big( \big( Q^{4r, \vartheta_j}_{j} (\bar{x}, \bar{t}) \big)^{\varsigma^0} \big)} \leqslant \nu
$$
then
\begin{align*}
u(x,t) \geqslant \frac{1}{2} \, c
\end{align*}
for a.e. $(x,t) \in \Big(B_{2r} (\bar{x}) \times [ t_o + \hat\vartheta (4r)^2 {\sf h}^{\ast}_j (x_o, t_o, 5r) ,
t_o + \vartheta_j (4r)^2 {\sf h}^{\ast}_j ( x_o, t_o, 5r) ] \Big) \cap \mathcal{Q}_j$. \\
Now we consider $\epsilon > 0$: we have then the existence of $\eta_o$ such that Lemma \ref{lemma2} holds. \\
Then taking $\epsilon$ such that ($b$ is the constant in Remark \ref{rimi2} and (C.10))
$$
\epsilon + \frac{(K_2 \, {\sf L} \, \epsilon )^{1/b}}{2 \upkappa} = \nu ,
$$
taking $c = \eta_o h$ we conclude taking $\lambda_j = \eta_o/2$ (corresponding to $a = 1/2$ in Proposition \ref{anninabella}). \\
Notice that the constant $\lambda_j$ depends on the same constants on which $\eta_o$ depends and, since $\eta_o$ depends also on $\epsilon$ which depends on $\nu$, 
consequently $\lambda_j$ depends 
on $K_2, b$ and on the constants on which $\nu$ and $\eta_o$ depend.
\finedimo

\section{The Harnack inequality}

Before stating the result we recall that ${\sf h}^{\ast}_j$ and ${\sf h}_{\ast}^j$ are defined in \eqref{funzioneacca}.

\bthm
\label{harnack}
Consider $j \in \{ 1, 2 \}$, $(x_o, t_o) \in \mathcal{Q}_j \cup \Upgamma$.
Consider $r > 0$ such that
$B_{5 r} (x_o) \times (t_o - (5r)^2 {\sf h}_{\ast}^j (x_o, t_o, 5r), t_o + (5 r)^2 {\sf h}^{\ast}_j (x_o, t_o, 5r)) \subset \Omega \times (0,T)$ and 
$B_{r} (x_o) \times (t_o , t_o + 2 \, r^2 {\sf h}^{\ast}(x_o, t_o, r)) \subset \Omega \times (0,T)$,
$$
5r \leqslant \bar{R} , \qquad \sup \, \{ (4 \varrho)^2 {\sf h}^{\ast}_j ( x, t, 5 \varrho) \, | \, 
	(x,t) \in Q \cap \mathcal{Q}_j , \varrho \leqslant r\}  < \updelta ,
$$
where $Q := B_{4 r} (x_o) \times (t_o - r^2 {\sf h}_{\ast}^j (x_o, t_o, r), t_o + r^2 {\sf h}^{\ast}_j (x_o, t_o, r))$.
Then there exists ${\sf c}_j$ such that 
for every $u \in DG (\Omega, T, \uprho, \gamma)$, $u \geqslant 0$, for every $(x_o, t_o) \in \mathcal{Q}_j \cup \Upgamma$ one has that
$$
u(x_o, t_o) \leqslant {\sf c}_j \, \inf_{B_{r}^j (x_o; t_o + r^2 {\sf h}^{\ast}_j (x_o, t_o, r))} u(x , t_o + r^2 {\sf h}^{\ast}_j (x_o, t_o, r)) .
$$
The constant ${\sf c}_j$ depends $($only$)$ on 
$\upgamma, \gamma, \kappa, r, \vartheta_j, c_{\uprho_j}, c_{\upchi_j} , \upkappa, {\sf L}, {\sf L}^{-1}, \upbeta, \frac{{\sf h}^{\ast}_j}{{\sf h}_{\ast}^j + {\sf h}^{\ast}_j}$.
\ethm

\boss
\label{scimmiettabella}
Notice that, in particular, ${\sf c}_j$
depends on $\frac{{\sf h}^{\ast}_j}{{\sf h}_{\ast}^j + {\sf h}^{\ast}_j}$.
\eoss

\boss
\label{maledetta_cervicale}
Following the proof one can realise that the constant ${\sf c}_j$ ($j = 1,2$) increases with increasing  $r$. For this reason, once fixed $r > 0$ as in the theorem above
and under the same assumptions of the theorem, one can consider ${\sf c}_j$ such that
$$
u(x_o, t_o) \leqslant {\sf c}_j \, \inf_{B_{\rho}^j (x_o; t_o + \rho^2 {\sf h}^{\ast}_j (x_o, t_o, r))} u(x , t_o + \rho^2 {\sf h}^{\ast}_j (x_o, t_o, r)) 
\qquad \text{for every } \rho \in (0, r]
$$
and then
$$
u(x_o, t_o) \leqslant {\sf c}_j \, \inf_{P_{r,  {\sf h}^{\ast}_j (x_o, t_o, r)} (x_o, t_o)} u
$$
where $P_{r,  {\sf h}^{\ast}_j (x_o, t_o, r)} (x_o, t_o)$, shown in the picture below, is
\begin{align}
\label{parabola}
P_{r,  {\sf h}^{\ast}_j (x_o, t_o, r)} (x_o, t_o) := \bigcup_{\rho \in (0, r]} B_{\rho}^j (x_o; t_o + \rho^2 {\sf h}^{\ast}_j (x_o, t_o, r)) .
\end{align}
\ \\
\ \\
\begin{center}
\begin{tikzpicture}
\begin{axis} [axis equal, xtick={1.5}, ytick={-1.2}, xticklabels={$x_o$}, yticklabels={$t_o$}, width=7cm, height=6cm, title={Figure 6}]
\addplot coordinates
{(1.5, -1.2)};
\addplot
[domain=0:3,variable=\t, smooth]
({t},{-1.2});
\addplot
[domain=3:4,variable=\t, smooth]
({t},{-1.2});
\addplot
[domain=-1.2:1.2,variable=\t, smooth]
({0}, {t});
\addplot
[domain=-1.2:1.2,variable=\t, smooth]
({4}, {t});
\addplot
[domain=0:3.85,variable=\t, smooth]
({t},{1.2});
\addplot
[domain=2.85:4,variable=\t, smooth]
({t},{1.2});
\addplot
[domain=1.2:1.6,variable=\t,smooth]
({0.6*t^(1/3)+2-0.64 + 0.85}, {t});
\addplot
[domain=0.05:1.2,variable=\t, smooth]
({0.6*t^(1/3)+1.36 + 0.85}, {t});
\addplot
[domain=0.05:0.6,variable=\t, smooth]
({-0.6*t^(1/3)+1.36 + 0.45 + 0.85}, {-t + 0.11});
\addplot
[domain=0.6:1.65,variable=\t, smooth]
({-0.6*t^(1/3)+1.36 + 0.45 + 0.85}, {-t + 0.11});
\addplot
[domain=0.41:2.09,variable=\t, smooth, very thin]
({t}, {1.2*(t-1.5)^2-1.2});
\addplot
[domain=0.42:2.55,variable=\t, very thick]
({t},{0.2});
%
\addplot
[domain=0.47:0.59,variable=\t, smooth]
({t},{t - 0.4});
\addplot
[domain=0.55:0.89,variable=\t, smooth]
({t},{t - 0.7});
\addplot
[domain=0.65:1.19,variable=\t, smooth]
({t},{t - 1});
\addplot
[domain=0.75:1.49,variable=\t, smooth]
({t},{t - 1.3});
\addplot
[domain=0.87:1.79,variable=\t, smooth]
({t},{t - 1.6});
\addplot
[domain=1:2.09,variable=\t, smooth]
({t},{t - 1.9});
\addplot
[domain=1.15:2.39,variable=\t, smooth]
({t},{t - 2.2});
\addplot
[domain=1.34:2.21,variable=\t, smooth]
({t},{t - 2.5});
\addplot
[domain=1.6:2.1,variable=\t, smooth]
({t},{t - 2.8});
\end{axis}
\end{tikzpicture}
\end{center}
\ \\
\eoss

\boss
\label{portugal}
Notice that if $(x_o, t_o) \in \Upgamma$ one has
$$
u(x_o, t_o) \leqslant {\sf c}_j \, \inf_{B_r^j (x_o; t_o + \rho^2 {\sf h}^{\ast}_j (x_o, t_o, r))} u(x , t_o + r^2 {\sf h}^{\ast}_j (x_o, t_o, r)) 
$$
both for $j=1$ and $j=2$. \\
Since $c_{\uprho_j}, c_{\upchi_j}$ can be chosen to be independent of $j$ (see (C.7))
the dependence on $j$ of ${\sf c}_j$ is confined to the ratio $\frac{{\sf h}^{\ast}_j}{{\sf h}_{\ast}^j + {\sf h}^{\ast}_j}$.
Clearly one can drop also the dependence on this term since
$$
\frac{\min \{ {\sf h}^{\ast}_1, {\sf h}^{\ast}_2 \}}{{\sf h}_{\ast}^j + {\sf h}^{\ast}_j} \leqslant \frac{{\sf h}^{\ast}_j}{{\sf h}_{\ast}^j + {\sf h}^{\ast}_j}
	 \leqslant \frac{\max \{ {\sf h}^{\ast}_1, {\sf h}^{\ast}_2 \}}{{\sf h}_{\ast}^j + {\sf h}^{\ast}_j} .
$$
\eoss

%

\noindent
With the same assumptions of Theorem \ref{harnack} with moreover $B_{r} (x_o) \times (t_o - 2 \, r^2 {\sf h}^{\ast}(x_o, t_o, r), t_o) \subset \Omega \times (0,T)$,
similarly one can prove the existence of $\tilde{\sf c}_j$ depending on
$\upgamma, \gamma, \kappa, r, \vartheta_j, c_{\uprho_j}, c_{\upchi_j} , \upkappa, {\sf L}, {\sf L}^{-1}, \upbeta$, $\frac{{\sf h}_{\ast}^j}{{\sf h}_{\ast}^j + {\sf h}^{\ast}_j}$
such that
\begin{align}
\label{cruccodim2}
\tilde{\sf c}_j \, \sup_{B_{r} (x_o)} u(x , t_o - r^2 {\sf h}_{\ast}^j(x_o, t_o, r)) \leqslant u(x_o, t_o)
\end{align}
and, as a consequence of this and Theorem \ref{harnack}, that
\begin{align}
\label{cruccodim}
\sup_{B_{r} (x_o)} u(x , t_o - r^2 {\sf h}_{\ast}^j (x_o, t_o, r)) \leqslant \hat{\sf c}_j \, \inf_{B_r^j (x_o; t_o + r^2 {\sf h}^{\ast}_j (x_o, t_o, r))} u(x , t_o + r^2 {\sf h}^{\ast}_j (x_o, t_o, r)) .
\end{align}

\noindent
In Figure 7 we show two examples, in the first we suppose that ${\sf h}^{\ast}_1 (x_o, t_o, r) > {\sf h}^{\ast}_2 (x_o, t_o, r)$, in the second
${\sf h}^{\ast}_1 (x_o, t_o, r) < {\sf h}^{\ast}_2 (x_o, t_o, r)$,
where $\mathcal{Q}_1$ is on the left hand side of the interface and $\mathcal{Q}_2$ is on the right one, the marked dot is $(x_o, t_o)$ and
we control the value of $u$ in $(x_o, t_o)$ by the values of $u$ in the {\em bold} set. \\
Notice that the two sets
$B_r^1 (x_o; t_o + r^2 {\sf h}^{\ast}_1 (x_o, t_o, r))$ and $B_r^2 (x_o; t_o + r^2 {\sf h}^{\ast}_2 (x_o, t_o, r))$ could have empty or not empty intersection. \\
In the following examples $(x_o, t_o) \in \Upgamma$, but it could belong also to $\mathcal{Q}_j$.
\ \\
\ \\
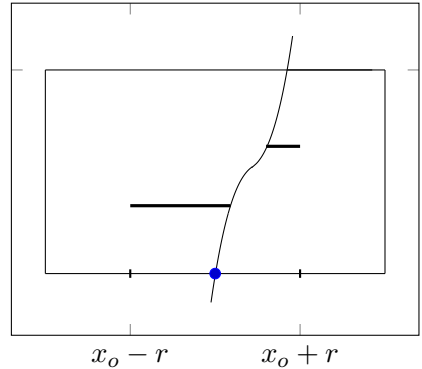
\begin{center}
\begin{tikzpicture}
\begin{axis} [axis equal, xtick={1,3}, ytick={1.2}, xticklabels={$x_o-r$, $x_o+r$}, yticklabels={$$}, width=7cm, height=6cm, title={Figure 7.a}]
\addplot coordinates
{(2, -1.2)};
\addplot
[domain=-1.25:-1.15,variable=\t, thick]
({3},{t});
\addplot
[domain=-1.25:-1.15,variable=\t, thick]
({1},{t});
\addplot
[domain=0:3,variable=\t, smooth]
({t},{-1.2});
\addplot
[domain=3:4,variable=\t, smooth]
({t},{-1.2});
\addplot
[domain=-1.2:1.2,variable=\t, smooth]
({0}, {t});
\addplot
[domain=-1.2:1.2,variable=\t, smooth]
({4}, {t});
\addplot
[domain=0:3.85,variable=\t, smooth]
({t},{1.2});
\addplot
[domain=2.85:4,variable=\t, smooth]
({t},{1.2});
\addplot
[domain=1.2:1.6,variable=\t,smooth]
({0.6*t^(1/3)+2-0.64 + 0.85}, {t});
\addplot
[domain=0.05:1.2,variable=\t, smooth]
({0.6*t^(1/3)+1.36 + 0.85}, {t});
\addplot
[domain=0.05:0.6,variable=\t, smooth]
({-0.6*t^(1/3)+1.36 + 0.45 + 0.85}, {-t + 0.11});
\addplot
[domain=0.6:1.65,variable=\t, smooth]
({-0.6*t^(1/3)+1.36 + 0.45 + 0.85}, {-t + 0.11});
\addplot
[domain=1:2.55,variable=\t, very thick]
({t},{0.2});
\addplot
[domain=2.15:3,variable=\t, very thick]
({t},{-0.5});
\end{axis}
\end{tikzpicture}
\qquad\quad
\begin{tikzpicture}
\begin{axis} [axis equal, xtick={1,3}, ytick={1.2}, xticklabels={$x_o-r$, $x_o+r$}, yticklabels={$$}, width=7cm, height=6cm, title={Figure 7.b}]
\addplot coordinates
{(2, -1.2)};
\addplot
[domain=-1.25:-1.15,variable=\t, thick]
({3},{t});
\addplot
[domain=-1.25:-1.15,variable=\t, thick]
({1},{t});
\addplot
[domain=0:3,variable=\t, smooth]
({t},{-1.2});
\addplot
[domain=3:4,variable=\t, smooth]
({t},{-1.2});
\addplot
[domain=-1.2:1.2,variable=\t, smooth]
({0}, {t});
\addplot
[domain=-1.2:1.2,variable=\t, smooth]
({4}, {t});
\addplot
[domain=0:3.85,variable=\t, smooth]
({t},{1.2});
\addplot
[domain=2.85:4,variable=\t, smooth]
({t},{1.2});
\addplot
[domain=1.2:1.6,variable=\t,smooth]
({0.6*t^(1/3)+2-0.64 + 0.85}, {t});
\addplot
[domain=0.05:1.2,variable=\t, smooth]
({0.6*t^(1/3)+1.36 + 0.85}, {t});
\addplot
[domain=0.05:0.6,variable=\t, smooth]
({-0.6*t^(1/3)+1.36 + 0.45 + 0.85}, {-t + 0.11});
\addplot
[domain=0.6:1.65,variable=\t, smooth]
({-0.6*t^(1/3)+1.36 + 0.45 + 0.85}, {-t + 0.11});
\addplot
[domain=1:2.18,variable=\t, very thick]
({t},{-0.4});
\addplot
[domain=2.6:3,variable=\t, very thick]
({t},{0.3});
\end{axis}
\end{tikzpicture}
\end{center}

\ \\

\noindent
Before going on we state a lemma needed to prove Theorem \ref{harnack}. The proof can be obtained adapting the proof of Lemma 2.14 in \cite{fabio10}
and for this reason we do not report it. Before we define the following set: for $(x_o, t_o) , (\upsilon, \sigma) \in \Omega \times (0,T)$ we define
\begin{align}
\label{luce}
C_{s}^j (\upsilon, \sigma) = \big( B_{s} (\upsilon) \times \left( \sigma - s^2 {\sf h}^{\ast}_j (x_o, 0, r) , \sigma \right) \big) \cap \big( \mathcal{Q}_j \cup \Upgamma \big) .
\end{align}

\begin{lemma}
\label{lemmaMisVar}
Consider $j \in \{1, 2 \}$, $(x_o, t_o) \in \mathcal{Q}_j \cup \Upgamma$, $r > 0$ such that $B_{2 r}(x_o) \subset \Omega$, 
$A, B , c > 0$. 
Consider $\sigma_o = ( r/2, c \, r^2 /4)$ and $\sigma_1 = ( r, c \, r^2)$.
Consider $\tilde{t} \in [t_o - \beta r^2, t_o]$ for some positive $\beta$.
Then, for every $a > 0$, $\varepsilon , \delta\in (0,1)$ there exists $\eta \in (0,1)$ such that for
every $v$ Lipschitz continuous function defined in the ball $B_{2 r} (x_o)$ satisfying
$$
\int_{(C_{2 r}^j (x_o, t_o) )^{\sigma_1} (\tilde{t}) } |Dv|^2 \, dx \leqslant A \, \frac{\big| \big(C_{2 r}^j (x_o, t_o) \big)^{\sigma_1} (\tilde{t}) \big|}{r^2},
$$
and
$$
\frac{\uprho (\tilde{t}) \big( \big\{ x \in \big(C_{r}^j (x_o, t_o) \big)^{\sigma_o} (\tilde{t}) \, \big| \, v (x) > a \big\} \big)}
	{\uprho (\tilde{t}) \big( \big( \big\{ x \in \big(C_{r}^j (x_o, t_o) \big)^{\sigma_o} (\tilde{t}) \big) } > B ,
$$
there exists $x^* \in \big(C_{r}^j (x_o, t_o) \big)^{\sigma_o} (\tilde{t})$ with $B_{\eta r}(x^*) \subset \big(C_{r}^j (x_o, t_o) \big)^{\sigma_o} (\tilde{t})$ such that
$$
\uprho (\tilde{t}) (\{ v > \varepsilon a \} \cap B_{\eta r} (x^*) ) > (1-\delta) \, \uprho (\tilde{t}) (B_{\eta r}(x^*) ) .
$$
\end{lemma}

\noindent
\dimoH
For the sake of simplicity we suppose $t_o = 0$, which is always possible up to a translation. 
Moreover we fix $j \in \{ 1, 2 \}$. \\
We may write $u(x_o, 0) = b \, r^{-\xi}$ for some $b, \xi > 0$ to be fixed later.
Define the functions
$$
\emm ( s ) = \sup_{C_{s}^j (x_o, 0)} u,		\qquad 	\enn (s) = b (r - s)^{-\xi}, 	\qquad s \in [0,r) ,
$$
where $C_{s}^j (\upsilon, \sigma)$ is defined in \eqref{luce}.
Let us denote by $s_o \in [0, r)$ the largest solution of $\emm (s) = \enn (s)$ (notice that $0$ is a solution). Define
$$
M := \enn (s_o) = b (r - s_o)^{-\xi} \, .
$$
We can find $(y_o, \tau_o) \in C_{s_o}^j (x_o, 0)$ such that 
\begin{equation}
\label{choicey0t0}
\frac{3}{4} M < \sup_{C_{\frac{r_o}{4}}^j (y_o, \tau_o)} u \leqslant M
\end{equation}
where $r_o \leqslant (r - s_o) / 2$: in this way we get
$\frac{r + s_o}{2} < r$
and
$$
C_{r_o}^j (y_o, \tau_o) \subset C_{\frac{r + s_o}{2}}^j (x_o, 0) \subset C_{r}^j (x_o, 0)
$$
and therefore (since $\frac{r + s_o}{2} > s_o$)
\begin{equation}
\label{bansky}
\sup_{C_{r_o}^j (y_o, \tau_o)} u \leqslant \sup_{C_{\frac{r + s_o}{2}}^j (x_o, 0)} u = \emm \left( \frac{r + s_o}{2} \right) < \enn \left(\frac{r + s_o}{2} \right) = 2^\xi M.
\end{equation}
We now proceed dividing the proof in seven steps. \\ [0.3em]
\textsl{\underline{Step 1} - }
In this step we want to show that there is $\bar{\nu} \in (0,1)$ such that
\begin{equation}
\label{mimancanogliamorimiei}
\begin{array}{c}
\uprho \left( \left\{ u > \frac{M}{2} \right\} \cap \big(C_{r_o/2}^j (y_o, \tau_o) \big)^{\varsigma_0} \right) 
	> \bar{\nu} \, \uprho \left( \big(C_{r_o/2}^j (y_o, \tau_o) \big)^{\varsigma_0} \right)
\end{array}
\end{equation}
where
\begin{align}
\label{sigmazero}
\varsigma_0 = \left( \frac{r_o}{4} , \frac{3 \, r_o^2}{16} \, {\sf h}_{\ast}^j (y_o, \tau_o, r) \right) .
\end{align}
To prove that, we prove first that there is $\bar{\bar{\nu}} \in (0,1)$
\begin{align}
\label{mimancanogliamorimiei_bis}
\frac{\uprho \left( \left\{ u > \frac{M}{2} \right\} \cap \big(C_{r_o/2}^j (y_o, \tau_o) \big)^{\varsigma_0} \right)}{\uprho \left( \big(C_{r_o/2}^j (y_o, \tau_o) \big)^{\varsigma_0} \right)} 
	+ \frac{\left| \left\{ u > \frac{M}{2} \right\} \cap \big(C_{r_o/2}^j (y_o, \tau_o) \big)^{\varsigma_0} \right|}{\left| \big(C_{r_o/2}^j (y_o, \tau_o) \big)^{\varsigma_0} \right|} > \bar{\bar{\nu}} \, .
\end{align}
Assume, by contradiction, that this is not true. Taking in Proposition \ref{carlettobello}
\begin{gather*}
(\bar{x}, \bar{t}) = (y_o, \tau_o) ,\quad s = \frac{r_o}{4}, \quad \tilde{s} = \frac{r_o}{2},  \quad R = r		,						\\
\theta = \tilde\theta \quad \text{s.t.} \quad \tilde\theta \, {\sf h}_{\ast}^j (y_o, \tau_o, r) = {\sf h}^{\ast}_j (x_o, 0, r)
	\quad \text{and} \quad  \theta \, {\sf h}_{\ast}^j (y_o, \tau_o, r) = {\sf h}^{\ast}_j (x_o, 0, r),									\\
\mu = \omega = 2^\xi M , \quad \sigma=1-2^{-\xi-1} \ \textrm{ and } \quad a = \sigma^{-1} \bigg( 1 - \frac{3}{2^{\xi+2}} \bigg) ,			\\
\end{gather*}
we obtain that
$$
u \leqslant \frac{3M}{4} \quad \textrm{in }\, C_{\frac{r_o}{4}}^j (y_o, \tau_o)
$$
which contradicts \eqref{choicey0t0}.
Notice that, according to Proposition  \ref{carlettobello}, $\bar{\bar{\nu}}$ depends on
$\upgamma, \gamma, \xi, \kappa, c_{\uprho_j} (r/r_o)$, $c_{\upchi_j} (r/r_o), {\sf h}^{\ast}_j (x_o, 0, r)/{\sf h}_{\ast}^j (y_o, \tau_o, r)$, $\upkappa, {\sf L}^{-1}$.
By (C.10), (C.4), \eqref{succofresco} (the constant $b$ is introduced in Remark \ref{rimi2}) we get
\begin{align}
\label{anninavainquarta}
\left( \frac{\left| \left\{ u > \frac{M}{2} \right\} \cap \big(C_{r_o/2}^j (y_o, \tau_o) \big)^{\varsigma_0} \right|}{\left| \big(C_{r_o/2}^j (y_o, \tau_o) \big)^{\varsigma_0} \right|} \right)^b
	\leqslant {\sf L} \, K_2 \, \frac{1 + \upkappa}{\upkappa^3} \, 
	\frac{\uprho \left( \left\{ u > \frac{M}{2} \right\} \cap \big(C_{r_o/2}^j (y_o, \tau_o) \big)^{\varsigma_0} \right)}{\uprho \left( \big(C_{r_o/2}^j (y_o, \tau_o) \big)^{\varsigma_0} \right)} .
\end{align}
Now, by \eqref{mimancanogliamorimiei_bis}, we derive that at least one of the two addends is greater than $\bar{\bar{\nu}}/2$. If
$$
\frac{\uprho \left( \left\{ u > \frac{M}{2} \right\} \cap \big(C_{r_o/2}^j (y_o, \tau_o) \big)^{\varsigma_0} \right)}{\uprho \left( \big(C_{r_o/2}^j (y_o, \tau_o) \big)^{\varsigma_0} \right)} > \frac{\bar{\bar{\nu}}}{2}
$$
we are done and $\bar{\nu} := \bar{\bar{\nu}}/2$. 
If, on the contrary, we have
$$
\frac{\left| \left\{ u > \frac{M}{2} \right\} \cap \big(C_{r_o/2}^j (y_o, \tau_o) \big)^{\varsigma_0} \right|}{\left| \big(C_{r_o/2}^j (y_o, \tau_o) \big)^{\varsigma_0} \right|} > \frac{\bar{\bar{\nu}}}{2}
$$
using \eqref{anninavainquarta} we consider
$$
\bar{\nu} := \left( \frac{\bar{\bar{\nu}}}{2} \right)^b \, \frac{\upkappa^3}{(1 + \upkappa) \, {\sf L} \, K_2} .
$$
\textsl{\underline{Step 2} - } In this step we want to prove that
\begin{equation}
\label{mimancanogliamorimiei-1}
\iint_{\big(C_{r_o}^j (y_o, \tau_o) \big)^{\varsigma_1}} |Du|^2 \, dx dt   \leqslant C \, (2^{\xi} M)^2 \, 
\min_{t \in [ \tau_o - r_o^2 \, {\sf h}^{\ast}_j (x_o, 0, r), \tau_o ]} \big| \big( C_{2 r_o}^j (y_o, \tau_o) \big)^{\varsigma_2} (t) \big|
\end{equation}
where $C = \frac{\gamma}{2} \, \frac{(1 + \upkappa) {\sf L}}{\upkappa^2} 
	\Big[ ( 4 + r_o^2) \, {\sf h}^{\ast}_j (x_o, 0, r) + \frac{4}{3} \, \frac{\uprho ( ( C_{r_o}^j (y_o, \tau_o) )^{\varsigma_1} )}
	{r_o^2 \, {\sf h}^{\ast}_j (x_o, 0, r) \max_{t \in [ \tau_o - r_o^2 \, {\sf h}^{\ast}_j (x_o, 0, r), \tau_o ]} | ( C_{r_o}^j (y_o, \tau_o) )^{\varsigma_1} (t) |} \Big]$. \\
and, once considered $\varsigma_0$ as above,
$$
\varsigma_1 = \left( \frac{r_o}{2} , \frac{3 \, r_o^2}{4} \, {\sf h}_{\ast}^j (y_o, \tau_o, r) \right) , \qquad
\varsigma_2 = \left( r_o , 3 r_o^2 \, {\sf h}_{\ast}^j (y_o, \tau_o, r) \right) .
$$
Consider a function $\zeta$ defined in $\Omega \times (0, \tau_o]$ such that
\begin{gather*}
\zeta \equiv 1 \hskip10pt \text{in } \big( C_{r_o}^j (y_o, \tau_o) \big)^{\varsigma_1} ,			\\	
\zeta \equiv 0 \hskip10pt \text{ outside of }  \big( C_{2 r_o}^j (y_o, \tau_o) \big)^{\varsigma_2} ,
	\qquad \text{for } t \leqslant \tau_o , 	\\
0 \leqslant \zeta \leqslant 1 \, , \qquad {\displaystyle |D \zeta| \leqslant \frac{1}{r_o} } \, , \qquad 
	{\displaystyle 0 \leqslant \zeta_t \leqslant \frac{1}{3} \, \frac{1}{r_o^2} \, \frac{1}{{\sf h}^{\ast}_j (x_o, 0, r)} }   \, .
\end{gather*}
Using this function $\zeta$ in \eqref{lavagna_bis}, taking $k = 0$, since $u \leqslant 2^{\xi} M$ in $C_{r_o}^j (y_o, \tau_o)$ (see \eqref{bansky}) 
and using \eqref{succofresco_0} we get
\begin{align*}
& \iint_{( C_{r_o}^j (y_o, \tau_o) )^{\varsigma_1}} |Du|^2 \, dx dt \leqslant 																				\\
& \hskip30pt \leqslant \frac{\gamma}{2} \iint_{( C_{2 r_o}^j (y_o, \tau_o) )^{\varsigma_2}} u^2 \big( |D \zeta|^2 + \zeta \zeta_t \uprho + \zeta^2 \big) dx dt		\leqslant		\\
& \hskip30pt \leqslant \frac{\gamma}{2} \bigg[ \frac{1}{r_o^2} (2^{\xi} M)^2 \bigg( 4 \big| \big( C_{2 r_o}^j (y_o, \tau_o) \big)^{\varsigma_2} \big| +
		\frac{4}{3} \, \frac{\uprho \big(\big( C_{2 r_o}^j (y_o, \tau_o) \big)^{\varsigma_2} \big)}{{\sf h}^{\ast}_j (x_o, 0, r)} \bigg) +
		(2^{\xi} M)^2 \big| \big( C_{2 r_o}^j (y_o, \tau_o) \big)^{\varsigma_2} \big| \bigg] \leqslant															\\
& \hskip30pt  \leqslant \frac{\gamma}{2} \, (2^{\xi} M)^2 \, \frac{(1 + \upkappa) {\sf L}}{\upkappa^2} 
	\min_{t \in [ \tau_o - 4 r_o^2 \, {\sf h}^{\ast}_j (x_o, 0, r), \tau_o ]} \big| \big( C_{2 r_o}^j (y_o, \tau_o) \big)^{\varsigma_2} (t) \big|
	\bigg[ ( 4 + r_o^2) \, {\sf h}^{\ast}_j (x_o, 0, r) \, + 																							\\
& \hskip70pt + \frac{4}{3} \, \frac{\uprho \big( \big( C_{2 r_o}^j (y_o, \tau_o) \big)^{\varsigma_2} \big)}
	{r_o^2 \, {\sf h}^{\ast}_j (x_o, 0, r) \max_{t \in [ \tau_o - 4 r_o^2 \, {\sf h}^{\ast}_j (x_o, 0, r), \tau_o ]} \big| \big( C_{2 r_o}^j (y_o, \tau_o) \big)^{\varsigma_2} (t) \big|} \bigg] .
\end{align*}
\\ [0.3em]
\textsl{\underline{Step 3} - }
The goal of this step is to show the existence of
$\bar{t} \in \big[\tau_o - \frac{r_o^2}{4} \, {\sf h}^{\ast}_j (x_o, 0, r) , \tau_o \big]$ such that
\begin{equation}
\label{carlettomio}
\begin{array}{c}
{\displaystyle
\frac{\uprho (\bar{t}) \big( \big\{ x \in \big(C_{r_o/2}^j (y_o, \tau_o) \big)^{\varsigma_0} (\bar{t}) \, \big| \, u(x,\bar{t}) > \frac{M}{2} \big\} \big)}
	{\uprho (\bar{t}) \big( \big(C_{r_o/2}^j (y_o, \tau_o) \big)^{\varsigma_0} (\bar{t}) \big) } > \frac{\bar\nu}{2} \, \frac{\upkappa^2}{(1 + \upkappa) {\sf L}}	}		\\	[2em]
{\displaystyle
\int_{\big(C_{r_o}^j (y_o, \tau_o) \big)^{\varsigma_1} (\bar{t})}  | Du(x,\bar{t}) |^2 \, dx \leqslant \frac{\alpha \, (2^{\xi} M)^2}{r_o^2} \, \min_{t \in [ \tau_o - r_o^2 \, {\sf h}^{\ast}_j (x_o, 0, r), \tau_o ]} 
	\big| \big(C_{r_o}^j (y_o, \tau_o) \big)^{\varsigma_1} (t) \big| \, , }
\end{array}
\end{equation}
where
\begin{equation}
\label{alfa}
\begin{array}{l}
{\displaystyle
\alpha = \frac{2^{n+ 2} \gamma}{\bar\nu} \, \frac{(1 + \upkappa)^2 {\sf L}^2}{\upkappa^5} \, \frac{1}{{\sf h}^{\ast}_j (x_o, 0, r)}
	\bigg[ ( 4 + r_o^2) \, {\sf h}^{\ast}_j (x_o, 0, r) +	}																\\	[1.5em]
{\displaystyle
\hskip40pt  + \frac{4}{3} \, \frac{\uprho ( ( C_{2 r_o}^j (y_o, \tau_o) )^{\varsigma_2} )}
	{r_o^2 \, {\sf h}^{\ast}_j (x_o, 0, r) \max_{t \in [ \tau_o - 4 r_o^2 \, {\sf h}^{\ast}_j (x_o, 0, r), \tau_o ]} | ( C_{2 r_o}^j (y_o, \tau_o) )^{\varsigma_2} (t) |} \bigg].   }
\end{array}
\end{equation}
and $\bar\nu$ has been determined in \textsl{Step 1}.
To do that we define the following sets (see \eqref{wlemedie} for the definition of $\big(C_{r_o/2}^j (y_o, \tau_o) \big)^{\varsigma_0} (t)$)
$$
\begin{array}{c}
A(t) = \big\{ x \in \big(C_{r_o/2}^j (y_o, \tau_o) \big)^{\varsigma_0} (t) \, \big| \, u(x,t) \geqslant \frac{M}{2} \big\}						\\	[4mm]
{\displaystyle
I = \bigg\{ t \in \big[ \tau_o - \tau_1 ,\tau_o \big] \, \Big| \, \uprho (t) ( A(t) ) > 
		\frac{\bar\nu}{2} \, \frac{\upkappa^2}{(1 + \upkappa) {\sf L}} \,  \uprho (t) \big( \big(C_{r_o/2}^j (y_o, \tau_o) \big)^{\varsigma_0} (t) \big) \bigg\},		}		\\	[4mm]
J_\alpha =  \bigg\{ t \in \big[ \tau_o - r_o^2 \, {\sf h}^{\ast}_j (x_o, 0, r), \tau_o \big] \, \Big| \,
{\displaystyle \int_{\big(C_{r_o}^j (y_o, \tau_o) \big)^{\varsigma_1} (t)}  | Du(x,t) |^2 \, dx \leqslant 		\hskip80pt	}										\\	[4mm]
	\hskip170pt	{\displaystyle \leqslant \frac{\alpha \, (2^{\xi} M)^2}{r_o^2} \, \min_{t \in [ \tau_o - \tau_1, \tau_o ]} 
	\big| \big(C_{r_o}^j (y_o, \tau_o) \big)^{\varsigma_1} (t) \big| } \bigg\} ,
\end{array}
$$
with $\alpha > 0$ and where, for the sake of simplicity, we denote by $\tau_1$ the quantity
$$
\tau_1 := \frac{r_o^2}{4} \, {\sf h}^{\ast}_j (x_o, 0, r) .
$$
By \eqref{mimancanogliamorimiei} and \eqref{succofresco_0} we get
\begin{align*}
\bar{\nu} \min_{t \in [\tau_o - \tau_1, \tau_o]} & \uprho (t) \big( \big(C_{r_o/2}^j (y_o, \tau_o) \big)^{\varsigma_0} (t) \big) \tau_1 <
		\int_{\tau_o - \tau_1}^{\tau_o} \uprho (t) \left( A(t) \right) \,dt =																			\\
& = \int_{I}  \uprho (t) \left( A(t) \right) \,dt + \int_{(\tau_o - \tau_1, \tau_o) \setminus I}  \uprho (t) \left( A(t) \right) \,dt \leqslant									\\
& \leqslant \max_{t \in [\tau_o - \tau_1, \tau_o]} \uprho (t) \big( \big(C_{r_o/2}^j (y_o, \tau_o) \big)^{\varsigma_0} (t) \big) \, |I| + 									\\
&		\hskip30pt + \frac{\bar\nu}{2} \, \frac{\upkappa^2}{(1 + \upkappa) {\sf L}} \, \max_{t \in [\tau_o - \tau_1, \tau_o]} 
		\uprho (t) \big( \big(C_{r_o/2}^j (y_o, \tau_o) \big)^{\varsigma_0} (t) \big) \left| (\tau_o - \tau_1, \tau_o) \setminus I \right| \leqslant						\\
& \leqslant \max_{t \in [\tau_o - \tau_1, \tau_o]} \uprho (t) \big( \big(C_{r_o/2}^j (y_o, \tau_o) \big)^{\varsigma_0} (t) \big) \tau_1
		\left[  \frac{|I|}{\tau_1} + \frac{\bar\nu}{2} \, \frac{\upkappa^2}{(1 + \upkappa) {\sf L}} \right] \leqslant												\\
& \leqslant \frac{(1 + \upkappa) {\sf L}}{\upkappa^2} \min_{t \in [\tau_o - \tau_1, \tau_o]} \uprho (t) \big( \big(C_{r_o/2}^j (y_o, \tau_o) \big)^{\varsigma_0} (t) \big) \tau_1
		\left[  \frac{|I|}{\tau_1} + \frac{\bar\nu}{2} \, \frac{\upkappa^2}{(1 + \upkappa) {\sf L}} \right]
\end{align*}
by which we derive the following lower bound on $I$:
\begin{align*}
|I| > \frac{\bar\nu}{2} \, \frac{\upkappa^2}{(1 + \upkappa) {\sf L}} \, \frac{r_o^2}{4} \, {\sf h}^{\ast}_j (x_o, 0, r) .
\end{align*}
On the other hand by  \textsl{Step 2} we get
\begin{align*}
& \frac{\alpha \, (2^{\xi} M)^2}{r_o^2} \, \min_{t \in [ \tau_o - r_o^2 \, {\sf h}^{\ast}_j (x_o, 0, r), \tau_o ]} \big| \big(C_{r_o}^j (y_o, \tau_o) \big)^{\varsigma_1} (t) \big|
		\big( r_o^2 \, {\sf h}^{\ast}_j (x_o, 0, r) - |J_{\alpha}| \big) \leqslant \iint_{\big(C_{r_o}^j (y_o, \tau_o) \big)^{\varsigma_1}} |Du|^2 \, dx dt \leqslant 							\\
& \hskip20pt  \leqslant \frac{\gamma}{2} \, \frac{(1 + \upkappa) {\sf L}}{\upkappa^2} \, (2^{\xi} M)^2 \, 
		\min_{t \in [ \tau_o - 4 r_o^2 \, {\sf h}^{\ast}_j (x_o, 0, r), \tau_o ]} \big| \big( C_{2 r_o}^j (y_o, \tau_o) \big)^{\varsigma_2} (t) \big| \cdot					\\
& \hskip40pt \cdot \bigg[ ( 4 + r_o^2) \, {\sf h}^{\ast}_j (x_o, 0, r) + \frac{4}{3} \, \frac{\uprho ( ( C_{2 r_o}^j (y_o, \tau_o) )^{\varsigma_2} )}
	{r_o^2 \, {\sf h}^{\ast}_j (x_o, 0, r) \max_{t \in [ \tau_o - 4 r_o^2 \, {\sf h}^{\ast}_j (x_o, 0, r), \tau_o ]} | ( C_{2 r_o}^j (y_o, \tau_o) )^{\varsigma_2} (t) |} \bigg] .
\end{align*}
Notice that, using (H.2) - $ii \, )$, we get that for each $t$
\begin{align*}
\big| \big( C_{2 r_o}^j (y_o, \tau_o) \big)^{\varsigma_2} (t) \big| \leqslant | B_{2 r_o} (y_o) | \leqslant 2^n  | B_{r_o} (y_o) | \leqslant
	\frac{2^n}{2 \upkappa} \upchi_j \big( B_{r_o} (y_o) \big) \leqslant
	\frac{2^n}{2 \upkappa} \big| \big( C_{r_o}^j (y_o, \tau_o) \big)^{\varsigma_1} (t) \big|.
\end{align*}
Then 
\begin{align*}
\min_{t \in [ \tau_o - r_o^2 \, {\sf h}^{\ast}_j (x_o, 0, r), \tau_o ]} \big| \big(C_{r_o}^j (y_o, \tau_o) \big)^{\varsigma_1} (t) \big| \geqslant
	\frac{2 \upkappa}{2^n}
	\min_{t \in [ \tau_o - 4 r_o^2 \, {\sf h}^{\ast}_j (x_o, 0, r), \tau_o ]} \big| \big( C_{2 r_o}^j (y_o, \tau_o) \big)^{\varsigma_2} (t) \big|	.
\end{align*}
Finally we derive
\begin{align*}
| J_{\alpha} | & \geqslant  r_o^2 \, {\sf h}^{\ast}_j (x_o, 0, r) - r_o^2 \, \frac{2^n}{2 \upkappa} \, \frac{\gamma}{2 \alpha} \frac{(1 + \upkappa) {\sf L}}{\upkappa^2} \, \cdot 		\\
& \hskip30pt \cdot \bigg[ ( 4 + r_o^2) \, {\sf h}^{\ast}_j (x_o, 0, r) + \frac{4}{3} \, \frac{\uprho ( ( C_{2 r_o}^j (y_o, \tau_o) )^{\varsigma_2} )}
	{r_o^2 \, {\sf h}^{\ast}_j (x_o, 0, r) \max_{t \in [ \tau_o - 4 r_o^2 \, {\sf h}^{\ast}_j (x_o, 0, r), \tau_o ]} | ( C_{2 r_o}^j (y_o, \tau_o) )^{\varsigma_2} (t) |} \bigg] .
\end{align*}
Then, since
\begin{align*}
| I \cap J_\alpha | & = | I | + | J_\alpha | - | I \cup J_\alpha | \geqslant 																					\\
	& \geqslant \frac{\bar\nu}{2} \, \frac{\upkappa^2}{(1 + \upkappa) {\sf L}} \, \frac{r_o^2}{4} \, {\sf h}^{\ast}_j (x_o, 0, r) + r_o^2 \, {\sf h}^{\ast}_j (x_o, 0, r)
			- r_o^2 \, \frac{2^n}{2 \upkappa} \frac{\gamma}{2 \alpha} \frac{(1 + \upkappa) {\sf L}}{\upkappa^2} \, \cdot											\\
	& \hskip30pt \cdot \bigg[ ( 4 + r_o^2) \, {\sf h}^{\ast}_j (x_o, 0, r) + \frac{4}{3} \, \frac{\uprho ( ( C_{2 r_o}^j (y_o, \tau_o) )^{\varsigma_2} )}
	{r_o^2 \, {\sf h}^{\ast}_j (x_o, 0, r) \max_{t \in [ \tau_o - 4 r_o^2 \, {\sf h}^{\ast}_ j (x_o, 0, r), \tau_o ]} | ( C_{2 r_o}^j (y_o, \tau_o) )^{\varsigma_2} (t) |} \bigg] +				\\
	& \hskip40pt - r_o^2 \, {\sf h}^{\ast}_j (x_o, 0, r) =																						\\
	& = \frac{\bar\nu}{2} \, \frac{\upkappa^2}{(1 + \upkappa) {\sf L}} \, \frac{r_o^2}{4} \, {\sf h}^{\ast}_j (x_o, 0, r)
			- r_o^2 \, \frac{2^{n-2} \gamma}{\alpha} \frac{(1 + \upkappa) {\sf L}}{\upkappa^3} \, \cdot																	\\
	& \hskip30pt \cdot \bigg[ ( 4 + r_o^2) \, {\sf h}^{\ast}_j (x_o, 0, r) + \frac{4}{3} \, \frac{\uprho ( ( C_{2 r_o}^j (y_o, \tau_o) )^{\varsigma_2} )}
	{r_o^2 \, {\sf h}^{\ast}_j (x_o, 0, r) \max_{t \in [ \tau_o - 4 r_o^2 \, {\sf h}^{\ast}_j (x_o, 0, r), \tau_o ]} | ( C_{2 r_o}^j (y_o, \tau_o) )^{\varsigma_2} (t) |} \bigg] ,
\end{align*}
taking (for instance)
\begin{align*}
\frac{1}{\alpha} = & \, \frac{\bar\nu}{4} \, \frac{\upkappa^5}{(1 + \upkappa)^2 {\sf L}^2} \, {\sf h}^{\ast}_j (x_o, 0, r) \frac{1}{2^{n} \gamma} \cdot			\\
& \cdot \bigg[ ( 4 + r_o^2) \, {\sf h}^{\ast}_j (x_o, 0, r) + \frac{4}{3} \, \frac{\uprho ( ( C_{2 r_o}^j (y_o, \tau_o) )^{\varsigma_2} )}
	{r_o^2 \, {\sf h}^{\ast}_j (x_o, 0, r) \max_{t \in [ \tau_o - 4 r_o^2 \, {\sf h}^{\ast}_j (x_o, 0, r), \tau_o ]} | ( C_{2 r_o}^j (y_o, \tau_o) )^{\varsigma_2} (t) |} \bigg]^{-1}
\end{align*}
one gets that $I \cap J_{\alpha}$ has positive measure and in particular
$$
| I \cap J_\alpha | = | I | + | J_\alpha | - | I \cup J_\alpha | \geqslant \frac{\bar\nu}{4} \, \frac{\upkappa^2}{(1 + \upkappa) {\sf L}} \, \frac{r_o^2}{4} \, {\sf h}^{\ast}_j (x_o, 0, r) .
$$
Then we get the existence of
$\bar{t} \in \big[\tau_o - {\textstyle \frac{r_o^2}{4}} \, {\sf h}^{\ast}_j (x_o, 0, r) , \tau_o \big]$ such that
\eqref{carlettomio} hold. \\ [0.3em]
\textsl{\underline{Step 4} - }
The goal of this step is to show that for every $\delta \in (0,1)$ there are $\eta \in (0,1)$ and $\hat{x} \in B_{r_o / 2} (y_o)$,
such that $B_{\eta \frac{r_o}{2}} (\hat{x}) \subset \big(C_{r_o/2}^j (y_o, \tau_o) \big)^{\varsigma_0} (\bar{t})$ 
and such that
\begin{equation}
\label{estMis}
\uprho (\bar{t}) \Big(\Big\{u(\cdot,\bar t) \leqslant \frac{M}{4} \Big\} \cap B_{\eta \frac{r_o}{2}} (\hat{x}) \Big) \leqslant
	\delta \, \uprho (\bar{t}) \big( B_{\eta \frac{r_o}{2}} (\hat{x}) \big).
\end{equation}
To show this it is sufficient to use \textsl{Step 3} and to apply Lemma \ref{lemmaMisVar} to the function $u (\cdot, \bar{t})$ 
with $(x_o, t_o) = (y_o, \tau_o)$, $\tilde{t} = \bar{t}$, $\sigma_o = \varsigma_o$, $\sigma_1 = \varsigma_1$, $A = 4 \, \alpha \, (2^{\xi} M)^2$,
$B = \frac{\bar\nu}{2} \, \frac{\upkappa^2}{(1 + \upkappa) {\sf L}}$,
$a = M/2$, $\varepsilon = 1/2$, $r = r_o/2$, $\beta = {\sf h}^{\ast}_j (x_o, 0, r)$,
and $\alpha$ is the value defined in \eqref{alfa}. Then we get \eqref{estMis}. \\ [0.3em]
\textsl{\underline{Step 5} - }
In this step we want to show that an estimate like \eqref{estMis} holds also in a cylinder. \\
Precisely one can show that for every $\bar{\delta} \in (0,1)$ there is
$\varepsilon \in (0,1)$ which will depend only on $\bar\delta$, and
$\bar{s} = (\varepsilon \, \eta \, r_o/4)^2 \, {\sf h}^{\ast}_j (\hat{x}, \bar{t}, \frac{\eta r_o}{4})$ such that
($\hat{x}, \bar{t}, \eta, r_o$ as in the previous step)
\begin{align}
\label{cilindro_brutto}
\uprho \left( \left\{ u \leqslant \frac{M}{8}\right \} \cap \big( B_{\eta \frac{r_o}{4}}(\hat{x}) \times [\bar{t}, \bar{t} + \bar{s}] \big) \right)
	\leqslant \bar\delta \, \uprho \big( B_{\eta \frac{r_o}{4}} (\hat{x}) \times [\bar{t}, \bar{t} + \bar{s}] \big) \, .
\end{align}
The proof of this point is the same of Theorem 5.1 in \cite{fabio19}, and for this reason we avoid to report it. \\
Notice that, choosing $\bar{\delta}$ and consequently $\varepsilon$ sufficiently small, the cylinder $B_{\eta \frac{r_o}{4}}(\hat{x}) \times [\bar{t}, \bar{t} + \bar{s}]$
is contained in
$$
\big( \mathcal{Q}_j \big)^{\varsigma_0} := \mathcal{Q}_j \ \cup \ \big\{ (y,s) \in \mathcal{Q} \setminus \mathcal{Q}_j \, \big| \,
		\text{dist}_{\varsigma_0} \big( (y,s) , \mathcal{Q}_j \big) \leqslant 1 \big\} .
$$
where $\varsigma_0$ is defined in \eqref{sigmazero} and $\text{dist}_{\varsigma}$ is defined before \eqref{svaroschi}.
We will choose $\varepsilon$ in such a way that this last requirement will be satisfied, therefore in particular we have that
$$
\big( B_{\eta \frac{r_o}{4}}(\hat{x}) \times [\bar{t}, \bar{t} + \bar{s}] \big)_j := \big( B_{\eta \frac{r_o}{4}}(\hat{x}) \times [\bar{t}, \bar{t} + \bar{s}] \big) \cap {\mathcal Q}_j
= B_{\eta \frac{r_o}{4}}(\hat{x}) \times [\bar{t}, \bar{t} + \bar{s}] .
$$
\textsl{\underline{Step 6} - }
Also this step can be proved following the analogous one in \cite{fabio19}. Using (C.10) and taking in Proposition \ref{anninabella}, point $iv \, )$,
$R = \tilde{s} = \eta r_o/4$, $ s= \eta r_o/8$, $\tilde\vartheta = \varepsilon^2$, $\hat\vartheta = \varepsilon^2/4$ one gets that
\begin{align*}
u (x, t) \geqslant \frac{M}{16} \qquad \text{ for a.e. } (x,t) \in \left( B_{\frac{\eta r_o}{8}} (\hat{x})
	\times \left[ \bar{t} + \frac{3 \varepsilon^2}{4} \left( \frac{\eta \, r_o}{4} \right)^2 {\sf h}^{\ast}_j \left( \hat{x}, \bar{t}, \frac{\eta \, r_o}{4} \right), \bar{t} + \bar{s} \right] \right) \cap \mathcal{Q}_j
\end{align*}
where $\bar{s} = (\varepsilon \, \eta \, r_o/4)^2 \, {\sf h}^{\ast}_j (\hat{x}, \bar{t}, \frac{\eta r_o}{4})$.
Notice that in real
$$
B_{\frac{\eta r_o}{8}} (\hat{x})
	\times \left[ \bar{t} + \frac{3 \varepsilon^2}{4} \left( \frac{\eta \, r_o}{4} \right)^2 {\sf h}^{\ast}_j \left(\hat{x}, \bar{t}, \frac{\eta r_o}{4} \right), 
		\bar{t} + \varepsilon^2 \left( \frac{\eta \, r_o}{4} \right)^2 {\sf h}^{\ast}_j \left( \hat{x}, \bar{t}, \frac{\eta r_o}{4} \right) \right] \subset \mathcal{Q}_j.
$$
Now, calling
$$
\hat{s} := \bar{t} + \bar{s} = \bar{t} + (\varepsilon \, \eta \, r_o/4)^2 \, {\sf h}^{\ast}_j \left( \hat{x}, \bar{t}, \frac{\eta r_o}{4} \right)
\qquad \textrm{and} \qquad \hat{r} = \frac{\eta r_o}{8}
$$
we can rewrite the estimate above as
\begin{align}
\label{passozero}
u (x, t) \geqslant \frac{M}{16} \qquad \text{ for a.e. } (x, t) \in B_{\hat{r}} (\hat{x})
	\times \left[ \hat{s} - \frac{1}{4} \left( \frac{\varepsilon \, \eta \, r_o}{4} \right)^2 {\sf h}^{\ast}_j \left( \hat{x}, \bar{t}, \frac{\eta \, r_o}{4} \right), \hat{s} \right] \subset \mathcal{Q}_j .
\end{align}
%
%
\textsl{\underline{Step 7} - }
Since \eqref{passozero} holds in particular in $\hat{s}$, by Theorem \ref{esp_positivita} there is $\vartheta_{j} \in (0, 1)$
such that for every $\hat{\vartheta} \in (0, \vartheta_{j} / 4)$ there is $\lambda_{j} \in (0,1)$ such that
\begin{align*}
& u \geqslant \lambda_{j} \frac{M}{16} \qquad \text{a.e. in } \Big( B_{2 \hat{r}} (\hat{x}) \times
	\big[ \hat{s} + \hat\vartheta (4 \hat{r})^2 {\sf h}^{\ast}_j (\hat{x}, \hat{s}, 5 \hat{r}) , \hat{s} + \vartheta_{j} (4 \hat{r} )^2 {\sf h}^{\ast}_j (\hat{x}, \hat{s}, 5 \hat{r}) \big]
						\Big) \cap \mathcal{Q}_j ,													\\
& \text{where we need } \quad (4 \hat{r} )^2 {\sf h}^{\ast}_j (\hat{x}, \hat{s}, 5 \hat{r}) < \updelta .
\end{align*}
Notice that by (C.3) one has that
\begin{align*}
& \uprho_j (t) \big( B_{\hat{r}} (\hat{x}) \big) \geqslant \upkappa \, \uprho (\hat{s})  \big( B_{\hat{r}} (\hat{x}) \big) \qquad \text{for } t \in [ \hat{s} , \hat{s} + \updelta ] ,
\end{align*}
i.e., in particular,
$$
B_{\hat{r}}^j (\hat{x}; t) \not= \emptyset \quad \text{ for each} \quad t \in 
	\big[ \hat{s} + \hat\vartheta (4 \hat{r})^2 {\sf h}^{\ast}_j (\hat{x}, \hat{s}, 5 \hat{r}) , \hat{s} + \vartheta_{j} (4 \hat{r} )^2 {\sf h}^{\ast}_j (\hat{x}, \hat{s}, 5 \hat{r}) \big] .
$$
%
%
%
%
In this way for every $t \in \big[ \hat{s} + \hat\vartheta (4 \hat{r})^2 {\sf h}^{\ast}_j (\hat{x}, \hat{s}, 5 \hat{r}) , \hat{s} + \vartheta_j (4 \hat{r} )^2 {\sf h}^{\ast}_j (\hat{x}, \hat{s}, 5 \hat{r}) \big]$
and $y \in B_{\hat{r}}^j (\hat{x}; t)$, by (C.3),
$$
B_{\hat{r}}^j (y; t) \subset B_{2 \hat{r}}^j (\hat{x}; t) 
$$
and then
\begin{align*}
\uprho_j (\hat{s}_1) \big( B_{2 \hat{r}} (\hat{x}) \big) \geqslant \uprho_j (\hat{s}_1) \big( B_{\hat{r}} (y) \big) 
		\geqslant \upkappa \, \uprho (\hat{s}_1)  \big( B_{\hat{r}} (y) \big)
\end{align*}
and in particular
\begin{align}
\label{brodino}
B_{2 \hat{r}}^j (\hat{x}; t) \not= \emptyset 
\quad \text{ for every} \quad t \in 
	\big[ \hat{s} + \hat\vartheta (4 \hat{r})^2 {\sf h}^{\ast}_j (\hat{x}, \hat{s}, 5 \hat{r}) , \hat{s} + \vartheta_j (4 \hat{r} )^2 {\sf h}^{\ast}_j (\hat{x}, \hat{s}, 5 \hat{r}) \big] .
\end{align}
Now we fix $\hat{\vartheta} \in (0, \vartheta_{j} / 4)$, for instance $\hat{\vartheta} = \vartheta_{j} / 8$, and define
$$
\hat{s}_1 := \hat{s} = \hat\vartheta (4 \hat{r})^2 {\sf h}^{\ast}_j (\hat{x}, \hat{s}, 5 \hat{r}) , \qquad
\hat{t}_1 := \hat{s} + \vartheta_{j} (4 \hat{r})^2 {\sf h}^{\ast}_j (\hat{x}, \hat{s}, 5 \hat{r}) 
$$
so that
\begin{align*}
u \geqslant \lambda_{j} \frac{M}{16} \qquad \text{a.e. in } \Big( B_{2 \hat{r}} (\hat{x}) \times [ \hat{s}_1, \hat{t}_1 ] \Big) \cap \mathcal{Q}_j .
\end{align*}
Thanks to \eqref{brodino} and
applying again Theorem \ref{esp_positivita} in $B_{2 \hat{r}}^j (2 \hat{x}; t)$ for some $t \in [ \hat{s}_1, \hat{t}_1 ]$ we get that
\begin{align*}
u \geqslant \lambda_{j}^2 \frac{M}{16} \qquad \text{a.e. in } \Big( B_{4 \hat{r}} (\hat{x}) \times [ \hat{s}_2, \hat{t}_2 ] \Big) \cap \mathcal{Q}_j .
\end{align*}
where
\begin{align*}
\hat{s}_2 := \hat{s}_1 + \hat\vartheta \, (8 \hat{r})^2 {\sf h}^{\ast}_j (\hat{x}, \hat{s}_1, 10 \hat{r}) , \qquad
\hat{t}_2 := \hat{t}_1 + \vartheta_{j} \, (8 \hat{r})^2 {\sf h}^{\ast}_j (\hat{x}, \hat{s}_1, 10 \hat{r}) .				
\end{align*}
As before, 
\begin{align*}
B_{4 \hat{r}}^j (\hat{x}; t) \not= \emptyset		\qquad \text{for every } t \in [ \hat{s}_2, \hat{t}_2 ] .										
\end{align*}
After iterating this argument $m$ times one reaches
\begin{align*}
& u \geqslant \frac{M}{16} \, \lambda_j^m \qquad \text{a.e. in }	\quad \Big( B_{2^m \hat{r}} (\hat{x}) \times \big[ \hat{s}_m , \hat{t}_m \big] \Big) \cap \mathcal{Q}_j .
\end{align*}
where
\begin{align*}
& \hat{s}_0 := \hat{s} ,			\qquad	\hat\vartheta = \vartheta_{j}/8	,								 								\\
& \hat{s}_k := \hat{s}_{k - 1} + \hat\vartheta \, (2^{k+1} \hat{r})^2 {\sf h}^{\ast}_j (\hat{x}, \hat{s}_{k-1}, 5 \, 2^{k-1} \hat{r}) ,	\qquad k \in \{ 1, 2, \ldots  m \}			\\
& \hat{t}_k := \hat{t}_{k - 1} + \vartheta_{j} \, (2^{k+1} \hat{r})^2 {\sf h}^{\ast}_j (\hat{x}, \hat{s}_{k-1}, 5 \, 2^{k-1} \hat{r})	,	\qquad k \in \{ 1, 2, \ldots  m \}	
\end{align*}
Here we recall that we need
\begin{align}
\label{inaddition}
(2^{k+1} \hat{r} )^2 {\sf h}^{\ast}_j (\hat{x}, \hat{s}_{k-1}, 5 \cdot 2^{k-1} \hat{r}) < \updelta \qquad \text{for } k =1 , \ldots m
\end{align}
which follows by assumption. We recall that in addition to \eqref{inaddition} we need
$$
2^m \hat{r} \leqslant \bar{R} ,
$$
but this is ensured by the fact that we can choose $m$ in such a way that
\begin{align}
\label{carlettocanta}
2 r \leqslant 2^m \hat{r} \leqslant 4 r;
\end{align}
the estimate $2 r \leqslant 2^m \hat{r}$ ensures that $B_r (x_o) \subset B_{2^m \hat{r}} (\hat{x}) \subset B_{5r} (x_o)$.
Now, according to the choice of $r_o$ made in \eqref{choicey0t0}, we have that
$$
\hat{r} = \frac{\eta r_o}{8} = \frac{\eta (r - s_o)}{16}
$$
and recalling that $M = b (r - s_o)^{-\xi}$ we have, using now the second inequality in \eqref{carlettocanta}, i.e. $2^m \hat{r} \leqslant 4r$,
\begin{align*}
& u \geqslant \frac{M}{16} \lambda_{j}^m = \frac{b (r - s_o)^{-\xi}}{16} \lambda_{j}^m =
\frac{b}{16} \left( \frac{\eta}{16 \, \hat{r}} \right)^{\xi} \lambda_{j}^m
	\geqslant 2^{\xi m} \lambda_{j}^m \, b \, \eta^{\xi} \, 2^{-6\xi - 4} r^{-\xi}
\end{align*}
a.e. in $B_{2^m \hat{r}} (\hat{x}) \times \big[ \hat{s}_m , \hat{t}_m \big]$.
First of all we get rid of the dependence of $m$, still to be chosen, choosing $\xi$ in such a way that
$$
2^{\xi} \lambda_{j} = 1,
$$
so that in particular we have, independently of $m$,
\begin{align}
\label{anacaterina}
u \geqslant b \, \eta^{\xi} \, 2^{-6\xi - 4} r^{-\xi} = c_o \, u (x_o, 0) \qquad  \text{a.e. in } B_{r} (x_o) \times \big[ \hat{s}_m , \hat{t}_m \big]
\end{align}
where $c_o := \eta^{\xi} \, 2^{-6\xi - 4}$ and $u (x_o, 0) = b \, r^{-\xi}$ for the choice made at the beginning of the proof. \\
Now to conclude we have to choose $\hat{s}_m$ and $\hat{t}_m$ in such a way that
\begin{align}
\label{mimancanoimieibimbi}
r^2 {\sf h}^{\ast}_j (x_o, 0, r) \in \big[ \hat{s}_m , \hat{t}_m \big]
\end{align}
where $c \in (0,1]$ is arbitrarily chosen. Since (here we use \eqref{carlettocanta} and the fact that $\hat{s} < 0$)
\begin{align*}
\hat{s}_m 	& = \hat{s} + \sum_{k = 1}^m \hat\vartheta \, (2^{k+1} \hat{r})^2 {\sf h}^{\ast}_j (\hat{x}, \hat{s}_{k-1}, 5 \, 2^{k-1} \hat{r}) \leqslant						\\
		& \leqslant	 \frac{16 \, (4^m - 1)}{3} \, \hat{r}^2 \max_{k = 1, \ldots m} \hat\vartheta \, {\sf h}^{\ast}_j (\hat{x}, \hat{s}_{k-1}, 5 \, 2^{k-1} \hat{r}) \leqslant		\\
		& \leqslant	 \frac{64}{3} \, r^2 \max_{k = 1, \ldots m} \hat\vartheta \max_{k = 1, \ldots m} {\sf h}^{\ast}_j (\hat{x}, \hat{s}_{k-1}, 5 \, 2^{k-1} \hat{r}) ,
\end{align*}
to get \eqref{mimancanoimieibimbi} first we require
$$
r^2 {\sf h}^{\ast}_j (x_o, 0, r) \geqslant \frac{64}{3} \, r^2 \max_{k = 1, \ldots m} \hat\vartheta \max_{k = 1, \ldots m} {\sf h}^{\ast}_j (\hat{x}, \hat{s}_{k-1}, 5 \, 2^{k-1} \hat{r}) 
$$
that is
$$
\hat\vartheta \leqslant \frac{3}{64} \, \frac{{\sf h}^{\ast}_j (x_o, 0, r)}{\max_{k = 1, \ldots m} {\sf h}^{\ast}_j (\hat{x}, \hat{s}_{k-1}, 5 \, 2^{k-1} \hat{r}) } .
$$
It remains to force $\hat{t}_m \geqslant r^2 {\sf h}^{\ast}_j (x_o, 0, r)$.
Clearly if this is true we are done choosing
$$
{\sf c}_j = c_o^{-1} \, .
$$
Otherwise, since $\hat{s}_m < r^2 {\sf h}^{\ast}_j (x_o, 0, r)$, by \eqref{anacaterina} there is $\tilde{t} > \hat{s}_m$, and $\tilde{t} > 0$, such that
$$
u (x, \tilde{t}) \geqslant c_o \, u (x_o, 0) \qquad  \text{a.e. in } B_{r} (x_o) .
$$
We can suppose, taking $\hat{\vartheta}$ smaller if necessary, that
\begin{align*}
\tilde{t} + \hat{\vartheta} \, r^2 {\sf h}^{\ast}_j (x_o, 0, r) \leqslant r^2 {\sf h}^{\ast}_j (x_o, 0, r) .
\end{align*}
We then get that,
by Theorem \ref{esp_positivita}, that for every $\hat\alpha \in (0, \vartheta_{j}/4)$ there is $\lambda_{j}' \in (0,1)$ such that
$$
\begin{array}{l}
u (x, t) \geqslant \lambda_{j}' \, c_o \, u (x_o, 0) 									\\		[0.7em]
\text{a.e. in } B_{2r} (x_o) \times
	\big[\tilde{t} + \hat\alpha \, (4 r)^2 {\sf h}^{\ast}_j (x_o, \tilde{t}, 5r) , \tilde{t} + \vartheta_{j} \, (4 r)^2 {\sf h}^{\ast}_j (x_o, \tilde{t}, 5r) \big]
\end{array}
$$
provided that
\begin{align*}
(4 r)^2 {\sf h}^{\ast}_j (x_o, \tilde{t}, 5r) < \updelta ,
\end{align*}
but this is true by assumption. For the same assumption notice that
$$
0 < \tilde{t} < r^2 {\sf h}^{\ast}_j (x_o, 0, r) \qquad \Longrightarrow \qquad \tilde{t} < \updelta
$$
and then by (C.5) we have that, since $\tilde{t} < \updelta$ and $5r \leqslant \bar{R}$,
$$
\left( \frac{\upkappa}{1 + \upkappa} \right)^2 {\sf h}^{\ast}_j (x_o, 0, 5r) \leqslant {\sf h}^{\ast}_j (x_o, \tilde{t}, 5r) \leqslant \left( \frac{1 + \upkappa}{\upkappa} \right)^2 {\sf h}^{\ast}_j (x_o, 0, 5r)
$$
then, taking a value for $\hat\alpha$ for which 
\begin{align}
\label{CaFi}
\hat\alpha < \frac{\vartheta_{j}}{4} \left( \frac{\upkappa}{1 + \upkappa} \right)^4 , \quad \text{for instance} \quad \hat\alpha = \frac{\vartheta_{j}}{8} \left( \frac{\upkappa}{1 + \upkappa} \right)^4 ,
\end{align}
we have (notice that we consider $B_r (x_o)$ and not $B_{2r} (x_o)$)
\begin{equation}
\label{anninapigrina}
\begin{array}{l}
u (x, t) \geqslant \lambda_{j}' \, c_o \, u (x_o, 0) 		\qquad \text{a.e. in } 						\\		[0.7em]
B_{r} (x_o) \times
	\big[ \tilde{t} + \hat\alpha \, \big( \frac{1 + \upkappa}{\upkappa} \big)^2 (4 r)^2 {\sf h}^{\ast}_j (x_o, 0, 5r) , 
		\tilde{t} + \vartheta_{j} \, \big( \frac{\upkappa}{1 + \upkappa} \big)^2 (4 r)^2 {\sf h}^{\ast}_j (x_o, 0, 5r) \big] .
\end{array}
\end{equation}
As before, if $\tilde{t} + \vartheta_{j} \, (4 r)^2 {\sf h}^{\ast}_j (x_o, 0, 5r) > r^2 {\sf h}^{\ast} (x_o, 0, r)$ we can conclude taking
$$
{\sf c}_j = (\lambda_{j}' \, c_o)^{-1} ,
$$
otherwise we go on. Applying Theorem \ref{esp_positivita} we get
\begin{align*}
& u (x, t) \geqslant ( \lambda_{j}' )^2 c_o \, u (x_o, 0) \qquad  \text{a.e. in } \\
& B_{r} (x_o) \times
	\Big[ \tilde{t} + \hat\alpha \, \Big( \frac{1 + \upkappa}{\upkappa} \Big)^2 (4 r)^2 {\sf h}^{\ast}_j (x_o, 0, 5r) 
	+ \hat\alpha (4 r)^2 {\sf h}^{\ast}_j \Big( x_o, \tilde{t} + \hat\alpha \, \Big( \frac{1 + \upkappa}{\upkappa} \Big)^2 (4 r)^2 {\sf h}^{\ast}_j (x_o, 0, 5r), 5r \Big) , 	\\
&	\qquad\qquad\qquad \tilde{t} + \vartheta_{j} \, \Big( \frac{\upkappa}{1 + \upkappa} \Big)^2 (4 r)^2 {\sf h}^{\ast}_j (x_o, 0, 5r) +								\\
&	\qquad\qquad\qquad +		\vartheta_{j} \, \Big( \frac{\upkappa}{1 + \upkappa} \Big)^2 (4 r)^2 
		{\sf h}^{\ast}_j \Big( x_o, \tilde{t} + \hat\alpha \, \Big( \frac{1 + \upkappa}{\upkappa} \Big)^2 (4 r)^2 {\sf h}^{\ast}_j (x_o, 0, 5r), 5r \Big) \Big] .
\end{align*}
Since, as above, by (C.5) we have
$$
\left( \frac{\upkappa}{1 + \upkappa} \right)^2 {\sf h}^{\ast}_j (x_o, 0, 5r) \leqslant
{\sf h}^{\ast}_j \Big( x_o, \tilde{t} + \hat\alpha \, \Big( \frac{1 + \upkappa}{\upkappa} \Big)^2 (4 r)^2 {\sf h}^{\ast}_j (x_o, 0, 5r), 5r \Big) \leqslant
\left( \frac{1 + \upkappa}{\upkappa} \right)^2 {\sf h}^{\ast}_j (x_o, 0, 5r) ,
$$
and since, by (C.9) and \eqref{valsusa},
\begin{align*}
& {\sf h}^{\ast}_j (x_o, 0, 5r) \leqslant 25 \, \frac{1 + \upkappa}{\upkappa} \, c_{\uprho_j} (5) \, {\sf h}^{\ast}_j (x_o, 0, r) ,			\\
& {\sf h}^{\ast}_j (x_o, 0, 5r) \geqslant \frac{2 \upkappa}{5^{n+2}} \, {\sf h}^{\ast}_j (x_o, 0, r) ,
\end{align*}
we derive
\begin{align*}
\left( \frac{\upkappa}{1 + \upkappa} \right)^2 \frac{2 \upkappa}{5^{n+2}} \, & {\sf h}^{\ast}_j (x_o, 0, r) \leqslant
{\sf h}^{\ast}_j \Big( x_o, \tilde{t} + \hat\alpha \, \Big( \frac{1 + \upkappa}{\upkappa} \Big)^2 (4 r)^2 {\sf h}^{\ast}_j (x_o, 0, 5r), 5r \Big) \leqslant		\\
&	\leqslant 25 \, \left( \frac{1 + \upkappa}{\upkappa} \right)^3 c_{\uprho_j} (5) \, {\sf h}^{\ast}_j (x_o, 0, r) .
\end{align*}
Then we get
\begin{align*}
& u (x, t) \geqslant (\lambda_{j}')^2 c_o \, u (x_o, 0) \qquad  \text{a.e. in } \quad  B_{r} (x_o) \times
	\Big[ \tilde{t} + \hat{t}^{\ast}_1 + \hat{t}^{\ast}_2 , \tilde{t} + t^{\ast}_1 + t^{\ast}_2 \Big]
\end{align*}
where
\begin{align*}
\hat{t}^{\ast}_1 = \hat{t}^{\ast}_2 :=	& \ \hat\alpha \, 25 \, \left( \frac{1 + \upkappa}{\upkappa} \right)^5 c_{\uprho_j} (5) \, (4 r)^2 {\sf h}^{\ast}_j (x_o, 0, r) ,			\\
t^{\ast}_1 = t^{\ast}_2 :=	& \ \vartheta_{j} \, \Big( \frac{\upkappa}{1 + \upkappa} \Big)^4 \frac{2 \upkappa}{5^{n+2}} (4 r)^2 {\sf h}^{\ast}_j (x_o, 0, r) .
\end{align*}
If we, in addition to \eqref{CaFi}, we require that
$$
\hat\alpha < \frac{1}{c_{\uprho_j} (5)} \left( \frac{\upkappa}{1 + \upkappa} \right)^5 \frac{2 \upkappa}{5^{n+4}} \, \vartheta_j ,
$$
for this choice of $\hat\alpha$ we also have
$$
\tilde{t} + \hat{t}^{\ast}_1 + \hat{t}^{\ast}_2 > \tilde{t} + t^{\ast}_1
$$
so that, by \eqref{anninapigrina},
$$
u (x, t) \geqslant (\lambda_{j}')^2 c_o \, u (x_o, 0) \qquad  \text{a.e. in } \quad  B_{r} (x_o) \times \Big[ \tilde{t} + \hat{t}^{\ast}_1  , \tilde{t} + t^{\ast}_1 + t^{\ast}_2 \Big]
$$
since $\lambda_{j}'< 1$. Now, if $\tilde{t} + t^{\ast}_1 + t^{\ast}_2 > r^2 {\sf h}^{\ast} (x_o, 0, r)$ we are done taking
$$
{\sf c}_j = ( (\lambda_{j}')^2 \, c_o)^{-1} ,
$$
otherwise we go on. Iterating the argument we derive the existence of $k \in \N$ such that
$$
u (x, t) \geqslant c_o (\lambda_{j}')^k \, u (x_o, 0) \qquad  \text{a.e. in } \quad  B_{r} (x_o) \times
	\Big[ \tilde{t} + \hat{t}^{\ast}_1  , \tilde{t} + \sum_{h = 1}^k t^{\ast}_h \Big]
$$
with
\begin{align}
\label{ipazia}
\tilde{t} + \sum_{h = 1}^k t^{\ast}_h > r^2 {\sf h}^{\ast}_j (x_o, 0, r) .
\end{align}
where
$$
t^{\ast}_h = \vartheta_{j} \, \Big( \frac{\upkappa}{1 + \upkappa} \Big)^4 \frac{2 \upkappa}{5^{n+2}} (4 r)^2 {\sf h}^{\ast}_j (x_o, 0, r) .
$$
Since $\tilde{t} > 0$, the condition \eqref{ipazia} is satisfied if
$$
\sum_{h = 1}^k t^{\ast}_h > r^2 {\sf h}^{\ast}_j (x_o, 0, r)  .
$$
Moreover we want $\tilde{t} + \sum_{h = 1}^k t^{\ast}_h$ to remain in the domain: since $c \leqslant 1$ we require that
$$
r^2 {\sf h}^{\ast}_j (x_o, 0, r) \leqslant \sum_{h = 1}^k t^{\ast}_h = 
	k \, \vartheta_{j} \, \Big( \frac{\upkappa}{1 + \upkappa} \Big)^4 \frac{2 \upkappa}{5^{n+2}} (4 r)^2 {\sf h}^{\ast}_j (x_o, 0, r) < 2 \, r^2 {\sf h}^{\ast}_j (x_o, 0, r) 
$$
i.e.
\begin{equation}
\label{bound}
1 \leqslant k \, \vartheta_{j} \, \Big( \frac{\upkappa}{1 + \upkappa} \Big)^4 \frac{32 \, \upkappa}{5^{n+2}} < 2 .
\end{equation}
Clearly, taking if necessary a suitable and smaller value for $\vartheta_{j}$, this is always possible and the number of steps $k$ does not depend on ${\sf h}^{\ast}_j (x_o, 0, r)$. \\
Now we conclude taking ${\sf c}_j = (c_o (\lambda_{j}')^k)^{-1}$.
Notice that increasing $k$ does increase ${\sf c}_j$, but the estimates \eqref{bound} put a bound on $k$, therefore in fact ${\sf c}_j$ is independent of $k$.
\finedimo

\section*{Part II - Application to mixed type equations}

\ \\ [-3em]

\section{A preliminary result}

In this section we state a result,
fundamental to obtain the Harnack inequalities for elliptic-parabolic equations and for parabolic equations of forward-backward type
obtained in the following sections (see Theorem \ref{harnack0+}).

\begin{lemma}
\label{lemmino}
Consider a sequence of measurable functions  $(v_n)_{n \in \N}$ defined in $A \subset \R^k$, $A$ measurable and bounded, $|A| > 0$,
and a measurable function $v$ defined in $A$ such that
$$
v_n \to v \qquad \text{a.e. in } A .
$$
Consider a sequence $(A_n)_n$ of measurable and bounded sets satisfying
$$
A_n \subseteq A , \quad A_{n} \subseteq A_{n+1}, \qquad \lim_{n \to +\infty}\big| A \setminus A_n \big| = 0 . 
$$
Then
$$
\limsup_{n \to +\infty} \inf_{A_n} v_n  \leqslant \inf_{A} v \qquad \text{and} \qquad \liminf_{n \to +\infty} \sup_{A_n} v_n  \geqslant \sup_{A} v .
$$
\end{lemma}
\boss
Clearly the lemma is true not only if the sequence $(A_n)_n$ is contained in $A$ but, even more so, if the sequence $(A_n)_n$ satisfies
$A_n \supseteq A$ and
$$
v_n \to v \qquad \text{a.e. in } \bigcup_n A_n .
$$
Indeed in this case, by the lemma, we get $\limsup_{n \to +\infty} \inf_{A} v_n  \leqslant \inf_{A} v$ and then, a fortiori,
$\limsup_{n \to +\infty} \inf_{A_n} v_n  \leqslant \inf_{A} v$. \\
\eoss
\noindent
\dimo
We show the lemma under the assumption $v \in L^{\infty} (A)$, being the general case very similar.
Moreover we confine to show only the first inequality, being the proof of the other very similar (and, in fact, equivalent).
We start assuming
$$
v_n \geqslant v
$$
and prove
\begin{align}
\label{lim}
\lim_{n \to +\infty} \inf_{A_n} v_n = \inf_A v .
\end{align}
Since $\inf_A v = \essinf_A v$ is finite we consider, just for the sake of simplicity,
\begin{align*}
& u_n := v_n - \inf_A v \quad \text{and} \quad u := v - \inf_A v 
\end{align*}
in such a way that
\begin{align*}
u_n \geqslant u, \qquad u_n \to u , \qquad u_n, u \geqslant 0  .
\end{align*}
We want to show that
$$
\lim_{n \to +\infty} \inf_{A_n} u_n = 0 ,
$$
that is for every $\delta > 0$ there is $N \in \N$ such that
$$
0 \leqslant \inf_{A_n} u_n < \delta \qquad \text{ for every } n \geqslant N .
$$
Fix $\delta > 0$. Since $\inf_A u = 0$ we have that
$$
\big| A_{\delta} := \big\{ x \in A \, \big| u(x) < \delta/2 \big\} \big| > 0 .
$$
Moreover, by assumption, $| A \setminus A_n| \to 0$ when $n \to +\infty$, then $|A_{\delta} \setminus A_n | \to 0$ when $n \to +\infty$.
Then we can find $N_{1,\delta} \in \N$ such that $A_{\delta} \cap A_{N_{1,\delta}} \not= \emptyset$; moreover we can find
$x_{\delta} \in A_{\delta} \cap A_{N_{1,\delta}} \subseteq A_{\delta} \cap A_n$ for every $n \geqslant N_{1,\delta}$ in such a way that
\begin{align*}
\lim_{n \to +\infty} u_n (x_{\delta} ) = u (x_{\delta} ) \qquad
\text{ and } \qquad 0 \leqslant u (x_{\delta} ) < \frac{\delta}{2} .
\end{align*}
By the first one we have that there is $N_{2,\delta} \in \N$, depending on $\delta$ and $x_{\delta}$, such that
$$
0 \leqslant u_n (x_{\delta} ) < u (x_{\delta} ) + \frac{\delta}{2} \qquad \text{ for every } n \geqslant N_{2,\delta} .
$$
Then we conclude that
$$
0 \leqslant \inf_{A_n} u_n \leqslant u_n (x_{\delta} ) < u (x_{\delta} ) + \frac{\delta}{2} < \delta 
\qquad \text{for every } n \geqslant N_{2,\delta}
$$
i.e. \eqref{lim}. Now we consider a generic sequence $(v_n)_n$ converging to $v$ almost everywhere in $A$. We define
$$
\bar{v}_n := \max \{ v_n, v \}
$$
in such a way that
$$
\bar{v}_n \to v \qquad \text{a.e. in } A , \qquad  \bar{v}_n \geqslant v .
$$
Then
$$
\inf_A v_n \leqslant \inf_A \bar{v}_n \leqslant \inf_{A_n} \bar{v}_n
$$
and, by what seen above,
$$
\displaylines{
\hfill
\limsup_{n \to + \infty} \inf_{A_n} v_n \leqslant \lim_{n \to + \infty} \inf_{A_n} \bar{v}_n = \inf_A v \, .
\hfill\llap{$\square$}}
$$
\fine

\section{A Harnack inequality for solutions of elliptic-parabolic equations}
\label{paragrafo7}

In this section we consider equations \eqref{equazionegenerale1} and \eqref{equazionegenerale2} where $\uprho$ is a function
$$\uprho \geqslant 0 . $$
These are equations of non-standard parabolic type, but in the region where $\uprho = 0$ we look at them as a family of elliptic equations
in the parameter $t$, and also the Harnack type inequality that their solutions satisfy is of elliptic type. For this reason we refer to them as elliptic-parabolic. \\
These equations have solutions (see \cite{fabio4}, \cite{fabio9.1}) and these solutions belong to the corresponding De Giorgi class,
defined by \eqref{lavagna_bis}, but with $\uprho \geqslant 0$. \\
We will suppose that there exist $\mathcal{Q}_1, \mathcal{Q}_2$ and $\Upgamma$ as in Section 2, but we will write
$$
\mathcal{Q}_+ := \mathcal{Q}_1 , \qquad \mathcal{Q}_0 := \mathcal{Q}_2 .
$$
For every $E \subset \Omega \times (0,T)$ and $A \subset \Omega$ we define
\begin{gather*}
E^+ := E \cap \mathcal{Q}_+ ,  \qquad E^0 := E \cap \mathcal{Q}_0 , 					\\
A^+ (t) := \big( A \times \{ t \} \big)\cap \mathcal{Q}_+ ,  \qquad
	A^0 (t) := \big( A \times \{ t \} \big)\cap \mathcal{Q}_0 .
\end{gather*}
Similarly we will denote by 
\begin{gather*}
B_r^+ (x_o; t) = \big( B_r ( x_o) \times \{ t \} \big) \cap \mathcal{Q}_+	,	\qquad
	B_r^0 (x_o; t) = \big( B_r ( x_o) \times \{ t \} \big) \cap \mathcal{Q}_0 .
\end{gather*}
We will suppose that there exists $\uprho_0$ such that
\begin{gather*}
\left\{
\begin{array}{cl}
\uprho > 0		&	\text{a.e. in } \mathcal{Q}_+	,	\\	[0.3em]
\uprho = 0		&	\text{a.e. in  } \mathcal{Q}_0 ,
\end{array}
\right.
\qquad
\left\{
\begin{array}{cl}
\uprho_0 = 0		&	\text{a.e. in } \mathcal{Q}_+	,	\\	[0.3em]
\uprho_0 > 0		&	\text{a.e. in  } \mathcal{Q}_0 ,
\end{array}
\right.
\ \\
\uprho	\quad \text{satisfies assumptions (H.1)} \qquad \text{and} \qquad \bar\uprho : = \uprho + \uprho_0 \quad \text{satisfies assumptions (H.1)-(H.4)} .
\end{gather*}
Now consider $\epsilon \in (0,1]$ and define
\begin{align*}
\uprho_{\epsilon} :=
\left\{
\begin{array}{cl}
\uprho			&	\text{ in } \mathcal{Q}_+	,	\\	[0.3em]
\epsilon \, \uprho_0 	&	\text{ in  } \mathcal{Q}_0 .
\end{array}
\right.
\end{align*}
Finally we define
$$
{\sf h}_{\ast}^+ (x_o, t_o, r) := {\sf h}_{\ast}^1 (x_o, t_o, r) , \qquad {\sf h}^{\ast}_+ (x_o, t_o, r) := {\sf h}^{\ast}_1 (x_o, t_o, r) ,
$$
$$
\begin{array}{c}
{\sf h}_{\ast}^{\epsilon} (x_o, t_o, r) := 
{\displaystyle \frac{ \epsilon \, \uprho_0 \big( \big( B_r (x_o) \times (t_o - r^2 , t_o) \big) \big)}
		{\big| \big( B_r (x_o) \times (t_o - r^2 , t_o) \big) \cap \mathcal{Q}_0 \big|}}		,						\\	[1.5em]
{\sf h}^{\ast}_{\epsilon} (x_o, t_o, r) := 
{\displaystyle \frac{\epsilon \, \uprho_0 \big( \big( B_r (x_o) \times (t_o, t_o + r^2) \big) \big)}
		{\big| \big( B_r (x_o) \times (t_o, t_o + r^2) \big) \cap \mathcal{Q}_0 \big| }}	 ,
\end{array}
$$
while ${\sf h}_{\ast}^{\epsilon} (x_o, t_o, r) = 0$ if $\big( B_r (x_o) \times (t_o - r^2 , t_o) \big) \cap \mathcal{Q}_0 = \emptyset$,
${\sf h}^{\ast}_{\epsilon} (x_o, t_o, r) = 0$ if $\big( B_r (x_o) \times (t_o, t_o + r^2) \big) \cap \mathcal{Q}_0 = \emptyset$. 
Finally for $r, R, \theta > 0$ with $r < R$
\begin{align*}
Q_{r, \theta}^+ (x_o, t_o) = \bigcup_{t \in [t_o - \theta \, r^2 \, {\sf h}_{\ast}^+ ( x_o, t_o, R ), t_o]} B_r^+ ( x_o ; t) ,		\\
Q_{r, \theta}^{\epsilon} (x_o, t_o) = \bigcup_{t \in [t_o - \theta \, r^2 \, {\sf h}_{\ast}^{\epsilon} ( x_o, t_o, R ), t_o]} B_r^0 ( x_o ; t) .
\end{align*}

\noindent
{\bf Assumptions - } Notice that $\uprho_{\epsilon}$ satisfies (H.2) and (H.3); since $\bar\uprho$ and $\uprho$ satisfy (H.1) also
$\uprho_{\epsilon}$ satisfies (H.1) for every $\epsilon > 0$ (this is in fact an assumption about $\Upgamma$ and for this we refer
to the examples in \cite{fabio4}, \cite{fabio9.1}). 
As regards (H.4), this assumption is required just for Theorem \ref{gutierrez-wheeden} to hold. Then assuming that $\bar\uprho$ satisfies (H.4)
is sufficient for Theorem \ref{gutierrez-wheeden} to hold with $\uprho_{\epsilon}$ in the place of $\upeta$. 
Indeed, if we confine to consider $\epsilon \in (0,1]$ (which is sufficient for our goal, since we will let $\epsilon$ go to zero) we have that
\begin{gather*}
\uprho \leqslant \uprho_{\epsilon} \leqslant \bar\uprho ,
\end{gather*}
and then taking in Theorem \ref{gutierrez-wheeden} (and fixing $j = 1$) $\uprho$ in the place of $\uprho_1$ (since we are considering $\mathcal{Q}_+ := \mathcal{Q}_1$),
$\uprho_{\epsilon}$ in the place of $\upeta$,  $\bar\uprho$ in the place of $\uprho$
we conclude that Theorem \ref{gutierrez-wheeden} holds in particular for $\uprho_{\epsilon}$, at least for $\epsilon \in (0,1]$.
As regards $j = 2$ we have that
\begin{gather*}
\uprho_0 \leqslant \frac{1}{\epsilon} \, \uprho_{\epsilon} \leqslant \frac{1}{\epsilon} \, \bar \uprho  \qquad \text{for every } \epsilon > 0 
\end{gather*}
Then if $\bar\uprho$ satisfies (H.4), also $\epsilon^{-1} \bar\uprho$ satisfies (H.4) and then
$\epsilon^{-1} \uprho_{\epsilon}$ can be taken in the place of $\upeta$ and $\epsilon^{-1} \bar\uprho$ in the place of $\uprho$ in Theorem \ref{gutierrez-wheeden}.
Since in the inequality of Theorem \ref{gutierrez-wheeden} we are taking the mean values, we conclude that Theorem \ref{gutierrez-wheeden} holds also
with $\uprho_{\epsilon}$ in the place of $\upeta$ and $\bar\uprho$ in the place of $\uprho$. \\ [0.4em]
Summing up
\begin{equation}
\label{cuscino}
\begin{array}{c}
\uprho \quad \text{satisfies (H.1)} \qquad \text{and} \qquad \bar\uprho \quad \text{satisfies (H.1), (H.2), (H.3), (H.4)} 			\\
\ \\	[-0.8em]
\qquad \Downarrow \qquad																			\\
\ \\	[-0.8em]
\uprho_{\epsilon} \quad \text{satisfies (H.1), (H.2), (H.3)} \quad \text{uniformly in } \epsilon								\\ [0.5em]
\text{and } \uprho_{\epsilon} \text{ can be taken in the place of } \upeta \text{ in Theorem } \ref{gutierrez-wheeden}
\end{array}
\end{equation}

\noindent
Consider
\begin{gather*}
\Sigma_p := \big( \partial \Omega \times (0,T) \big) \cup \big( \Omega  \times \{ 0 \} \big) , \qquad 
\Sigma := \big( \partial \Omega \times (0,T) \big) \cup \big( \Omega_+ (0)  \times \{ 0 \} \big) ,
\end{gather*}
where $\Omega_+ (0) = \{ x \in \Omega \, | \, \uprho (x,0) > 0 \}$.
Denote by $(u_{\epsilon})_{\epsilon \in (0,1]}$ a family of solution to equation
\begin{gather*}
\left\{
\begin{array}{l}
{\displaystyle \uprho_{\epsilon} (x,t) \frac{\partial u}{\partial t} } - \textrm{div} \big(A (x,t,u,Du) \big) + B (x,t,u,Du) + C(x,t,u) = 0 	\qquad \text{in } \Omega \times (0,T)	,	\\	[0.5em]
u = \phi \qquad \text{in } \Sigma_p
\end{array}
\right.
\end{gather*}
for some boundary datum $\phi \in L^{\infty} (\Sigma_p)$ and by $u$ the solution of equation 
\begin{gather*}
\left\{
\begin{array}{l}
{\displaystyle \uprho (x,t) \frac{\partial u}{\partial t} } - \textrm{div} \big(A (x,t,u,Du) \big) + B (x,t,u,Du) + C(x,t,u) = 0 	\qquad \text{in } \Omega \times (0,T)	,	\\	[0.5em]
u = \phi \qquad \text{in } \Sigma ,
\end{array}
\right.
\end{gather*}
but the result holds also with different boundary condition.
Then (see Theorem 3.8 in \cite{fabio6})
\begin{gather}
\label{convergenza0}
u_{\epsilon} \to u	\qquad \text{in } L^2 (0,T; H^1 (\Omega)) \, \cap \, L^2 (\Omega \times (0,T); \uprho)
\end{gather}
and in particular (for the sake of simplicity we suppose the family $(u_{\epsilon})_{\epsilon \in (0,1]}$, and not a selected sequence, is pointwise converging)
\begin{gather}
\label{convergenza}
u_{\epsilon} \to u	\qquad \text{a.e. in } \Omega \times (0,T) .
\end{gather}
The same happens if $(u_{\epsilon})_{\epsilon \in (0,1]}$ is a family of solution to equation
\begin{gather*}
{\displaystyle \frac{\partial}{\partial t} } \big( \uprho_{\epsilon} (x,t) u \big) - \textrm{div} \big( A (x,t,u,Du) \big) + B (x,t,u,Du) + C(x,t,u) = 0 \qquad \text{in } \Omega \times (0,T)	
\end{gather*}
with the same boundary datum $\phi$ in $\Sigma_p$ and by $u$ the solution of equation \eqref{equazionegenerale2} with the same boundary datum restricted to $\Sigma$. \\
Now we state the Harnack type inequality for solutions of elliptic-parabolic equations.

\bthm
\label{harnack0+}
For every $r > 0$ such that
$B_{5 r} (x_o) \times (t_o - (5r)^2 {\sf h}_{\ast}^+ (x_o, t_o, 5r), t_o + (5 r)^2 {\sf h}^{\ast}_+ (x_o, t_o, 5r)) \subset \Omega \times (0,T)$ and 
$B_{r} (x_o) \times (t_o , t_o + 2 \, r^2 {\sf h}^{\ast}_+ (x_o, t_o, r)) \subset \Omega \times (0,T)$,
$$
5r \leqslant \bar{R} , \qquad \sup \, \{ (4 \varrho)^2 {\sf h}^{\ast}_+ ( x, t, 5 \varrho) \, | \, 
	(x,t) \in Q \cap \mathcal{Q}_+ , \varrho \leqslant r\}  < \updelta ,
$$
where $Q := B_{4 r} (x_o) \times (t_o - r^2 {\sf h}_{\ast}^+ (x_o, t_o, r), t_o + r^2 {\sf h}^{\ast}_+ (x_o, t_o, r))$, 
there exists $\upeta$ such that
for every non-negative solution $u$ of \eqref{equazionegenerale1} or of \eqref{equazionegenerale2}
one has that
\begin{itemize}
\item[$i \, )$]
if $(x_o, t_o) \in \Upgamma$ then for almost every $t_o \in (0,T)$
\begin{align*}
& \max \{ \sup_{B_r^+ (x_o; t_o - c \, r^2 {\sf h}_{\ast}^+ (x_o, t_o, r))} u(x , t_o - c \, r^2 {\sf h}_{\ast}^+ (x_o, t_o, r)) ,
	\sup_{B_r^0 (x_o; t_o)} u(x , t_o) \}	\leqslant							\\
& \qquad  \leqslant \upeta \, \min\{ \inf_{B_r^+ (x_o; t_o + c \, r^2 {\sf h}^{\ast}_+ (x_o, t_o, r))} u(x , t_o + c \, r^2 {\sf h}^{\ast}_+ (x_o, t_o, r)) ,
	\inf_{B_r^0 (x_o; t_o)} u(x , t_o) \} ;
\end{align*}
\item[$ii \, )$]
if $(x_o, t_o) \in {\mathcal Q}_+$ and $P_{r,  {\sf h}^{\ast}_+ (x_o, t_o, r))} (x_o, t_o) \cap \Upgamma \not= \emptyset$ $($defined analogously as \eqref{parabola}$)$,
then there is $(x_1, t_1) \in \Upgamma$ with $t_1 > t_o$ and $| x_1 - x_o | = \bar\rho < r$, such that
\begin{align*}
u(x_o, t_o) \leqslant \upeta \, \min \big\{ \inf_{B_{r}^+ (x_o; t_o + r^2 {\sf h}^{\ast}_+ (x_o, t_o, r))} u (x , t_o + r^2 {\sf h}^{\ast}_+ (x_o, t_o, r)) ,  
		\inf_{B_{r - \bar\rho}^0 (x_1; t_1)} u(x , t_1) \big\} ;
\end{align*}
\item[$iii \, )$]
if $(x_o, t_o) \in {\mathcal Q}_0$ and $\big( B_r (x_o) \times \{ t_o \} \big) \cap \Upgamma \not= \emptyset$
then for almost every $t_o \in (0,T)$ there is $(x_1, t_o) \in \Upgamma$ with $| x_1 - x_o | = \bar\rho < r$, such that
\begin{align*}
u (x_o, t_o) \leqslant	\upeta \, \min \{ \inf_{B_{r - \bar\rho}^+ (x_1; t_o + (r - \bar\rho)^2 {\sf h}^{\ast}_+ (x_1, t_o, r))} u(x , t_o + (r - \bar\rho)^2 {\sf h}^{\ast}_+ (x_1, t_o, r)) ,
	\inf_{B_r^0 (x_o; t_o)} u(x , t_o) \} .
\end{align*}
\end{itemize}
The constant $\upeta$ depends $($only$)$ on 
$\upgamma, \gamma, \kappa, r, \vartheta_j, c_{\uprho_j}, c_{\upchi_j} , \upkappa, {\sf L}, {\sf L}^{-1}, \upbeta$.
\ethm
\boss
Notice that the inequality in point $iii)$ holds for {\em almost} every $t_o$ and not for {\em every} $t_o$. This holds in fact only for the inequality
\begin{align*}
\sup_{B_r^0 (x_o; t_o)} u(x , t_o) \leqslant \upeta \, \inf_{B_r^0 (x_o; t_o)} u(x , t_o) .
\end{align*}
\eoss
\noindent
\dimo
By what observed above we have a family of functions $(u_{\epsilon})_{\epsilon \in (0,1]}$ satisfying \eqref{convergenza}. \\ [0.3em]
$i \, )$ - First suppose $(x_o, t_o) \in \Upgamma$ and fix $r > 0$.
By Theorem \ref{harnack}, using first \eqref{cruccodim2} and then \eqref{cruccodim}, and since $(x_o, t_o) \in \Upgamma$ (see also Remark \ref{portugal}), we have that
\begin{align*}
& \max  \big\{ \sup_{B_r^+ (x_o; t_o + r^2 {\sf h}^{\ast}_+ (x_o, t_o, r))} u_{\epsilon} (x , t_o - r^2 {\sf h}_{\ast}^+ (x_o, t_o, r)) , 										\\
& 	\hskip100pt		\sup_{B_r^0 (x_o; t_o + r^2 {\sf h}^{\ast}_{\epsilon} (x_o, t_o, r))} u_{\epsilon} (x , t_o - r^2 {\sf h}_{\ast}^{\epsilon} (x_o, t_o, r)) \big\} \leqslant		\\
& \hskip30pt \leqslant \min \{ \hat\upeta_{+} , \hat\upeta_{\epsilon} \}  \, 
	\min \big\{ \inf_{B_r^+ (x_o; t_o + r^2 {\sf h}^{\ast}_+ (x_o, t_o, r))} u_{\epsilon} (x , t_o + r^2 {\sf h}^{\ast}_+ (x_o, t_o, r)) ,										\\
& 	\hskip150pt \inf_{B_r^0 (x_o; t_o + r^2 {\sf h}^{\ast}_{\epsilon} (x_o, t_o, r))} u_{\epsilon} (x , t_o + r^2 {\sf h}^{\ast}_{\epsilon} (x_o, t_o, r)) \big\} ,
\end{align*}
for some positive constants $\hat\upeta_+, \hat\upeta_{\epsilon}$.
Observe that 
\begin{gather*}
\hat\upeta_+ \quad \text{depends, in particular, on} \quad  \frac{{\sf h}^{\ast}_+}{{\sf h}_{\ast}^+ + {\sf h}^{\ast}_+}
		\quad \text{and} \quad \frac{{\sf h}_{\ast}^+}{{\sf h}_{\ast}^+ + {\sf h}^{\ast}_+}  ,							\\
\hat\upeta_ {\epsilon} \quad \text{depends, in particular, on} \quad  \frac{{\sf h}^{\ast}_{\epsilon}}{{\sf h}_{\ast}^{\epsilon} + {\sf h}^{\ast}_{\epsilon}}
		\quad \text{and} \quad \frac{{\sf h}_{\ast}^{\epsilon}}{{\sf h}_{\ast}^{\epsilon} + {\sf h}^{\ast}_{\epsilon}}  ,
\end{gather*}
but in fact
\begin{align}
\label{bottarga}
\frac{{\sf h}^{\ast}_{\epsilon}}{{\sf h}_{\ast}^{\epsilon} + {\sf h}^{\ast}_{\epsilon}} 
\quad \text{and} \quad \frac{{\sf h}_{\ast}^{\epsilon}}{{\sf h}_{\ast}^{\epsilon} + {\sf h}^{\ast}_{\epsilon}}
\qquad \text{do not depend on } \epsilon ,
\end{align}
therefore we can chose $\upeta$ for which
\begin{align}
\label{sardefritte}
\begin{array}{l}
\max  \big\{ \sup_{B_r^+ (x_o; t_o + r^2 {\sf h}^{\ast}_+ (x_o, t_o, r))} u_{\epsilon} (x , t_o - r^2 {\sf h}_{\ast}^+ (x_o, t_o, r)) , 										\\
 	\hskip100pt		\sup_{B_r^0 (x_o; t_o + r^2 {\sf h}^{\ast}_{\epsilon} (x_o, t_o, r))} u_{\epsilon} (x , t_o - r^2 {\sf h}_{\ast}^{\epsilon} (x_o, t_o, r)) \big\} \leqslant		\\
 \hskip30pt \leqslant \upeta  \, 
	\min \big\{ \inf_{B_r^+ (x_o; t_o + r^2 {\sf h}^{\ast}_+ (x_o, t_o, r))} u_{\epsilon} (x , t_o + r^2 {\sf h}^{\ast}_+ (x_o, t_o, r)) ,									\\
	\hskip130pt \inf_{B_r^0 (x_o; t_o + r^2 {\sf h}^{\ast}_{\epsilon} (x_o, t_o, r))} u_{\epsilon} (x , t_o + r^2 {\sf h}^{\ast}_{\epsilon} (x_o, t_o, r)) \big\} .
\end{array}
\end{align}
Now observe what follows: by \eqref{convergenza0} and Theorem 1 in \cite{simon} we get that for every $t_1, t_2 \in (0,T)$, $t_1 < t_2$,
\begin{align}
\label{freschino}
\lim_{h \to 0} \int_{t_1}^{t_2} \| u_{\epsilon} (\sigma + h) - u_{\epsilon} (\sigma) \|_{H^1 (\Omega)} \, d\sigma = 0
\qquad \text{uniformly in } \epsilon .
\end{align}
Then we derive that
\begin{gather*}
\lim_{\epsilon \to 0^+} \int_{t_1}^{t_2} \| u_{\epsilon} (t - r^2 {\sf h}_{\ast}^{\epsilon} (x_o, t_o, r) ) - u_{\epsilon} (t) \|_{H^1 (\Omega)} \, d t	= 0 ,		\\
\lim_{\epsilon \to 0^+} \int_{t_1}^{t_2} \| u_{\epsilon} (t + r^2 {\sf h}^{\ast}_{\epsilon} (x_o, t_o, r) ) - u_{\epsilon} (t) \|_{H^1 (\Omega)} \, d t	= 0 
\end{gather*}
and in particular that
\begin{gather*}
\lim_{\epsilon \to 0^+} \big( u_{\epsilon} (x , t - r^2 {\sf h}_{\ast}^{\epsilon} (x_o, t_o, r) ) - u_{\epsilon} (x, t) \big) = 0	\qquad \text{a.e. in } \Omega \times (0,T) 		\\
\lim_{\epsilon \to 0^+} \big( u_{\epsilon} (x , t + r^2 {\sf h}^{\ast}_{\epsilon} (x_o, t_o, r) ) - u_{\epsilon} (x, t) \big) = 0	\qquad \text{a.e. in } \Omega \times (0,T) .
\end{gather*}
Then for almost every $t_o \in (0,T)$ we get
\begin{gather*}
\lim_{\epsilon \to 0^+} \big( u_{\epsilon} (x , t_o - r^2 {\sf h}_{\ast}^{\epsilon} (x_o, t_o, r) ) - u_{\epsilon} (x, t_o) \big) = 0	\qquad \text{a.e. in } \Omega 		\\
\lim_{\epsilon \to 0^+} \big( u_{\epsilon} (x , t_o + r^2 {\sf h}^{\ast}_{\epsilon} (x_o, t_o, r) ) - u_{\epsilon} (x, t_o) \big) = 0	\qquad \text{a.e. in } \Omega
\end{gather*}
and since for almost every $t_o \in (0,T)$
\begin{gather*}
\lim_{\epsilon \to 0^+} u_{\epsilon} (x, t_o) = u (x, t_o)		\qquad \text{a.e. in } \Omega 
\end{gather*}
we finally get
\begin{gather}
\label{attenzioneacarletto}
\begin{array}{l}
u_{\epsilon} (x , t_o - r^2 {\sf h}_{\ast}^{\epsilon} (x_o, t_o, r) ) \to u (x, t_o) , 		\\	[0.5em]
u_{\epsilon} (x , t_o + r^2 {\sf h}^{\ast}_{\epsilon} (x_o, t_o, r) ) \to u (x, t_o) ,
\end{array}
\qquad \text{a.e. in } \Omega 	\quad \text{and} \quad \text{for a.e. } t_o \in (0,T)
\end{gather}
and in particular \eqref{attenzioneacarletto} holds in $A \subset \Omega$ such that
$$
A \supset B_r^0 (x_o; t_o + \tau)  \quad \text{for every } \tau \in [- r^2 {\sf h}^{\ast}_{\epsilon} (x_o, t_o, r), r^2 {\sf h}^{\ast}_{\epsilon} (x_o, t_o, r) ], \ \epsilon \in (0,1].
$$
Now using \eqref{attenzioneacarletto} in $A$ as above and Lemma \ref{lemmino} in \eqref{sardefritte} we get
\begin{align*}
& \liminf_{\epsilon \to 0^+} \inf_{B_r^+ (x_o; t_o + r^2 {\sf h}^{\ast}_+ (x_o, t_o, r))} u_{\epsilon} (x , t_o + r^2 {\sf h}^{\ast}_+ (x_o, t_o, r)) \leqslant						\\
& \hskip30pt \leqslant \limsup_{\epsilon \to 0^+} \inf_{B_r^+ (x_o; t_o + r^2 {\sf h}^{\ast}_+ (x_o, t_o, r))} u_{\epsilon} (x , t_o + r^2 {\sf h}^{\ast}_+ (x_o, t_o, r)) \leqslant 		\\
& \hskip30pt \leqslant \inf_{B_r^+ (x_o; t_o + r^2 {\sf h}^{\ast}_+ (x_o, t_o, r))} u (x , t_o + r^2 {\sf h}^{\ast}_+ (x_o, t_o, r)) ,										\\
& \liminf_{\epsilon \to 0^+} \inf_{B_r^0 (x_o; t_o + r^2 {\sf h}^{\ast}_{\epsilon} (x_o, t_o, r))} u_{\epsilon} (x , t_o + r^2 {\sf h}^{\ast}_{\epsilon} (x_o, t_o, r)) \leqslant			\\
& \hskip30pt \leqslant \limsup_{\epsilon \to 0^+} 
						\inf_{B_r^0 (x_o; t_o + r^2 {\sf h}^{\ast}_{\epsilon} (x_o, t_o, r))} u_{\epsilon} (x , t_o + r^2 {\sf h}^{\ast}_{\epsilon} (x_o, t_o, r))  \leqslant 		\\
& \hskip30pt \leqslant \inf_{B_r^0 (x_o; t_o)} u (x , t_o)
\end{align*}
Analogous estimates for the left hand side in \eqref{sardefritte} lead to the following inequality
\begin{align*}
& \max \{ \sup_{B_r^+ (x_o; t_o - r^2 {\sf h}_{\ast}^+ (x_o, t_o, r))} u(x , t_o - r^2 {\sf h}_{\ast}^+ (x_o, t_o, r)) ,
	\sup_{B_r^0 (x_o; t_o)} u(x , t_o) \}	\leqslant							\\
& \qquad  \leqslant \upeta \, \min\{ \inf_{B_r^+ (x_o; t_o + r^2 {\sf h}^{\ast}_+ (x_o, t_o, r))} u(x , t_o + r^2 {\sf h}^{\ast}_+ (x_o, t_o, r)) ,
	\inf_{B_r^0 (x_o; t_o)} u(x , t_o) \}
\end{align*}
which holds for {\em almost} every $t_o \in (0,T)$. ``{\em Almost}'' comes from \eqref{attenzioneacarletto} and is essentially due to the inequality
\begin{align*}
\sup_{B_r^0 (x_o; t_o)} u(x , t_o) \}	\leqslant \upeta \, \inf_{B_r^0 (x_o; t_o)} u(x , t_o)
\end{align*}
which cannot hold for {\em every} $t_o \in (0,T)$.  \\ [0.3em]
$ii \, )$ -  Now consider $(x_o, t_o) \in {\mathcal Q}_+$, fix $r > 0$ and suppose that (the set $P_{r,  {\sf h}^{\ast}_j (x_o, t_o, r))} (x_o, t_o)$ is defined in \eqref{parabola})
$$
P_{r,  {\sf h}^{\ast}_j (x_o, t_o, r))} (x_o, t_o) \cap {\mathcal Q}_0 \not= \emptyset .
$$
Then there is
\begin{gather*}
\bar\rho \in (0, r] \qquad \text{such that} \qquad
\overline{B_{\bar\rho} (x_o) \times \{ t_o + {\bar\rho}^2 {\sf h}^{\ast}_+ (x_o, t_o, r) \}} \cap \Upgamma \not= \emptyset 		\\
\text{and } \overline{B_{\rho} (x_o) \times \{ t_o + {\rho}^2 {\sf h}^{\ast}_+ (x_o, t_o, r) \}} \cap \Upgamma = \emptyset 		\quad \text{for every } \rho \in (0, \bar\rho).
\end{gather*}
By Theorem \ref{harnack}, Remark \ref{maledetta_cervicale} and Remark \ref{portugal}
there is $\upeta$ such that
$$
u(x_o, t_o) \leqslant \upeta \, \inf_{P_{r,  {\sf h}^{\ast}_+ (x_o, t_o, r)}} u
$$
and in particular
\begin{equation}
\label{nonnepossopiu`}
\begin{array}{c}
{\displaystyle
u(x_o, t_o) \leqslant \upeta \, \inf_{B_{\bar\rho}^+ (x_o; t_o + \bar\rho^2 {\sf h}^{\ast}_+ (x_o, t_o, r))} u (x , t_o + \bar\rho^2 {\sf h}^{\ast}_+ (x_o, t_o, r)) ,		}	\\	[0.5em]
{\displaystyle
u(x_o, t_o) \leqslant \upeta \, \inf_{B_{r}^+ (x_o; t_o + r^2 {\sf h}^{\ast}_+ (x_o, t_o, r))} u (x , t_o + r^2 {\sf h}^{\ast}_+ (x_o, t_o, r)).	}
\end{array}
\end{equation}
Notice that $B_{\bar\rho}^+ (x_o; t_o + \bar\rho^2 {\sf h}^{\ast}_+ (x_o, t_o, r)) = B_{\bar\rho} (x_o)$.
This is represented in Figure 8.a.
Now consider a point
$$
(x_1, t_1) \in \overline{B_{\bar\rho} (x_o) \times \{ t_o + {\bar\rho}^2 {\sf h}^{\ast}_+ (x_o, t_o, r) \}} \cap \Upgamma , \quad \text{where } t_1 = t_o + \bar\rho^2 {\sf h}^{\ast}_+ (x_o, t_o, r) .
$$
Now suppose, for the the sake of simplicity, that for the same constant $\upeta$ we have (using point $i \, )$)
\begin{align*}
u(x_1, t_1) \leqslant \upeta \, \inf_{B_{\rho}^0 (x_1; t_1)} u(x , t_1)  \qquad \text{ for every } \rho \in (0, r].
\end{align*}
If this inequality were not true with the constant $\upeta$, since by point $i \,)$ this is true for almost every $t_1$,
we could find $(x_2, t_2) \in \overline{B_{\bar\rho} (x_o) \times \{ t_o + {\bar\rho}^2 {\sf h}^{\ast}_+ (x_o, t_o, r) \}} \cap \Upgamma$ with $t_2 > t_1$
for which the above inequality is true. \\
Then, taking $\rho = r - \bar\rho$ and using \eqref{nonnepossopiu`}, we get
\begin{align*}
u(x_o, t_o) & \leqslant \upeta \, \min \big\{ \inf_{B_{r}^+ (x_o; t_o + r^2 {\sf h}^{\ast}_+ (x_o, t_o, r))} u (x , t_o + r^2 {\sf h}^{\ast}_+ (x_o, t_o, r)) ,  u(x_1 , t_1) \big\} \leqslant			\\
	& \leqslant \upeta^2 \, \min \big\{ \inf_{B_{r}^+ (x_o; t_o + r^2 {\sf h}^{\ast}_+ (x_o, t_o, r))} u (x , t_o + r^2 {\sf h}^{\ast}_+ (x_o, t_o, r)) ,  
		\inf_{B_{r - \bar\rho}^0 (x_1; t_1)} u(x , t_1) \big\} \, .
\end{align*}
\ \\
\ \\
In the following pictures $\bar{t} = t_o + r^2 {\sf h}^{\ast}_+ (x_o, t_o, r)$.
\ \\
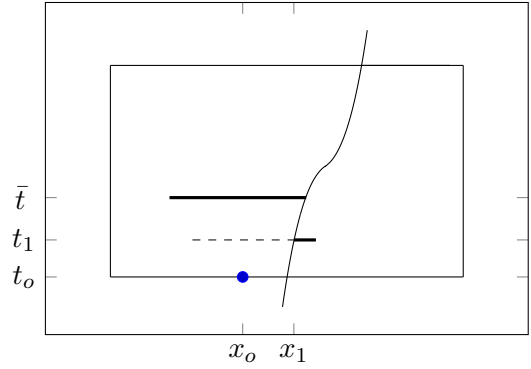
\begin{center}
\begin{tikzpicture}
\begin{axis} [axis equal, xtick={1.5}, ytick={-1.2, -0.78, -0.3}, xticklabels={$x_o$}, yticklabels={$t_o$, $t_1$, $\bar{t} \ $}, width=8cm, height=6cm, title={Figure 8.a}]
\addplot coordinates
{(1.5, -1.2)};
\addplot
[domain=0.5:2.5,variable=\t, smooth, ultra thin, dashed]
({t},{1.3*(t-1.5)*(t-1.5)-1.2});
\addplot
[domain=0:3,variable=\t, smooth]
({t},{-1.2});
\addplot
[domain=3:4,variable=\t, smooth]
({t},{-1.2});
\addplot
[domain=-1.2:1.2,variable=\t, smooth]
({0}, {t});
\addplot
[domain=-1.2:1.2,variable=\t, smooth]
({4}, {t});
\addplot
[domain=0:3.85,variable=\t, smooth]
({t},{1.2});
\addplot
[domain=2.85:4,variable=\t, smooth]
({t},{1.2});
\addplot
[domain=1.2:1.6,variable=\t,smooth]
({0.6*t^(1/3)+2-0.64 + 0.85}, {t});
\addplot
[domain=0.05:1.2,variable=\t, smooth]
({0.6*t^(1/3)+1.36 + 0.85}, {t});
\addplot
[domain=0.05:0.6,variable=\t, smooth]
({-0.6*t^(1/3)+1.36 + 0.45 + 0.85}, {-t + 0.11});
\addplot
[domain=0.6:1.65,variable=\t, smooth]
({-0.6*t^(1/3)+1.36 + 0.45 + 0.85}, {-t + 0.11});
\addplot
[domain=0.67:2.22,variable=\t, very thick]
({t},{-0.3});
\addplot
[domain=0.93:2.08,variable=\t, very thick]
({t},{-0.78});
\end{axis}
\end{tikzpicture}
\qquad \quad
\begin{tikzpicture}
\begin{axis} [axis equal, xtick={1.5, 2.08}, ytick={-1.2, -0.78, -0.3}, xticklabels={$x_o$, $x_1$, $\bar{t} \ $}, yticklabels={$t_o$, $t_1$, $\bar{t} \ $}, width=8cm, height=6cm, title={Figure 8.b}]
\addplot coordinates
{(1.5, -1.2)};
\addplot
[domain=0:3,variable=\t, smooth]
({t},{-1.2});
\addplot
[domain=3:4,variable=\t, smooth]
({t},{-1.2});
\addplot
[domain=-1.2:1.2,variable=\t, smooth]
({0}, {t});
\addplot
[domain=-1.2:1.2,variable=\t, smooth]
({4}, {t});
\addplot
[domain=0:3.85,variable=\t, smooth]
({t},{1.2});
\addplot
[domain=2.85:4,variable=\t, smooth]
({t},{1.2});
\addplot
[domain=1.2:1.6,variable=\t,smooth]
({0.6*t^(1/3)+2-0.64 + 0.85}, {t});
\addplot
[domain=0.05:1.2,variable=\t, smooth]
({0.6*t^(1/3)+1.36 + 0.85}, {t});
\addplot
[domain=0.05:0.6,variable=\t, smooth]
({-0.6*t^(1/3)+1.36 + 0.45 + 0.85}, {-t + 0.11});
\addplot
[domain=0.6:1.65,variable=\t, smooth]
({-0.6*t^(1/3)+1.36 + 0.45 + 0.85}, {-t + 0.11});
\addplot
[domain=0.93:2.08,variable=\t, ultra thin, dashed]
({t},{-0.78});
\addplot
[domain=2.08:2.33,variable=\t, very thick]
({t},{-0.78});
\addplot
[domain=0.67:2.22,variable=\t, very thick]
({t},{-0.3});
\end{axis}
\end{tikzpicture}
\end{center}
\ \\
\ \\
$iii \, )$ -  Now suppose $(x_o, t_o) \in {\mathcal Q}_0$ and, for $r > 0$, $\overline{B_r (x_o)} \times \{ t_o \} \cap \Upgamma \not= \emptyset$. Then there is
\begin{gather*}
\bar\rho \in (0, r] \qquad \text{such that} \qquad \Big( \overline{B_{\bar\rho} (x_o)} \times \{ t_o \} \Big) \cap \Upgamma \not= \emptyset 		\\
\text{and } \Big( \overline{B_{\rho} (x_o)} \times \{ t_o \} \Big) \cap \Upgamma = \emptyset 		\quad \text{for every } \rho \in (0,\bar\rho).
\end{gather*}
Then, arguing as in point $i \, )$, one can show that there is $\upeta$ such that
$$
\sup_{x \in B_r^0 (x_o; t_o)} u(x) \leqslant \upeta \, \inf_{x \in B_r^0 (x_o; t_o)} u(x) .
$$
Now consider $x_1 \in \overline{B_{\bar\rho} (x_o)}$ such that $(x_1, t_o) \in \big( \overline{B_{\bar\rho} (x_o)} \times \{ t_o \} \big) \cap \Upgamma$. Then
\begin{align*}
u (x_1, t_o) \leqslant \upeta \, \inf_{B_{\rho}^+ (x_1; t_o + {\rho}^2 {\sf h}^{\ast}_+ (x_1, t_o, r))} u(x , t_o + {\rho}^2 {\sf h}^{\ast}_+ (x_1, t_o, r)) \qquad \text{for every } \rho \in (0, r],
\end{align*}
where we considered $\upeta$ just for the the sake of simplicity (if the constant were different one can take the minimum between this constant and $\upeta$). \\
Taking $\rho = r - \bar\rho$ we finally get
\begin{align*}
u (x_o, t_o) \leqslant	\upeta^2 \, \min \{ \inf_{B_{r - \bar\rho}^+ (x_1; t_o + (r - \bar\rho)^2 {\sf h}^{\ast}_+ (x_1, t_o, r))} u(x , t_o + (r - \bar\rho)^2 {\sf h}^{\ast}_+ (x_1, t_o, r)) ,
	\inf_{B_r^0 (x_o; t_o)} u(x , t_o) \} .
\end{align*}
This concludes the proof.
\finedimo
\ \\ [-0.5em]
\subsection{H\"older continuity}
\label{sottoparagrafo7.1}
Using the classical argument due to Moser (see \cite{moser61}, but also, e.g., Section 7.9 in \cite{giusti}), one can prove that the solution of
(notice that $C \equiv 0$)
\begin{gather*}
\left\{
\begin{array}{l}
{\displaystyle \uprho (x,t) \frac{\partial u}{\partial t} } - \textrm{div} \big(A (x,t,u,Du) \big) + B (x,t,u,Du) = 0 	\qquad \text{in } \Omega \times (0,T)	,	\\	[0.5em]
u = \phi \qquad \text{in } \big( \partial \Omega \times (0,T) \big) \cup \big( \Omega_+(0) \times \{ 0 \} \big) .
\end{array}
\right.
\end{gather*}
is locally H\"older continuous in $\mathcal{Q}_+ \cap \Upgamma$ and $\Omega \ni x \mapsto u (x,t)$
is locally H\"older continuous for almost every $t$ in $\Omega_0 (t)$. \\
Notice that local H\"older continuity holds in particular in $\Upgamma$, while inside $\mathcal{Q}_+$ it follows, for instance, by the result in \cite{fabio19}. \\
To get local continuity in $\mathcal{Q}_0$ (with respect to $t$) cannot be expected in general.


\ \\


\end{document}